# LA CORRISPONDENZA TRA CELESTE CLERICETTI E LUIGI CREMONA DAL 1871 AL 1887 NELLE LETTERE CONSERVATE PRESSO L'ISTITUTO MAZZINIANO DI GENOVA

a cura di Paola Testi Saltini

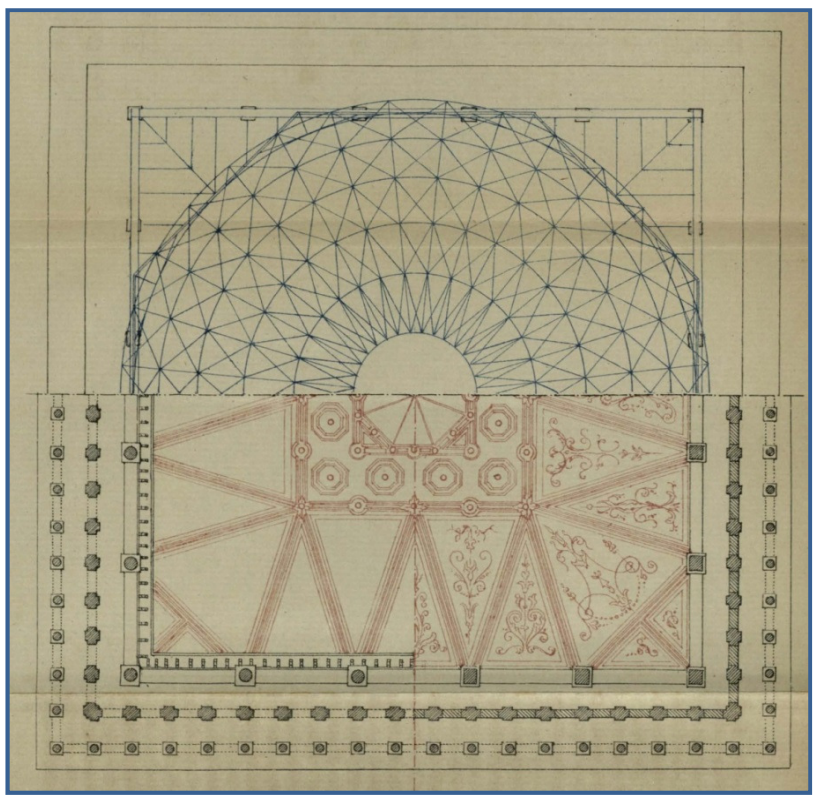

# Indice

In copertina: disegno tratto dalla Relazione presentata da L. Tatti e C. Clericetti in: *Nuovo tempio israelitico in Torino*, Tipografia C. Favale e Compagnia, Torino, 1874, Tavola IV; firma autografa di Celeste Clericetti.

## Presentazione della corrispondenza

La corrispondenza qui riprodotta è composta da 48 lettere scritte da Clericetti a Cremona e 10 lettere scritte da Cremona a Clericetti che trattano perlopiù di questioni personali. Da esse si evince soprattutto il profondo legame di amicizia che esisteva tra i due. Frequenti sono le considerazioni e i giudizi su situazioni e persone (per esempio su Francesco Brioschi) destinati evidentemente a restare nell'ambito della più stretta confidenza. Di più, le due famiglie si trovavano spesso per trascorrere insieme le vacanze e Cremona, quando nel 1873 si trasferì a Roma per dirigere la locale Scuola di Applicazione per gli Ingegneri, fu proprio ai Clericetti che "affidò" il figlio Vittorio rimasto in collegio a Milano.

Clericetti fu ingegnere, appassionato di monumenti antichi e attento osservatore della realtà. Nelle sue lettere spesso racconta - arricchendo la narrazione con osservazioni molto puntuali - di scoperte archeologiche e di cantieri visitati, tanto che, da queste note, si può ricostruire la nascita di parecchi monumenti, opere pubbliche ed eventi (da Palazzo Montecitorio a Palazzo Madama, dal traforo del Gottardo al Teatro Massimo di Palermo, dalla Mole Antonelliana al primo forno crematorio nel Cimitero Monumentale di Milano, dall'Esposizione Nazionale a Milano del 1881 a quella di Roma del 1883, ecc.).

Altrettanto interessanti sono le notizie e i commenti sulla vita dell'Istituto Tecnico Superiore di Milano: l'avvicendamento dei professori, le gite con gli studenti del 3° anno, le discussioni con Brioschi e così via.

Alcune lettere trattano in particolare della situazione nella quale Cremona si trovò come direttore dell'Istituto di Applicazione per gli Ingegneri di Roma e del desiderio di Brioschi di far chiudere quella scuola.

Parecchie sono inoltre le notizie biografiche dei personaggi che gravitavano attorno al Politecnico o al Reale Istituto Lombardo: Clericetti infatti ne dava pronta notizia a un lontano Cremona. Così, accanto agli annunci di morte, si trovano accenni alle difficoltà nelle quali versò Brioschi nel 1885 o considerazioni sull'operato di diversi Ministri.

Anche dalle lettere di Cremona appare la grande confidenza esistente tra i due, basata su una profonda stima. Cremona scrive in inglese con l'esplicito intento di fare esercizio di lingua e racconta all'amico i suoi progetti di vacanza o commenta alcuni avvenimenti legati all'Istituto di cui è direttore. Buona parte della corrispondenza è dedicata alla conta dei voti nelle adunanze del R. Istituto Lombardo in una prima occasione per concedere a Cremona la pensione e in seguito per far nominare Clericetti membro effettivo.

# Criteri di edizione

**a)** <u>lettere scritte da Clericetti</u>
Poiché le lettere sono conservate in buono stato, i criteri adottati per la trascrizione hanno mirato a rispettare il più possibile il testo originale sia nella grafia (ad es. molte parole iniziano con la lettera maiuscola e così sono state lasciate), sia nella punteggiatura. Le uniche due eccezioni riguardano: la trascrizione di parole maschili precedute da articolo indeterminativo: Clericetti usava inserire l'apostrofo che invece è stato ignorato; la trascrizione di parole con "ì" finale (ad es: così, costì): Clericetti non metteva sempre l'accento sulla lettera i, che invece è sempre stato inserito.
Si è mantenuta la sottolineatura "a parola" del manoscritto.
Si sono mantenuti termini oggi desueti, ma in uso all'epoca, sono stati segnalati invece con *[sic!]* termini e periodi che appaiono scorretti o comunque di senso poco chiaro. I nomi propri sono stati riprodotti come scritti dall'autore, anche se inesatti e sono stati corretti soltanto in nota.
Non è stato possibile decifrare alcune parole e in tali casi si è lasciata tra parentesi quadre un'interpretazione plausibile, oppure si è sostituita la parola con il segno *[...]*. La lettera 089-20888 risulta molto danneggiata: si sono trascritte le poche frasi leggibili e si è deciso di sostituire le parti non recuperabili con la scritta *[illeggibile]*.
Clericetti era nato a Londra e alcune lettere sono scritte in inglese (sono sei, dal 1879 al 1881): da quanto si evince dai testi, si trattava di un modo per fare esercizio concordato con Cremona.
Tutte le scritte in corsivo tra parentesi quadre sono state inserite dal curatore.

**b)** <u>lettere scritte da Cremona</u>
Si tratta di minute che Cremona conservò. Sono scritte in inglese e corrette dalla persona - Mary Isabella Taddeucci - che probabilmente a Roma insegnava la lingua a Cremona: la scelta è stata quella di trascrivere i testi corretti, come si suppone che Cremona li abbia inviati al suo corrispondente o, nel caso siano presenti più opzioni, si è scelto di trascriverle tutte.

# Celeste Clericetti (Londra, 20/11/1835 - Milano, 30/5/1887)

COMMEMORAZIONE
DELLA VITA E DELLE OPERE DI CELESTE CLERICETTI

**LETTURA**

fatta al Collegio degli Ingegneri ed Architetti dal Prof. A. SAYNO
nell'Adunanza del 17 Marzo 1889
in occasione della inaugurazione del ricordo monumentale che venne eretto
per pubblica sottoscrizione nella Sede del R. Istituto Tecnico Superiore di Milano.*

Signori!

La presidenza del Collegio degli Ingegneri ed Architetti volle affidarmi l'onorevole incarico di commemorare la vita dì CELESTE CLERICETTI. Molti e meglio di me avrebbero potuto dire di quella nobile esistenza, la quale fu una sintesi di alte aspirazioni e di amore intenso allo studio e al lavoro; io, che non ho potuto esimermi dall'accettare il non facile compito, ravviverò come meglio potrò nella vostra memoria il caro maestro e collega perduto; vogliate voi ascoltarmi benignamente.

Nel giorno primo del Giugno 1887 uno stuolo di persone di ogni classe accompagnava all'ultima dimora la salma di Celeste Clericetti. Milano aveva perduto uno de' suoi Cittadini valenti ed operosi; Milano che studia e lavora, riconoscente, rendeva alla sua memoria un ben degno tributo di affetto e di stima[1]; allora mi era sembrato che l'industriale e l'operaio, lo studente e lo scienziato, il tecnico e l'artista, i quali tutti nei loro sodalizi avevano conosciuto ed apprezzato il compianto professore, seguissero quel feretro non per le usate ragioni di convenienza, ma che ivi li guidasse un intimo convincimento, quasi vorrei dire un bisogno del cuore. Celeste Clericetti struggendosi per lo studio con una febbrile attività aveva dianzi tempo finito la sua giornata, dopo una vita tutta spesa nello strenuamente volere e nel largamente donare per l'altrui beneficio quanto ha potuto, pago per sé solo di avere sempre e fortemente sentita la religione del dovere e nobilmente aspirato, colla sua mente di scienziato e di artista, agli altissimi ideali del vero e del bello che furono tormento all'animo suo e voto insoddisfatto delle sua breve esistenza.

Egli ci lasciava per sempre, ma colla fede di sopravivere caramente nella memoria di chi ebbe da Lui istruzione ed affetti, aiuto e consigli. Triste sarebbe per gli uomini quel giorno in cui si affievolisse questo culto altamente educatore e morale che noi sentiamo di dover rendere a chi ha sacrificato ogni potenzialità della propria mente al duplice scopo di arricchire e diffondere la pubblica cultura, la quale è tanta parte del patrimonio vivo delle nazioni. Ma voi che foste discepoli, colleghi ed amici del compianto professore avete oggi dimostrato che questo culto è in voi vivissimo inaugurando in questo recinto dedicato agli studi - dove un giorno risuonava la concitata e convinta parola del maestro perduto - il bronzo che ritrae con mirabile fattura le di Lui sembianze, dalle quali sembra che irradi l'animo suo, quale era, franco e leale.

Celeste Clericetti nacque il 20 Novembre 1835 a Londra, nella quale città la famiglia allora trovavasi per ragioni di commercio; Egli, compiuti gli studi tecnici secondari in Milano, nel 1853 passava all'Università di Pavia, dove consacrava il tempo che non gli veniva assorbito dalle matematiche discipline nell'affinare la sua cultura letteraria e storica, nell'apprendere le lingue moderne e nel coltivare lo studio dei classici in cui fu maestro a sé stesso. Laureato in matematica nel 1856, si diede tosto alla pratica dell'ingegneria, non trascurando in pari tempo di perfezionarsi negli studi architettonici presso l'Accademia di Belle Arti in Milano, ove nella composizione venne distinto con medaglie d'argento,

AI giovane Clericetti che aveva sì bene incominciato erano facilmente prevedibili i compensi di una brillante carriera tecnica; ma in Lui dominava imperioso il bisogno dello studio e per meglio assecondare sì nobile impulso, alieno da ogni personale interesse, troncava la pratica per dedicarsi all'insegnamento, nel quale esordiva come assistente alla cattedra di architettura civile nella R. Università di Pavia, posto che

---

* *Il Politecnico - Giornale dell'ingegnere architetto civile ed industriale*, v. 21, 1889, pp. 354-368.
[1] Giornale *La Perseveranza* del 2 Giugno 1887: *In Memoria del* Prof. Ing. CELESTE CLERICETTI «*Discorsi pronunciati davanti al feretro e tributi di omaggio resi alla di Lui memoria*», Milano tipografia di C. Rebeschini e C., 1887.

però lasciava dopo un biennio per stabilirsi a Milano come professore di disegno geometrico e di macchine nel R. Istituto Tecnico a S. Marta, al quale incarico nel 1863 vi aggiungeva quello di insegnante la geometria descrittiva e il disegno industriale presso la Società di Incoraggiamento di Arti e Mestieri,

Il 1863 segnava una data memorabile per gli studi nella città di Milano. Lo Stato spinto dalla tenace iniziativa del senatore Brioschi fondava in quell'anno il Politecnico milanese, e tra i pochi e distinti cooperatori dei quali volle circondarsi questo nostro illustre concittadino per dare vita all'opera sua e che con lui divisero la fede e compresero l'alto intendimento di questa istituzione troviamo Celeste Clericetti, al quale il Direttore della nuova scuola affidava da prima il disegno di costruzioni e in seguito nel 1864 l'insegnamento della scienza delle costruzioni, posto che Egli coprì con grande onore e che troppo presto dalla morte doveva venirgli strappato! Celeste Clericetti comprese l'importanza dell'arduo compito assunto e chi ha presente lo stato della coltura tecnica in Italia 25 anni or sono, può ora giudicare ed apprezzare quanto sia stata utile e feconda di bene l'opera sua.

Le tradizioni dell'arte edificatrice italiana, della cui storia il compianto professore era tanto erudito, l'ambiente in cui ebbe la sua educazione scientifica, l'indirizzo ricevuto dai suoi maestri ed una certa predilezione artistica verso la quale si sentiva attratto, tutto consigliava al prof. Clericetti di continuare nella nuova scuola i metodi didattici già in uso nelle Accademie di Belle Arti d'allora, nelle quali, la scienza delle costruzioni racchiusa in brevissimi confini, senza una base razionale, trovavasi vagante, confusa coll'estetica ed impotente ad affermarsi con una vita propria. Ma il prof. Clericetti ebbe la percezione dei tempi nuovi e capì che si doveva cambiare strada. Egli, nei suoi frequenti viaggi in Europa e segnatamente in Inghilterra, poté constatare il sorgere gigante delle industrie metallurgiche e il rapido diffondersi dell'uso del ferro, della ghisa e dell'acciajo nelle costruzioni di ogni genere, per cui si creano di continuo nuove ed ardite strutture nell'arte edilizia e nelle opere attinenti alla pubblica viabilità, imperiosamente volute dai bisogni della Società moderna, la quale, schiava del tornaconto, domanda che si estenda sempre più l'uso del ferro anche in quelle costruzioni nelle quali l'impiego di questo metallo, non da necessità richiesto, porta con sè l'obblio e la decadenza dell'arte pura.

In quello stesso tempo nelle scuole tecniche superiori della Francia, Germania, Svizzera e nei laboratori di diversi sperimentatori di quei paesi e di Inghilterra prendevano largo sviluppo alcuni studi teorici-sperimentali sulle resistenze dei materiali intesi a trovare delle formule, le quali fossero atte a rappresentare, almeno per approssimazione, le dipendenze che esistono tra le forze esterne applicate alle costruzioni, le deformazioni e le resistenze interne che vi corrispondono, e ciò allo scopo di dare al tecnico una guida nel calcolo della stabilità delle opere di resistenza e segnatamente delle grandi costruzioni metalliche, per le quali l'esperienza del passato a nulla serviva, nulla diceva.

Il prof. Clericetti ebbe una pronta e retta intuizione dei bisogni della sua scuola e, lasciato il posto di professore a S. Marta, si gettò a capo fitto in questa nuova corrente di studi, già poderosa all'estero e quasi nulla da noi, coi quali gli stranieri avevano già saputo speculare nel nostro paese, invaso da progettisti e costruttori di ogni genere, i quali ci imponevano i loro lavori perché allora le nostre cognizioni tecniche nulla di meglio sapevano o potevano produrre. La grande tettoia alla stazione centrale, la copertura della Galleria Vittorio Emanuele ed altre opere di Milano che datano da quelle epoche, sono frutti dell'industria straniera[2], alla quale però, oggi si sa contraporre, dopo un breve corso di anni, col sussidio della scienza e dell'industria italiana opere grandiose come il ponte sul Ticino a Sesto Kalende, il ponte di Paderno e quello che ora si sta costruendo sul Po nella provincia di Cremona.

Fu nei primi anni di insegnamento al Politecnico che il Prof. Clericetti dovette reggere ad un lavoro intensissimo, capace di fiaccare tempre ben più robuste della sua. Egli, privo di tradizioni scolastiche e solo colla scorta di lavori e studi sparsi all'estero dovette organizzare due insegnamenti, uno di resistenza dei materiali e l'altro di ponti e insegnare a numerose scolaresche le moderne teorie delle quali si arricchiva ogni giorno la sua mente con uno studio indefesso. E come ciò fosse poco sovrintendeva con un orario pesantissimo a due scuole di disegno di costruzioni civili e di ponti frequentate in alcuni anni da più di 80 allievi per corso, e non aveva finito: alla sera, alle domeniche, altre lezioni orali, altre scuole di disegno da dirigere con numerosi allievi presso la Società di incoraggiamento di arti e mestieri: incarico quest'ultimo al quale attendeva con zelo ma che volentieri avrebbe rinunciato per meglio dedicarsi agli studi da Lui prediletti se le sue condizioni finanziarie lo avessero permesso. Sembrava che col crescere del lavoro si invigorisse sempre più la lena in quell'uomo, nella cui mente eravi una preoccupazione continua ed esclusiva per lo studio e la scuola.

---

[2] La tettoia della stazione centrale di Milano venne calcolata dall'Ufficio d'arte della Parigi-Lione- Mediterraneo (1862) e costruita a Parigi. La grande copertura metallica della Galleria Vittorio Emanuele venne costruita dal Sig. Henry Joret di Parigi nel 1867.

Sono molti i giovani ingegneri che devono ricordare con affetto il nome di sì caro maestro che li istruì con vera passione e li amò sotto le apparenze di un burbero severo, mentre in Lui dettavano legge un cuore d'oro ed un animo gentile. Celeste Clericetti non tardò molto a rendersi pienamente edotto del moderno movimento scientifico relativo al suo insegnamento ed a presentare all'esame degli studiosi i frutti delle sue ricerche. Tra le numerose memorie da Lui lasciate che si riferiscono alla Scienza delle Costruzioni sonvene alcune pregevoli per erudizione e per ricerche originali, le quali hanno attinenze colla parte sperimentale e colla teoria.

Dietro iniziativa del prof. Clericetti, il Gabinetto di Costruzioni nel R. Istituto Tecnico Superiore, da Lui coordinato, si arricchiva nel 1866 di una macchina Clair, la prima introdottasi in Italia, per sperimentare la resistenza dei materiali alla trazione ed allo schiacciamento. In quell'epoca tale macchina era considerata come una delle migliori, analoga a quella che aveva servito al Signor Michelot nel Conservatorio di Parigi per eseguire 3000 esperienze, molto accreditate, sulle pietre naturali della Francia[3]. Il Prof. Clericetti coadiuvato dall'egregio Prof. Loria si proponeva di sperimentare la resistenza dei principali materiali italiani onde stabilire dei criteri generali sul modo di rottura, sul grado di elasticità, sull'influenza del numero dei pezzi sovrapposti e del rapporto fra le loro altezze e la base compressa, in vista delle diverse applicazioni ai muri degli ordinarii edificii, ai piedritti delle volte, ecc. A queste esperienze dovevano fare seguito altre ricerche relative alla durabilità dei materiali rispetto alle diverse azioni degradatrici degli agenti atmosferici onde, infine, poterli classificare convenientemente nell'interesse della pratica. Ognuno può comprendere quanto fosse importante questo programma di esperienze, le quali dovevano recare tante notizie utili, in allora sconosciute o confusamente note ai nostri costruttori, i quali, in merito alle proprietà di resistenza dei materiali che usavano nelle loro opere, non avevano che nozioni incerte, apprese da tradizioni ereditarie. Ma tale programma per la sua vastità non venne realizzato che in parte, sebbene le prove venissero continuate attivamente per ben tre anni di seguito. Sono circa 900 le esperienze note al pubblico, che allora si eseguirono per constatare la resistenza allo schiacciamento delle pietre naturali dell'Alta Italia e dei composti artificiali a base di cemento e di calce[4], e tra queste esperienze sono da annoverarsi come originali e di speciale importanza quelle relative all'elasticità dei calcestruzzi. I risultati sperimentali sulla resistenza dei calcestruzzi pubblicati con tanta ricchezza di particolari e di confronti relativi alla loro composizione, stagionatura e gradi di resistenza se non possono avere una importanza rigorosamente scientifica, perché i mezzi dei quali disponeva il Prof. Clericetti non potevano rispondere a tale scopo, hanno però fornito delle cognizioni utilissime alla pratica e contribuito a diffondere tra noi l'uso di questi materiali artificiali, a far conoscere ai costruttori le attitudini di resistenza degli ottimi calcari idraulici dei quali è ricca la Lombardia, ad attivare diverse industrie su vasta scala ed a liberarci dalla importazione di molti cementi esteri coi quali gareggiano, e con pieno successo, i nostri prodotti.

A complemento di queste esperienze il Prof. Clericetti eseguì con molta cura e con speciali apparecchi da Lui ideati diverse prove sopra un volto di calcestruzzo della corda di 8 metri; la monta di 2, lo spessore in chiave di met. 0,10 prove, colle quali Egli volle dimostrare quanto siano adatti i calcestruzzi per simili opere per la grande resistenza di cui sono capaci, il loro modico prezzo e la speditezza con cui si impiegano. Il Prof. Clericetti, fatta l'ipotesi di una uniforme pressione in chiave dimostrò, ammesse alcune speciali condizioni, quanto fosse soddisfacente l'accordo tra la spinta teorica e quella determinata sperimentalmente[5]. Altre numerose esperienze eseguì il Prof. Clericetti sopra diversi materiali da

---

[3] BELGRAD. - «*Rapporto sulle dette esperienze pubblicate negli Annales des Ponts et Chaussées*», 1855, 2.⁰ sem. pag. 189.

 MICHELOT P. - «*Espériences sur la résistance des matériaux à l'écrasement*», Paris, 1863. *Annales du Ponts et Chaussées*. 1863, Tom. V.

[4] CLERICETTI C. e LORIA L. - «*Esperienze sui materiali da costruzione*», Prima serie, Anno 1869, Milano, Tip. e Lit degli Ingegneri.

 CLERICETTI Prof. CELESTE. - «*Esperienze sui materiali da costruzione*», Seconda serie «*Esperienze sui calcestruzzi*», Anno 1871, Milano Tip. e Lit. degli Ingegneri.

[5] *Politecnico*; anno 1871.

costruzione dietro incarico avuto da varie amministrazioni e da privati, iniziando così un pubblico servizio pei costruttori, della cui utilità nessuno può dissentire[6].

Il problema sperimentale della determinazione della resistenza dei materiali, per quanto si voglia limitare ai principali bisogni della pratica, è estremamente complicato e bastano alcune insensibili variazioni nel metodo di prova, o nello stato e forma dei pezzi per ottenere risultati tra loro disparatissimi per gli identici materiali. Le esperienze del Prof Clericetti, al pari di quelle eseguite da altri illustri sperimentatori moderni tra i quali il Bauschinger[7] ed il Tetmajer[8], hanno messo in evidenza queste difficoltà e quindi la necessità di studiare con maggior diligenza e con processi più esatti l'importante problema, dal quale molta luce può irradiare anche nel campo della fisica molecolare.

I continui perfezionamenti ed anche le novità di ogni genere che la tecnica moderna sa introdurre nella fabbricazione dei laterizi, delle pietre artificiali, dei ferri, acciai, ghise, leghe metalliche, ecc. e la facilità con cui attualmente le Provincie e gli Stati si scambiano tra loro i propri materiali da costruzione naturali ed artificiali, tutto ciò ingenera una certa confusione nel campo pratico dei costruttori, i quali difficilmente possono difendersi dalle frodi abilmente mascherate e conoscere con esattezza le proprietà di resistenza e di durabilità di questi materiali, pel cui accertamento si richiedono, come già si disse, procedimenti non facili, mezzi meccanici ed istrumenti di precisione i quali non possono essere alla portata di tutti. L'arte edilizia, segnatamente nei grandi centri di popolazione, inceppata dalla speculazione, avara dello spazio, tende continuamente a ridurre le sezioni delle principali membrature di resistenza nelle strutture degli edifici, ed anche per questo motivo si rende evidente la necessità di dare un vigoroso e nuovo impulso alle ricerche sperimentali relative alla resistenza e durabilità dei materiali, mettendo questi in quelle condizioni statiche e dinamiche che meglio si avvicinino alle condizioni effettive che si raggiungono nella messa in opera, e ciò allo scopo di racchiudere in più ristretti confini le indeterminazioni sui limiti di resistenza che contribuiscono a rendere incerto e molte volte arbitrario il giudizio del tecnico sul grado di stabilità di una data costruzione, anche in quei casi semplici nei quali il problema teorico può ritenersi completamente determinato.

A capo della schiera di filosofi e tecnici che si sono occupati degli studi sperimentali e teorici sulla resistenza dei materiali sta il nome del sommo Galileo, del Fabri, del Marchetti, del Padre Grandi e di altri italiani del 17.° e 18.° secolo[9]. I loro studi però rimasero per molto tempo isolati e la preziosa eredità venne più tardi raccolta dagli stranieri: Hodgkinson, Barlow, Fairbairn, Clark, Kirkaldy, ecc. in Inghilterra[10]; Gerstner in Germania[11]; Rondelet, Vicat, Love, Morin, Tresca, ecc. in Francia[12] diedero, prima del 1869, un grande sviluppo a queste ricerche sperimentali con speciale indirizzo ai bisogni della pratica. In Italia questi studi furono debolmente ripresi solo verso il 1841, ed in un periodo di pressoché 20 anni non possiamo

---

[6] È di particolare interesse il rapporto: «*Intorno ad alcune esperienze sulla resistenza delle volte in cemento*» pubblicato dagli ingegneri Clericetti e Tatti nel *Politecnico*, anno 1875.

[7] *Mittheilung des mech. Tech. Laboratoriums*, Munchen.

[8] *Mittheilungen der Anstalt zur Prüfung von Baumaterialien am eid. Polytechnikum in Zürich,* Zürich, 1884, und f.

[9] *Sunto delle lezioni del Navier; prefazione storica del Saint Venant*. Paris, 1864.

[10] Di «Eaton Hodgkinson» vi sono molti rapporti sulle resistenze dei ferri e delle ghise, tubi, colonne, ecc. inseriti nelle «*Philosophical transactions*» anno 1840 e seg. È interessante la memoria «*Esperimental researches on the strength and the other properties of cast iron*» London, 1842. Altri lavori furono pubblicati dalla «*Associazione Britannica sul progresso delle scienze*» anno 1846 e seguenti.

Barlow Peter - «*A treatise on the strength of timber, cast and malleable iron and other materials*» London, 1851.

Le esperienze di *Clarck*, *Fairbain*, alcune delle quali datano sino dal 1838, non che quelle di *Hodgkinson* e di altri sperimentatori inglesi sono ricordate dal *Morin* e dal *Love* nelle loro opera pubblicate anteriormente al 1869.

Nell'*Enginer* 1869 e seguenti, si trovano inseriti i risultati di parecchie esperienze fatte a Londra nell'officina *Kirkaldi*.

[11] Citato dal *Weisbach* nel «*Leherbuch des Ingenieur und Maschinen Mechanik*».

[12] Rondelet - «*Traité théorique et pratique de l'art de bâtir*» Paris, 1812-1814.

Vicat – Vi sono molti rapporti di esperienze sui metalli, i tubi, le pietre, le malte, ecc. inseriti negli «*annales de Ponts e Chaussées*» anno 1831 e seg. e memorie distinte pubblicate a Parigi negli anni 1828, 1830, 1836, 1847 ed a Grenoble nel 1856.

G.H. Love - «*Des diverses Résistances et autres propriétés de la fonte, su fer et de l'acier*», Paris, 1852, Lacroix, Editeur.

Morin - «*Résistance des matériaux*» Paris, 1862.

citare come interessanti che alcune esperienze del Prof. Giulio[13], al quale fecero seguito le ricerche dei signori Paccinotti e Peri[14], quelle dell'Ufficio Tecnico del Real Cantiere di Castellamare[15], del generale Cavalli e dell'Ing. Noè[16]. Ma il primo che in Italia diede un largo sviluppo a queste esperienze, le quali contribuirono a richiamare maggiormente l'attenzione dei tecnici sulla loro importanza, fu il Prof. Clericetti, ed è certo che la sua iniziativa giovò a diffondere questi studi sperimentali nelle principali Scuole d'Applicazione dello Stato dove sono in continuo progresso, e in alcuni grandi cantieri da costruzione delle amministrazioni pubbliche e dei privati.

Il Collegio degli Ingegneri, l'Accademia di Belle Arti di Milano ed il Ministero della Pubblica Istruzione che apprezzavano altamente la competenza del Professor Clericetti lo chiamavano sempre a far parte di speciali commissioni delle quali era quasi sempre il relatore, ogni qual volta, tra le altre, si dovevano risolvere questioni di interesse scientifico, le quali avessero attinenze col ramo degli studi sperimentali da Lui coltivato. Citeremo ad esempio: a) lo studio sull'influenza del gelo sui materiali da costruzione; b) l'ordinamento della raccolta dei materiali da costruzione per l'esposizione di Milano, che ebbe un importante successo; c) gli studi risguardanti il consolidamento dei materiali da costruzione degli antichi edifici deperiti per opera del tempo e d) finalmente gli studi per formare una collezione di prodotti minerali italiani per uso delle Arti Edilizie.

Le esperienze di Wöhler incominciate a Berlino sino dal 1858 e continuate dallo Spangenberg[17] relative all'influenza della ripetizione rapida degli sforzi di tensione e di pressione sopra sbarre di ferro e di acciaio, esperienze che destarono un vivissimo interesse nel campo dei tecnici e degli studiosi, diedero occasione a diversi distinti autori, quali il Launhardt, il Gerber, il Weyrauch, il Ritter[18], ecc., di proporre alcune formole empiriche atte a rappresentare con approssimazione le leggi di Wöhler, lasciando quasi supporre che si trattasse di una nuova legge fisica non ancora determinata. Il Prof. Clericetti basandosi sopra alcune recenti esperienze fatte dal Tresca al Conservatorio di Parigi, e confermate dal Bauschinger, relative all'aumento del limite primitivo di elasticità che può venire spinto sin presso al modulo di rottura, pur mantenendosi quasi costante il modulo di elasticità[19], spiega i resultati di Wöhler come conseguenza della nota proprietà che una tensione istantanea applicata ad una barra vi produce delle oscillazioni tali che la più ampia corrisponde ad un allungamento doppio di quello che vi produrrebbe la stessa tensione se aumentasse per gradi insensibili da zero sino a raggiungere il valore finale; e di fatti la formula semplicissima[20] che il Prof. Clericetti dedusse da tale legge rappresenta in generale e con sufficiente approssimazione i numeri di Wöhler. I coefficienti specifici pel ferro calcolati colla formula del Prof Clericetti, la quale ha una base teorica, nella maggior parte dei casi differiscono al più del 4 % dai valori medi che si ottengono da nove formule empiriche proposte dai detti autori e da altri fisici e tecnici[21]. Sarebbe forse troppo l'asserire che la formula di Clericetti rappresenta il vero fenomeno di meccanica

---

[13] GIULIO - «*Expériences sur la résistance à la flexion et sur la résistance à la rupture des fers forgés dont on fait le plus d'usage en Piémont*», Accademia delle Scienze di Torino, Serie II, Tomo III, anno 1840. Altri rapporti furono inseriti negli atti della stessa accademia negli anni 1841, 1842. Le stesse esperienze sono citate dal BORGUIS nei suoi «*Elementi di statica architettonica*», Milano 1842.

[14] Rapporti di esperienze sui legnami inseriti nel «*Cimento*» 3.°(1845). Sono anche citate dal GABUSSI nel trattato «*L'Arte del Costruttore*» 1866, Milano.

[15] PULLINO GIACINTO - «*Corso di resistenza dei materiali* » Castellamare, 1866.

[16] CURIONI - «*Resistenza dei Materiali*» Torino, 1863.

[17] SPANGENBERG - «*Uber des Verhalten der Metalle bei Wiederholten Austrengungen*» Berlin, 1875.

[18] LAUNHARDT - «*Die Inaspruchnahme des Eisen, ecc. Zaitschr. Des Hanov. Arch. Und Ing. Vereins*» 1873.

GERBER - «*Bestimmung der Zulässingen Spannungen in Eisen constructionen. Zeitschrift der Bair. Arch. Und Ing. Vereins*», 1874.

J. WEYRAUCH - «*Festigkeit und Dimensionen Berechung der Eisen und Stahl Constructionen, ecc.*» Leipzig, 1876.

[19] TRESCA - «*Résultats des expériences de flexion faites sur des rails en fer et en acier au de là de la limite d'élasticité et jusqu'à la rupture*» Mémoires de la Société des Ingénieur Civils, 1880.

[20] CELESTE CLERICETTI - «*Sulla determinazione dei coefficienti di sforzo specifico, dietro le esperienze di Wöhler*» Politecnico. Anno XXIX, Ottobre-Novembre, 1881.

[21] Oltre alle formule degli autori citati alla nota (18), Clericetti considerò altre formule citate dal Mohr e pubblicate nel *Civilingenieur*, 1881, fasc.1°.

molecolare che trovasi implicato nelle esperienze di Wöhler, ma a noi pare che questa formula abbia un certo valore e tale almeno da reggere il confronto colle formule empiriche suaccennate[22].

Ecco, in breve, quanto ha fatto il Prof. Clericetti pel progresso degli studi sperimentali sulla resistenza dei materiali, e se a Milano, città che ha il vanto di tante iniziative, si potesse istituire una stazione di prova, un vero laboratorio, per eseguire su vasta scala e coi mezzi perfezionati e potenti che ora si conoscono le esperienze relative alla resistenza dei materiali da costruzione, come già da tempo si fa nelle principali capitali di Europa, sarà sempre dovuta una speciale benemerenza alla memoria del Prof. Clericetti, il quale iniziò nella nostra città queste ricerche con una ricca serie di lavori, che ben meritano un cenno distinto nella storia contemporanea degli studi sperimentali nel nostro paese.

Quando il Prof. Clericetti iniziava nel R. Istituto Tecnico Superiore l'insegnamento della scienza delle costruzioni, la teoria dell'equilibrio dei sistemi elastici che ne è la base professavasi nei principali Politecnici di Europa seguendo metodi che il Lamè ebbe a definire come «semianalitici e semi empirici» perché dedotti da ipotesi dimostrate attendibili in pochi casi e da esperienze istituite in condizioni particolarissime e non atte a rappresentare con sufficiente rigore il fenomeno complesso delle deformazioni elastiche, dalle quali dipendono le interne resistenze dei solidi; allora nei corsi scolastici si era ancora ben lontani dal ritenere che la teoria matematica dell'elasticità si potesse considerare come la naturale introduzione al corso teorico di scienza delle costruzioni, perché i mirabili risultati di questa teoria, frutto degli studi di tante menti superiori[23], presentavano difficoltà analitiche insuperabili o conducevano a calcoli laboriosi e troppo complessi se si volevano portare nel campo delle applicazioni, anche attenendosi ai casi i più semplici della pratica; difficoltà però che ora vanno mano mano appianandosi in seguito ai recenti studi del Castigliano, seguiti da quelli del Crotti, del Weyrauch, del Müller-Breslau, del Melan e di altri distinti autori[24]. Anche gli studi di Carlo Culmann che arricchirono la scienza delle costruzioni di metodi grafici utilissimi per la risoluzione di molti problemi complessi di statica e che ebbero a diffondersi rapidamente in Europa e negli Stati Uniti di America, 25 anni or sono erano appena iniziati nel Politecnico di Zurigo[24bis].

Il Prof. Clericetti nello svolgere il suo corso di resistenza dei materiali e quello di applicazione ai ponti non si scostò dai metodi allora in uso ed illustrati dal Morin, dal Bresse, dallo Scheffler, dal Collignon, dal Ritter e da altri autori; però in qualche argomento Egli seppe generalizzare diverse teorie, in altri mettere in evidenza gli accordi e le equivalenze degli svariati processi di calcolo seguiti dagli autori e sotto questo riguardo alcune monografie del Professor Clericetti hanno reso un utile servizio agli allievi della sua

---

[22] Questi studi del prof. Clericetti vennero apprezzati anche all'Estero e segnatamente a Bruxelles dove furono pubblicati nella Revue Universelle des Mines (1882), A Parigi il Collignon li rammenta nella sua opera «*Résistance des Matériaux*»; III Ediz., 1885.

Il CASTIGLIANO nel suo «*Manuale pratico per gli Ingegneri*» limita d'assai l'importanza delle formule empiriche tedesche; anche la Società degli Ingegneri Civili di Parigi nella seduta del 21 luglio 1881 opinò che l'interpretazione data dagli autori tedeschi alle esperienze di Wöhler lascia molti dubbi e non ne accetta nella loro generalità, le conclusioni.

[23] TODHUNTER AND PERSON - «*A history of the elasticity (1839-1850)*» Cambridge, at the University press. An 1866.

Vedere inoltre pei più recenti studi le opera del LAMÉ, BARRÉ DE D. VENANT, GRASHOF, CLEBSCH, THOMSON e TAIT, KIRCHOFF, NEUMANN e le pubblicazioni degli illustri professori italiani BETTI, BELTRAMI, TURAZZA, CERRUTI.

[24] A. CASTIGLIANO. - «*Théorie de l'équilibre des systèmes élastique et ses applications*», Torino A.F. Negro.

F. CROTTI. - «*La teoria dell'elasticità nei suoi principii fondamentali e nelle sue applicazioni pratiche alle costruzioni*» Ulrico Hoepli, Milano, 1888.

WEYRAUCH. - «*Theorie elasticher Körper*» Lipsia Teubner (1884).

MÜLLER-BRESLAU - «*Uber die Anwendung des Princips der Arbeit in der Festigkeitslehre, Wochenblatt f. Arch. und Ing. An.*» 1883.

MELAN. - «*Uber den Einfluss der Würme auf elastiche Systeme*» Wochenschrift des oester Ing. und arch. Vereins, 1881.

- «*Beitrag zur Berechnung statisch unbestimmter stabsysteme*» Zeitschrift des Oester Arch. und Ingenieur, Vereins, 1884.

CANEVAZZI. - «*Sulla teoria delle travature. Monografia, Atti della R. Accademia delle Scienze dell'Istituto di Bologna* (1886).

[24bis] MEYER. - *Le Dr. Charles Culmann ingénieur et professeur à l'École polytechnique fedérale à Zurich* (1882).

FAVARO. - «*Della vita e degli scritti di Carlo Culmann*» Venezia (1882). La prima edizione della «*graphische Statik*» venne pubblicata verso la fine del 1865.

scuola ed ai tecnici; pei quali, in generale, non è facile cosa trovare la soluzione di molti problemi complessi che si presentano nella pratica sulle pubblicazioni sparse dei giornali scientifici; perché la sola differenza di metodo può essere d'imbarazzo alla pronta intelligenza di un dato argomento.

Il corso autolitografato di scienza delle costruzioni pubblicato dal Prof. Clericetti sino dal 1866, sebbene non fosse privo di inesattezze, scusabili nella compilazione affrettata di un lavoro di una certa importanza, fu per molti anni una buona guida per la scuola e la pratica[25]. In allora le lezioni del Prof. Clericetti ed il corso di resistenza dei materiali del Prof. Curioni si può dire che compendiassero tutto quanto di più recente si conosceva in Italia in ordine a questi studi teorici sulla stabilità delle costruzioni.

Omettendo per brevità di parlare di alcune note scientifiche di minore importanza[26], non possiamo dimenticare i due studi interessanti del Prof. Clericetti che si riferiscono alle travature reticolari ed alla teoria delle volte. In una prima pubblicazione fatta nel 1866 sulla teoria elementare delle travature reticolari[27], il prof. Clericetti volle riassumere, coordinare tra loro i diversi metodi di calcolo seguiti dal Ritter, dal Lentz dal Culmann e dal Goudard[28] nel risolvere i problemi che interessano questo ramo importante della Scienza delle Costruzioni ed appianare ogni difficoltà, sicché, come dice Lui, i tecnici e gli studiosi fossero in grado di poter calcolare da soli le travature nei casi particolari della pratica. Questa monografia del Prof. C!ericetti è ricca di molti esempi e nella parte che riguarda l'influenza dei carichi mobili vi sono alcune osservazioni di rilievo, le quali possono servire a semplificare in alcuni casi la ricerca dei massimi e minimi sforzi nelle singole membrature di contorno e di collegamento che compongono le dette travature. A questo studio teorico sulle travature reticolari fecero seguito più tardi altri lavori affatto originali i quali si riferiscono al calcolo delle travature composte e alle applicazioni relative ai moderni ponti sospesi americani[29].

Si tratta di quelle meravigliose ed ardite opere complesse, create dal senso pratico degli ingegneri americani, i quali combinando tre strutture semplici costituite da un sistema funicolare, da uno elastico e da un terzo articolato, hanno saputo costruire delle travate di piccola altezza e di considerevoli lunghezze, delle quali è splendido esempio il ponte che attraversa la Riviera dell' Est fra New-York e Brooklyn, della luce libera di m.493. Il Prof Clericetti dettò una teoria ingegnosa per questi sistemi la quale, sebbene non sia che approssimata, perché basata sopra ipotesi sulla cui influenza non si possono fare apprezzamenti completamente determinati, dà però dei risultati, i quali si accordano in un modo soddisfacente colle prove sperimentali citate dal Malezieux, il quale in un rap- porto al Ministero dei Lavori Pubblici in Francia, descrisse i principali manufatti degli Stati Uniti d'America esistenti nel 1870[30]. In questo rapporto si asserisce, ad esempio, che il ponte alla Niagara-Falls, costruito sino dal 1855, della luce libera di Metri 250 e che appartiene al tipo dei sistemi complessi indicati, quando trovasi caricato per tutta la sua lunghezza da un treno merci dei più pesanti, dà una saetta nel mezzo che non supera i m. 0,25.

---

[25] La seconda edizione ampliata e corretta di queste lezioni venne pubblicata nel 1872 dalla Litografia Ronchi, Milano.
[26] *Il Metodo dell'Area Momento nella determinazione delle condizioni di resistenza delle travi elastiche*, R. Istituto Lombardo, adunanza 22 Gennajo 1881, Milano, Tip. Bernardoni.
   *Sulla determinazione dei Momenti Massimi dovuti ai pesi vincolati sopra una trave appoggiata*. R. Istituto Lombardo, adunanza 10 Marzo 1881. Milano, Tip. Bernardoni.
[27] Il Politecnico; fascicoli di Giugno, Luglio, Ottobre, Novembre 1866.
[28] RITTER. - «*Elem. Theorie und Berechnung eiserner Dach und Brücken Construc.*», Hannover, 1863.
   LENTZ. - «*Die Balkenbrücken von Schmiedeeisen*», Berlin, 1865.
   CULLMANN, - «*Die Graphische Statik*» Zurigo, 1866.
   GOUDARD, - «*Etude comparative de divers systèmes des ponts en fer*» Paris, 1865.
[29] C. CLERICETTI. - «*Sopra i moderni ponti americani e sulle fondazioni tubolari*» Atti del Collegio degli Ingegneri ed Architetti in Milano. Anno V, fasc. II, 1873.
   - *Teoria dei sistemi composti in generale e specialmente dei moderni ponti sospesi americani. Atti del R. Istituto Lombardo di Scienze e Lettere*, seduta del 12 Aprile 1877.
   - *Teoria dei sistemi composti, ecc. Influenza dei carichi accidentali, Atti c.s.*, seduta del 27 Giugno 1878.
   - «*I ponti sospesi rigidi*» Rendiconti c.s. seduta del 3 Aprile 1879.
[30] *Les Travaux publics aux Ètats-Unis d'Amerique en 1870*. Paris, 1875.

Se al calcolo di questo ponte si applica la teoria del Prof. Clericetti, la saetta risulta di m. 0,213. Questi studi sulle travature reticolari complesse vennero apprezzati anche all'estero e segnatamente a Londra ed a New-Jork[31] dove vennero riprodotti in alcuni giornali tecnici con parole lusinghiere per l'autore.

Anche nello studio della stabilità delle volte le ricerche del Prof. Clericetti hanno recato un contingente di nozioni teoriche e sperimentali, le quali, se non raggiunsero lo scopo di rendere completamente determinato l'importante problema, fanno però fede della profonda coltura tecnica che l'autore aveva di questo argomento e dell'acume con cui Egli passava alla disanima dei numerosi fatti sperimentali da Lui raccolti e sapeva dedurne molte utili conseguenze. Il problema teorico sull'equilibrio delle volte, che venne studiato per la prima volta dal La Hire[32] nel 1772, interessò in seguito molti dei più potenti ingegni di matematici e fisici di questo secolo, i quali dettarono sopra questa questione teorie e critiche di grandissimo valore. Bisogna però confessare che a tuttora l'arduo quesito rimane insoluto e nessuna delle ingegnose e svariatissime teorie che si conoscono, può accettarsi come l'esatta sintesi di tutti i fenomeni di deformazione e di resistenze che si manifestano allorché una volta viene abbandonata a sè stessa in seguito al disarmo e si assetta nello stato finale di equilibrio. Il Prof. Clericetti volle affrontare il grave problema; Egli accettò il principio della cerniera del Dupuit[33], però come un fatto transitorio durante il disarmo, e nell'ipotesi che la curva delle pressioni esca dal terzo medio alla chiave ed in un altro giunto inferiore, il che può ritenersi ammissibile nella maggior parte dei casi, conclude: che l'equilibrio definitivo della volta non può avverarsi se non quando le pressioni massime unitarie sono uguali in tre punti.

La dimostrazione colla quale il Prof Clericetti sanzionava il principio da Lui stabilito dei punti di eguale pressione non è rigorosamente matematica né parte da una ben determinata legge fisica, ma risulta da un complesso di deduzioni con cui si interpretano alcuni fatti sperimentali, deduzioni che quasi sempre convincono ma che in alcune parti, troppo arrischiate, non possono distruggere completamente tutti i dubbi che si presentano a chi analizza la teoria del Prof. Clericetti, la quale in ogni modo, rimarrà sempre un documento interessantissimo per gli studi sulle volte e tale da reggere per importanza il confronto colle ricerche sullo stesso argomento che vennero fatte dal Clapeyron, dallo Scheffler, dal Drouets, dal Durand-Claye dal Debauve e da altri autori riputatissimi che prima di Lui si illustrarono con tali studi: e ben fece l'Istituto Lombardo ad assegnare al Clericetti nel 1875 il premio ordinario di L. 1200, incoraggiando così i nobili conati del distinto professore, che con tanto amore cercava collo studio di alcune delle più scabre questioni tecniche di essere sempre più utile alla sua scuola, onorando in pari tempo sè ed il proprio paese.

Tanta intensità di applicazione non arrivava a fare stanca la mente di Celeste Clericetti, Egli, inclinato per natura all'arte, con uno studio profondo, cercava e trovava in questo campo sconfinato del bello riposo e conforto alle sgradite lotte della vita e nuova lena per riprendere il faticoso cammino della sua carriera. Non aveva ancora raggiunta la meta dei suoi studi universitari che Egli, erudito e dotato di un fine sentimento indagatore, percorreva solitario le nostre terre lombarde per studiare, anche nei più dimenticati angoli del contado, i monumenti ed i ruderi che ricordano, della nostra medioevale civiltà, la cultura, le credenze, gli intendimenti artistici, le gioie ed i dolori del popolo.

Il Prof Clericetti sino dal 1862 si fece conoscere tra i più cospicui ricercatori dell'Arte colla memoria sulle ricerche dell'architettura religiosa in Lombardia dal V° al XI° secolo[34], alla quale fecero seguito altre pubblicazioni sull'architettura Lombarda[35] che furono accolte con gran favore dagli studiosi. Questi lavori, i quali hanno meritatamente posto il nome del Prof. Clericetti tra i primi di quanti distinti cultori di questi studi onorino l'Italia, procacciarono a Lui nuove ed importanti occupazioni, dalle quali Egli non mai rifuggiva, ogni qual volta l'opera sua disinteressata veniva richiesta a pubblico beneficio. Come membro, e in seguito vice presidente della Commissione Conservatrice dei Monumenti ed Antichità della Provincia di Milano, nel breve periodo di 4 anni decorsi dal 1882 al 1886 il Prof Clericetti seppe spiegare una, tale attività e competenza nel disbrigo delle non facili mansioni del suo Ufficio da meritarsi il plauso dei colleghi. A questo riguardo. ecco le splendide parole che il Prof. Tito Vignoli, attuale vice presidente della

---

[31] *Abstract of Papers. Institut of Civil Engineers, London. - Van Nostrandis Enqineering Magazin*, New-York, 1881.
[32] *Histoire de l'Académie des Sciences* - Paris, 1772.
[33] *Traité sur l'équilibre des voûtes* - Paris, 1870.
[34] *Politecnico*, 1862.
[35] C. CLERICETTI - «*Ricerche sull'Architettura Lombarda*». Milano, coi tipi della *Perseveranza*, 1869.

Commissione Conservatrice dei Monumenti pronunciava davanti al feretro del compianto professore[36]: «L'acume, la sana critica, la chiarezza, che è l'eleganza del vero, rendevano i tuoi lavori non solo proficui all'arte e al paese, ma erano splendide norme ai compagni, che ascoltavano le tue Relazioni con animo reverente, e con quasi estetico sentimento di un'opera d'arte. Sono più di sedici le tue Relazioni che illustrano l'Archivio della Commissione, e tutte di molto rilievo, risguardanti monumenti della città di Milano e della Provincia[37]. Nelle quali non si potrebbe dire se sia più pregievole il fine gusto, e la sicura intuizione dell'arte, o la dottrina storica, vasta e profonda, che rintegra e fa risorgere sovente da ruderi un monumento, che il tempo, la barbarie, e l'infausta peste dei restauri corrosero, o imbastardirono. Né si dimentichi che tanta parte Egli ebbe nella resurrezione assennata dell'antica Basilica di S. Vincenzo in Prato, di cui si occupò in un dotto lavoro sull'architettura Lombarda sino dal 1869».

Ad altri uffici faticosi e scabri e ad altre iniziative escogitate dalla sua mente, allo scopo di rendersi sempre più utile agli uomini, Celeste Clericetti legò il suo nome che dobbiamo ricordare con reverenza ed amore.

Autore del primo progetto tecnico per la cremazione dei cadaveri[38] e strenuo apostolo di questa riforma igienica che ancora molti contrastano, Egli visse però abbastanza per vedere coronata di successo la sua iniziativa e forse non molto lontano il tempo in cui tutti si potranno convincere che il lampo delle fiamme purificatrici non è di offesa a sentimento alcuno e ben vale il verme che lentamente decompone e trasforma la materia nelle viscere della terra.

Celeste Clericetti amatissimo di Milano, l'opera sua prestava a vantaggio di molte pubbliche istituzioni, e là dove impegnava il suo nome eravi con certezza un largo tributo di attività, di sapere, non

---

[36] Vedi la nota (1).

[37] Comprendendo le relazioni che si riferiscono a lavori secondari, il numero totale di queste ammonta a 40 delle quali le più importanti sono qui ricordate in ordine cronologico:

1882    *S. Lorenzo di Lodi. - Sulle opere in corso e sui ristauri da farsi.*
"    *Portoni di Porta Nuova. Milano. Sui ristauri da farsi.*
"    *Bernate Ticino. - Sulla lunetta marmorea di quella Chiesa che si voleva vendere e proposte relative*
"    *S. Francesco di Lodi. - Diverse relazioni sui ristauri da farsi.*
"    *Duomo di Lodi. - Sulle scale proposte per salire al Presbiterio, e proposte di rilievo e studi sulla vecchia Cripta.*
"    *Castello di Trezzo. - Sulla sua importanza, sul suo stato presente, e proposte relative.*
1883    *S. Marco di Milano - Sulla proposta aggiunta del Pinacolo al Campanile*; due relazioni.
"    *Oratorio di Varedo presso Monza. - Sullo stato ed importanza sua.*
"    *Castello di Trezzo. Schema di regolamento per la sua custodia.*
"    *Torre di S. Giovanni in Conca Milano (in concorso dell'Arch. Archinti e del Consigliere Zerbi). Proposte al Municipio.*
1884    *S. Vincenzo in Prato in Milano (in concorso del Prof. Sacchi) sulla sua importanza e stato suo.*
"    *S Stefano in Vimercate. - Sullo stato del Pronao e proposte sul modo di armarlo.*
"    *Chiesetta ottagonale nel centro del Lazzaretto di Milano. - Sulla proposta sua riforma e dell'aggiunta progettata.*
"    *Torre di S Giovanni in Conca. - Relazione con schizzi e proposte al Municipio di completarne la demolizione.*
"    *Oratorio di S. Eusebio in Cinisello. - Sulla sua importanza.*
1885    *Battistero di Agliate (in concorso dell'Ing. Carcano del Genio Civile) sull'isolamento di quel Battistero.*
"    *Torre di S. Giovanni in Conca. - Sui frammenti di lapidi, iscrizioni, avanzi antichi, ecc. rinvenuti nella demolizione definitiva.*
"    *Villanova Sillaco. – Relazione su quella chiesa campestre e sugli stalli di quel Coro, opera del 16° secolo.*
"    *Palazzo dell'Arrengo in Monza. - Relazione e collaudo del pilone ricostruito.*
1886    *S. Bernardino alle Monache. - Sull'importanza artistica dei pochi dipinti che vi rimangono.*
"    *Castello di Milano. - Proposte e studi per la sua conservazione e condizioni di cessione al Comune.*

Vi sono altre pubblicazioni di interesse artistico del Prof. Clericetti che portano i seguenti titoli:
*Sopra una lapide Longobarda a Beolco (Brianza).* Rivista archeologica di Como, 1876.
*Il Ponte Acquedotto di Spoleto, detto Ponte delle Torri.* Atti del Collegio degli Ingegneri ed Architetti, 1883.

[38] Il Prof, Clericetti nel 1874 costruì per incarico degli Eredi di Alberto Keller la prima Urna ed il primo tempio Crematorio nel Cimitero Monumentale di Milano. Nel 1876 il Prof. Clericetti ed il Dottor Giuseppe Polli ottennero la medaglia d'oro all'Esposizione di Igiene a Bruxelles per l'apparecchio Crematorio eretto a Milano. Nel 1878 il Prof. Clericetti venne nominato Membro Corrispondente straniero della «Société Royale de Médicine publique de Belgique a Bruxelles» in benemerenza della sua ardita riforma igienica.

mai disgiunto dal vivo desiderio di donare generosamente tutto quanto potesse produrre il suo ingegno e la sua assidua applicazione. Molte classi di cittadini lo ebbero nei propri sodalizi come maestro e consigliere valente; sempre tra i primi nel caricarsi di lavoro, molto di sovente era dimenticato o rimaneva tra gli ultimi nel ricevere qualsiasi morale compenso. Sino dal 1868 Egli aveva studiato la questione delle case economiche per gli operai; fu uno dei fondatori della nostra scuola dei Capomastri, la quale ha ricevuto un così giusto indirizzo da riuscire una delle migliori istituzioni scolastiche professionali; contribuì efficacemente a mantenere in fiore la Scuola di Incoraggiamento di Arti e Mestieri, a creare la Scuola degli Orefici, ed anche per la Società cooperativa degli Impiegati trovava tempo e modo di dedicarsi con assiduità nell'interesse di questa filantropica istituzione. Fu giurato operoso, critico valente e relatore all'esposizione di Milano nel 1881 e quale presidente del Collegio degli Ingegneri ebbe la direzione della «Milano Tecnica», nella quale venne inserita una sua pregievole memoria sulla «archeologia milanese». Per altri incarichi né lievi né facili lavorò Celeste Clericettl nell'interesse del Comune, della Provincia e dello Stato, non mai scostandosi da quei sentimenti di onestà e di delicatezza che irradiarono da tutti gli atti della sua vita. Molto alto nelle aspirazioni, molto severo con sè stesso, esatto in tutto, dignitoso allo scrupolo, Celeste Clericetti, solo collo studio indefesso volle salire ad uno ad uno i difficili gradini che lo hanno portato alla sua distinta posizione sociale, alla quale arrivò, ma pur troppo, quasi sfinito di forze e colla mente esausta!

Il povero Clericetti nel giorno 23 Luglio 1866 davanti alla tomba del collega amatissimo Archimede Sacchi[39] pronunciava queste parole «È triste cosa il vedere sparire tante elette intelligenze anche nell'autunno della loro vita; ma è ben più doloroso e straziante il vederci strappati d'intorno, i migliori di noi, nei giorni più attivi dell'esistenza e quando dopo un brillante passato l'attività fisica e la mentale, raggiunto il punto più elevato, tendono a realizzare le più felici previsioni». Queste parole, che a guisa di angoscioso lamento uscivano dal suo cuore addolorato, furono per Lui di ben triste presago! Dieci mesi dopo il povero Clericetti veniva per sempre tolto agli affetti della sua adorata famiglia e gli amici e colleghi di Lui, commossi, gettavano fiori sulla sua bara!

Celeste Clericetti fu una eletta intelligenza, un uomo di cuore, un robusto lavoratore, la cui potenzialità tutta si diffuse modestamente e placidamente ammaestrando nel culto del vero e del bene e creando per noi e pei giovani che dobbiamo educare ed istruire, un raro esempio di quelle virtù che sono base, scudo di ogni civile progresso e che molto si avvicinano all'ideale di quanto dagli uomini può sperare e desiderare la famiglia e la patria!

---

[39] *Politecnico*, 1886. *Discorsi pronunciati in onore alla memoria del Prof. Archimede Sacchi.*

## Il carteggio

**1.**

Roma 30 Aprile 1871.

Caro amico

Ti mando un saluto dall'eterna città, dolente di non averti in compagnia a visitarne le maestose rovine. È un tal cumulo di memorie storiche, di monumenti e di meraviglie artistiche che ci vorrebbe qualche mese a visitarla parte a parte.
Ma quello che sorprende più di tutto è il vedere gli avanzi di due mondi, le reliquie di due civiltà sovrapposte le une alle altre: il mondo romano e il mondo cristiano primitivo, delle Catacombe.
Nulla di più imponente, da una parte, delle colossali rovine delle Terme di Caracalla, del Tempio della Pace, del Coliseo e del Foro Romano, ma dall'altra nulla di più commovente d'una visita alle catacombe fra quelle migliaia di tombe modeste con quelle semplici iscrizioni "Hic requiescit in pace, bona memoria...".
Ieri mattino abbiamo visitato i lavori per la Camera dei Deputati al Monte Citorio[40] e quelli pel Senato, al Palazzo Madama:[41] sono abbastanza inoltrati e diretti con assai intelligenza i primi dall'Ing$^{re}$ Comotti, i secondi dell'Ing$^{re}$ Gabbet. Si sono fatte delle sperienze al Tevere, alle quali non ho assistito per tema che il sole, piuttosto caldo, non mi avesse a far male al capo. Oggi si parte per Velletri alla visita delle Paludi pontine, dove trovasi già il Direttore,[42] recatosi colà ieri mattina. Si ritorna a Roma domani sera e di qui si partirà per Napoli Mercoledì venturo. Abbiamo visitato le meraviglie del Vaticano, cioè il Museo le gallerie, le Logge di Raffaello, la Cappella Sistina col Giudizio di Michelangiolo. Sono veramente sbalorditi da tanti capi d'arte usciti di sotterra freschi e interi come fatti da ieri: il Sofocle, il Demostene, l'Arianna abbandonata, il Fauno e tante altre statue sono capolavori che tutto il mondo ammira: sono fra le più elette manifestazioni del genio romano. Ma non voglio annoiarti con cacchiere *[sic!]* e addio. La mia salute è abbastanza buona e la mantengo tale con un severo regime di astensioni (vino, caffè etc). Godo in anticipazione la visita notturna al Vesuvio collo spettacolo delle lave ardenti. Che peccato che tu non sia qui! Ho ricevuto lettera da Milano: mia moglie ed i ragazzi continuano bene. Altrettanto spero della tua Signora[43] e di tutta la famiglia. Se ne avrò la opportunità ti manderò un saluto anche da Napoli, ma non lo assicuro. Addio. Ricevi una cordiale stretta di mano

dall'aff$^{mo}$ Amico
Celeste Clericetti

**2.**

Caro Amico

27 luglio 1872.

Andiamo oggi a bagnarci nella fresca Acqua del Ticinetto? Sarò a casa tua alle 4 pomerid$^{e}$ a prenderti: se non puoi lascia un tuo biglietto da visita alla posta, con un semplice "non posso". Intanto una stretta di mano all'inglese.

Aff$^{mo}$ Amico
C. Clericetti

---

[40] Il 27 novembre 1871 venne inaugurato Palazzo Montecitorio e vi fu insediata la Camera dei Deputati. I lavori erano stati seguiti da Paolo Comotto, ingegnere dei lavori pubblici.
[41] Il Senato del Regno si riunì a Palazzo Madama per la prima volta il 28 novembre 1871. L'aula era stata progettata dall'ingegnere Luigi Gabet.
[42] Francesco Brioschi, che fu Direttore del Regio Istituto Tecnico Superiore di Milano - poi Politecnico - dalla sua fondazione nel 1863 fino alla morte, avvenuta nel 1897. Suo successore sarà l'allievo Giuseppe Colombo.
[43] Elisa Ferrari.

**3.**

Tremezzo 8 Agosto 1872.

Caro Amico

Ti mando un saluto da queste amene sponde[44] dove vi dimentichiamo con tanto piacere i pensieri e le cure di Milano, ma non gli amici però, come tu vedi. Sono occupato... a far nulla con un'intensità ed un abbandono degno di miglior causa, come si suol dire. Ma sfido io: Alla vista di questo lago, di questi monti, cullati dal rumore di un torrente vicino che mi canta la Ninna-Nanna, che altro posso fare se non l'indiano? A pensare che v'è tanto spazio, tanta apertura di cielo e che noi ci ammucchiamo a migliaia nelle città! Evviva l'età della pietra! Ma pur troppo debbo tornare presto costì e se Brioschi vorrà proprio che si faccia la gita al Brennero vi andrò, ma sarà tempo rubato alle vacanze e perciò da rimpiangere. La Cecilia[45] sta benino: però il viaggio le dà sempre molta stanchezza e non ha ancora passeggiato. Dei ragazzi poi sono poco contento. Pallidi sparuti e dimagrati non sembrano più i nostri due birbanti. Li conduco a prendere la Doccia qui nel Torrente, in un punto dove l'acqua fa un salto d'un 3 metri. Spero che con questo bel cielo e quest'aria balsamica, rifioriranno presto. Ho detto bel cielo, ma l'espressione è un po' iperbolica e convenzionale. Iersera avemmo un temporale fortissimo con la solita scorta di lampi, tuoni e sbuffi di vento. Vorrei pregare la Provvidenza di tenerla corta questa storia, perché come ha detto il Lamartine, <u>anche</u> il <u>Sublime</u> <u>finisce</u> <u>ad</u> <u>annoiare</u>.[46] Stamane il cielo era sereno, ma ora è già ingrugnato e minaccia una seconda edizione della paura di iersera. Rilevo dai giornali che gli scioperi durano tuttavia e si vanno anzi estendendo costì.[47] È un giochetto assai bello, quello degli operai. Quando vogliono berne un litro in più al giorno, si mettono in sciopero e domandano una lira di più ed un'ora di meno di lavoro. È una bellezza: anche facendone uno solo all'anno, il che è assai poco, fra otto o nove anni, pretenderanno una ventina di lire al giorno e, al più, un quarto d'ora di lavoro giornaliero! I professori straordinarii delle Scuole Tecniche e degli Istituti secondarii, hanno (tenuto conto delle tante deduzioni e trattenute e per averne un'idea vedi "El Milanes in Mar"[48]) circa £ 3.00 al giorno e per giungere a tante hanno [stancato] vent'anni e consumato un capitale. I Muratori, che sono muratori pretendono £ 4.0 al giorno! Ma benone: io consiglio da senno, i professori di mettersi in isciopero od almeno di prendere la <u>Cazzuola</u> ed <u>il fratoccio</u>. Addio: con grandissimo piacere: dirigi <u>Como per Tremezzo</u>

Tanti rispetti alla tua signora ed una cordiale stretta di mano per te

dall'aff.mo Amico
Celeste Clericetti

---

[44] Tremezzo si trova sul lago di Como.
[45] Moglie di Clericetti; i figli sono Emilio e Guido.
[46] Alphonse de Lamartine nel prefazio dell'opera *Les confidences* scrive: "Le sublime lasse, le beau trompe, le pathétique seul est infaillible dans l'art".
[47] Si tratta dell'ondata di scioperi che durò per tutta l'estate del '72. Le manifestazioni interessarono le città del settentrione; gli operai scesero in piazza per ridurre l'orario di lavoro, per aumentare lo stipendio e per le elezioni a suffragio allargato. A proposito si veda ad esempio A. Nascimbene, "Manifestazioni popolari e scioperi a Milano dal 1870 al 1872", *Il Politico*, v. 39, 1974, pp. 622-639.
[48] Forse si riferisce alla commedia *On milanes in mar* di Carlo Righetti (alias Cletto Arrighi), esponente della Scapigliatura milanese.

**4.**

Tremezzo 9 Ottobre 1872.

Caro Amico.

La tua venuta è stata per noi una vera festa e tu sei troppo cortese di ringraziarci per quel pochissimo che abbiamo potuto fare. I tuoi sono complimenti che onorano assai più chi li manda che non chi li riceve. La passeggiata nella Val d'Intelvi in vostra compagnia, impedita da questo tempaccio che imperversa tuttora, m'è rimasta nella strozza come una resca di pesce: speriamo almeno per l'anno venturo. Ieri, la Cecilia ed io, ebbimo un'altra grata sorpresa nel ricevere una visita del Prof. Tardy e della sua signora che si sono fermati a Bellagio alla Serbelloni[49] e vi resteranno per una quindicina di giorni, se pure il brutto tempo non li scaccierà. S'è parlato molto di te. Tardy mi disse che deve trovarsi nel venturo mese con te a Torino per quegli Esami di aspiranti ad una Cattedra, di cui mi hai fatto cenno. Sono profondamente addolorato per le notizie che ricevemmo sulla salute del mio povero amico Zaverio Tagliasacchi. È agli estremi, né s'ha alcuna speranza di salvarlo: sarei corso a Milano, e ne ho un vivo desiderio, se egli potesse ricevermi: ma so che non lo può. Povero Zaverio! La nostra amicizia data da vent'anni in questi medesimi giorni: fummo sempre più che fratelli e lui muore; così presto! Sono così inquieto, così rattristato, che può darsi che ad onta della certezza di non poterlo vedere, oggi stesso faccia una corsa a Milano. Addio caro amico: tanti complimenti alla signora ed una calda stretta di mano

dall'aff$^{mo}$
C. Clericetti

**5.**

24 Maggio 1873.

Carissimo Amico

Non esito a proporti il Prof$^e$ Luigi Pagani per l'insegnamento del disegno di fiori, di Paesaggio etc alla signora Elena. Riunisce tutti gli elementi desiderabili: e le informazioni assunte concordano perfettamente col mio giudizio personale sul Conto del prof Pagani. Troverai anche una persona assai simpatica e di buone maniere, il che è già qualche cosa trattandosi specialmente di insegnare ad una ragazza. Il Prof Pagani da me interrogato se accetterebbe l'incarico, rispose affermativamente: ma siccome ha altri impegni bisogna convenire sui giorni e l'ora dell'insegnamento. Assicuro che il Prof Pagani non sarà molto esigente sul prezzo che tu concerterai del resto con lui. Aspetterò la tua risposta per avvertire il Prof Pagani e mandarlo a casa tua. A domattina pel viaggio di <u>circumvelocipedazione</u>. Colla più cordiale stretta di mano credimi sempre

l'aff$^{mo}$ Amico
Clericetti

*[Sul retro appunti e disegni di Cremona]*

---

[49] Hotel Villa Serbelloni a Bellagio.

**6.**

Lucino[50] 3 7bre 1873.

Sono giunto ieri da Milano ove credevo di trovarti per l'Adunanza del Consiglio Com̃le[51]: ma vi ho inteso con molto dispiacere mio e della Cecilia che s'è ammalato il Vittorio:[52] speriamo vivamente che non sia una grave malattia e desideriamo averne notizie. Saprai che la Giunta fu riconfermata quasi per intero non essendo rimasto fuori fra gli Assessori che il Camperio e fra i Supplenti il Fano ed un altro.
Il Sebregondi non ha però ancora accettato definitivamente. Rilevo dai giornali che ti fu offerto il Segretariato Generale del Municipio d'Agricoltura e Commercio: peccato che non sia quello dell'Istruzione! Hai tu accettato?[53] A suo tempo verrà anche l'altro. Il Guido trovasi finalmente in Collegio fin dal gño[54] 30 scorso. È un forte dispiacere per noi che avendo solo due ragazzi fummo obbligati di allontanarne uno dalla famiglia. La prima impressione del collegio sul Guido sembra sia stata buona perché disse e scrisse che è contento: sperava di trovarvi il tuo Vittorio ma confido che lo vedrà presto. Ho passato l'Agosto andando innanzi e indietro di qui a Milano e viceversa. Non ho mai potuto dare le Classi agli Allievi perché non avevano terminati i progetti e il termine della consegna fu successivamente protratto fino all'ultimo del mese. Ho finito domenica. Domani giov 4 si terranno gli Esami generali: e così un'altra schiera di giovinotti abbandona la nostra tutela per slanciarsi nella vita, confidenti in un brillante avvenire. Qui i giorni passano e si assomigliano. È la calma del romitorio, fra poggi ridenti ed ubertose campagne. E il velocipede? Io non l'ho dimenticato: è qui con me e quasi ogni mattina faccio una corsa solitaria su queste strade campestri che, per altro, hanno brusche salite e rapide discese. Quante volte ti desidero compagno, l'aria va facendosi fresca, il sole meno caldo e la vera stagione del Biciclo s'avanza: quando verrete a trovarci? La Cecilia è stata sofferente e malaticcia per una buona parte del mese: ma ora s'è rimessa bene e fa buone passeggiate. La gita pedestre di qui a Como che ho già fatto parecchie volte, è quanto si può dire simpatica: si sale a San Fermo e si discende al Lago entro un'amenissima valle ombrosa: in mancanza d'un amico, un libro qualunque, un romanzo, mi tiene compagnia. E tu, contami come hai passato il mese costì: anfibio come sei, innamorato del sole e del mare, non ti sarai certo risparmiato. Milano è deserta, specialmente di professori e perciò vi ho visti assai pochi amici e colleghi, ma il Martelli è stato malato e l'ho trovato triste e sparuto: ha sofferto una pericardite che gli è capitata d'improvviso e benché siesi rimesso non ha potuto far parte della Commissione che terrà gli Esami generali.[55]
Addio: in attesa di tue notizie, presenta i miei rispetti alla signora Elisa e a tutta la famiglia. Con una calda stretta di mano credimi sempre

Aff.mo Amico
C. Clericetti

Dirigi "Camerlata per Lucino.

---

[50] In provincia di Como, oggi è il comune di Montano Lucino.
[51] Comunale.
[52] Figlio di Cremona.
[53] Luigi Cremona era stato eletto nel Consiglio comunale il 27 luglio 1873, ma non risultano altre cariche.
[54] Giugno.
[55] Probabilmente si tratta dell'ing. Giuseppe Martelli, docente di Lavori in terra e Disegno di opere stradali ed idrauliche al Politecnico.

## 7.

Carissimo Amico

Lucino 13 Settembre 1873.

Ho ricevuto con vivissimo piacere notizie di te che aspettavo da tanto tempo: le nostre prime lettere si sono incrociate. Godo di sentire dalla tua seconda che non vi sia nulla di grave nella malattia del tuo Vittorio: poveretto! Ha goduto ben poco della sua scampagnata. A dirti il vero, sono lieto che tu abbia rifiutato il posto eminente a cui ti chiamava il Ministro Finali:[56] lieto non per te ma per me che amo averti vicino. Del resto se si trattasse dell'Istruzione pubblica sarei il primo a spingerti ad accettare: ad ogni modo è ormai quistione di tempo. La Cecilia è assai fiacca di salute: vi sono giorni in cui si sente in forze: allora si azzarda a fare qualche passeggiata un po' lunga che sconta poi sempre amaramente. Ma essa ama questi ridenti paesi, questo placido Lucino dove abbiamo già passate molte Vacanze.

Abbiamo un panorama dinanzi, così lieto quando il cielo è sereno, così ridente e così vario che innamora! È curioso che quando mi trovo qui mi par di vivere non in questo ma nello scorso secolo, o più indietro, al tempo dei poeti arcaici e degli splendidi Mecenati Lombardi che tenevano Accademia nelle sontuose sale delle loro ville campestri, alcune delle quali rimangono tuttora, testimonii della prosperità economica di un tempo, ma ora deserte. Penso al Parini, al Metastasio al Frugoni, al Giani. Siamo in mezzo ai vigneti e la vendemmia non è lontana.

Mi sono ormai abituato anche a questa strada col velocipede e passo delle ore deliziose a percorrerla in ogni senso sull'obbediente Biciclo: ma vorrei che tu fossi qui meco. Mio fratello Pietro è fuori con noi da parecchi giorni: ha preso a nolo un velocipede a Como: ma incomincia appena a montarlo. Tu riderai che io ti venga innanzi con degli idilii, costì in Roma dove avrai certo gravi occupazioni, ma insomma che vuoi; "Chi va al mulino s'infarina" e in questi paesi ti assicuro che si diventa calmi e sereni come un bel cielo. La strada di qui a Como passando dallo storico San Fermo è assai piacevole e spero che la faremo assieme. Io m'intasco un Romanzo inglese, prendo il bastone e passo passo mi incammino cullando i miei pensieri fra le innocenti emozioni della Novella e quelle ancor più innocenti che desta il verde dei colli, il sorriso della natura. Ecco che mi casca l'asino ancora come diceva il Poeta in quella celebre lettera al Broglio allora Ministro.[57] Non credere però ch'io me ne stia affatto ozioso: leva l'<u>affatto</u> e lascia il resto. Sono circondato (quando sono al tavolo) da libracci di Meccanica, perché debbo preparare un Corso di Lezioni di Meccanica pratica pei Capomastri da tenere alla Società d'Incoragg.$^{to}$.[58] Ho ricevuto lettera dal Paolo Mantegazza: il suo Casino vicino a Lerici è quasi compiuto e ne è proprio contento: anzi mi ha chiesto il permesso di

---

[56] Nel 1873 Gaspare Finali (Ministro dell'Agricoltura) offrì a Cremona il posto di Segretario generale nel Ministero dell'Agricoltura. Cremona inizialmente accettò l'incarico per poi rinunciarvi.

[57] Alessandro Manzoni, presidente della commissione - voluta nel 1868 dal Ministro Emilio Broglio - che si occupò dell'annosa "questione della lingua". Nella sua relazione finale, Manzoni sosteneva che l'unica lingua adatta a diventare la lingua comune fosse il fiorentino, ma tale relazione non contiene la frase citata.

[58] Dalla *Guida di Milano* del 1882: "Scuola dei Capi-Mastri – Questa scuola venne aperta nell'anno 1872 coll'approvazione del Ministero di agricoltura, industria e commercio. Essa è sussidiata dal Municipio e dalla Provincia di Milano; consta di tre anni di studi; e gli insegnamenti prescritti, che si danno ripartitamente presso il R. Istituto tecnico, la R. Accademia di Belle Arti e la Società d'Incoraggiamento d'arti e mestieri, sono: l'aritmetica teorico-pratica, la geometria, il disegno geometrico, le costruzioni, la meccanica pratica, gli elementi di architettura, e la planimetria ed altimetria. [...] Alla fine del terzo corso si danno gli esami di licenza [...]; in base ai risultati dei medesimi la Giunta municipale conferisce le patenti pel libero esercizio della professione di Capomastro." Clericetti qui insegnava geometria, disegno geometrico e meccanica pratica tenendo i corsi presso la Società d'Incoraggiamento d'arti e mestieri che così è descritta a p. 431 della medesima Guida: "Questa instituzione, fondata nel 1838 con offerte dei negozianti di Milano per promuovere le arti utili in queste provincie, e tutelata e promossa con annuale contributo dalla Camera di Commercio, si propagò poi, con generale soscrizione, a tutte le classi. Nel 1841 essa fondò presso di sé una Commissione di chimica e una di meccanica, e nel 1843 una d'agricoltura e una di commercio, incaricate di studiare le condizioni industriali del paese, propagare le utili cognizioni e illuminare coi loro voti la Società, nel 1854 una di meccanica industriale, nel 1857 quella di geometria, nel 1863 disegno geometrico, nel 1871 disegno di macchine, nel 1872 scienze fisiche, e nel 1881 la scuola di frutticoltura d'istituzione agraria Eredi Ponti residente in Varese." Per questa scuola Clericetti teneva i corsi di geometria descrittiva e disegno geometrico.

dedicarvi l'Almanacco Igienico[59] del venturo Anno... È la migliore consolazione d'un architetto, quella di sapere contento il suo cliente! Io dovrò recarmi a Milano pel San Michele, perché come sai cambio d'alloggio; tu potrai venire qui prima o dopo a tua scelta. Faremo qualche scorsa sulle montagne che dominano il Lago: forse per tali passeggiate è preferibile l'Ottobre, come è anche migliore pel velocipede perché ora sudo tuttavia come un puledro. Addio, caro Amico: riceverò sempre con grandissimo piacere tue notizie. Manda i nostri rispetti alla signora Elisa, alla signora Elena e un bacio alla cara Itala. Con una stretta di mano, vigorosa all'inglese credimi sempre

Aff.mo tuo
C. Clericetti.

**8.**

Milano 17 Ottobre 1873.

Carissimo Amico

Brioschi doveva essere a Milano Martedì: mezz'ora prima di partire da Lucino, un telegramma di mio fratello mi avverte che egli aveva differito la partenza e che aspettassi un altro avviso. Mercoledì ricevo una lettera coll'annunzio che Brioschi arrivava in quel giorno e che il successivo (Giovedì) avrebbe presieduto un'Adunanza del Consiglio d'Amm.e della Banca.[60] Venni dunque a Milano Giovedì mattina, cioè ieri ed ecco che Brioschi era già partito per Varese la sera prima e non tornerà che domani Sabato! Mi recai dall'amico Loria[61] per avere qualche notizia su ciò che poteva aver detto Brioschi nelle poche ora in cui rimase qui. Seppi allora che da quasi 15 gñi[62] le lettere dirette al Direttore non erano state spedite a Vienna perché ogni giorno egli telegrafava che era in atto di partenza, e che al suo arrivo gli furono consegnate a un tratto N <u>40</u> lettere che egli si proponeva di leggere in viaggio verso Varese. Non sapeva nulla dell'Istituto, nulla affatto. Quando gli fu riferita la tua decisione[63] e gli si disse che anch'io sarei probabilmente partito e oltre a ciò che anche il Paladini ed il Saviotti[64] erano stati da te richiesti, non disse altro che queste parole: "Ebbene, se ne vadano pure: ne troveremo degli altri!". Era dunque corrucciato. Il Paladini trovò campo di parlargli un istante, un solo istante: gli parlò di se: Brioschi gli disse con un fare pensoso "Vedo che si tratta della sua Carriera, ed io non vorrei comprometterglielа." Mi sorprese quanto poi mi soggiunse il Loria: che cioè tu avevi richiesto il Cerradini, <u>prima di me</u>, per farlo professore ordinario, e che era affare quasi combinato.[65] Gli dissi che non lo credevo: che piuttosto avrai promesso al Cerradini, pel caso in cui io rimanessi a Milano. E <u>persisto</u> <u>a</u> <u>non</u> <u>credere</u> <u>che</u> <u>tu</u> <u>abbia</u> <u>cercato</u> <u>lui</u> <u>prima</u> <u>di</u> <u>me</u>. Ad ogni modo il Loria non contesta il mio diritto di priorità morale sopra suo cognato e si ripromette per Cerradini il mio posto qui, pel caso in cui io venga a Roma, pregandomi in pari tempo di voler fare in modo che egli mi venga pur anco sostituito nei miei uffici presso la società di Incoragg.to.[66] Questi Ebrei non lascierebbero il tempo di morire, ad un cristiano!! È arrivato il Segretario Giovannini[67] ma non gli [punto] parlato, e neppure lo desiderio. Egli crede che si rimedierà a tutto, per quanto riguarda l'avvenire dell'Istituto. Povera testa! Anzi io aspiro che se l'Istituto avesse a ridarsi ad una ventina di assistenti e ad un Direttore-Segretario nella sua persona, tale

---

[59] *Almanacco Igienico Popolare*: serie di pubblicazioni a carattere divulgativo a cura di Paolo Mantegazza.
[60] La Banca di costruzioni di Milano era stata fondata nel 1871 per favorire la realizzazione di opere pubbliche o di pubblica utilità, quali ferrovie, ponti e canali; Francesco Brioschi era presidente del suo Consiglio di Amministrazione. Chiuderà nel 1876; si veda anche la lettera 17.
[61] Leonardo Loria.
[62] Giorni.
[63] Cremona, su invito del Ministro Antonio Scialoja, si era spostato a Roma con l'incarico di riordinare e dirigere la locale Scuola di Applicazione per gli Ingegneri e stava chiamando a sé alcuni colleghi di Milano.
[64] Per l'anno scolastico 1873-74, Ettore Paladini rimarrà a Milano come assistente alle cattedre di Idraulica e Costruzioni idrauliche e ponti, mentre Carlo Saviotti si sposterà a Roma.
[65] Cesare Ceradini si sposta a Roma per l'anno 1873-74, ma come professore straordinario di Statica delle costruzioni e incaricato di Costruzioni stradali. Diverrà ordinario due anni più tardi.
[66] Si veda la lettera precedente.
[67] Ing. Giuseppe.

insieme formerebbe un vero ideale per lui! Egli comanderebbe a bacchetta, regolando i soldi a norma delle ore di lezione, scontando le assenze, dando un'indennità di <u>una lira</u> al giorno per le gite cosidette scientifiche, facendo chiamare professori gli assistenti e riducendo allo stato di assistenti i professori… Domani finalmente spero di poter vedere la Fenice Europea,[68] e regolerò la mia condotta su quanto sarà per dirmi. L'interesse scientifico mi spinge a Roma, non v'ha dubbio: l'interesse materiale del momento, specialmente le considerazioni di famiglia mi spingono a restare . Sono incertezze che confido a te come a carissimo amico. Sono spaventato dal prezzo degli affitti costì e dalla probabilità di dover condurre mia moglie ad un quarto piano!
Confido però che tu saprai risolvere la quistione dell'onorario in modo che almeno a numeri, io non perda, perché le probabilità sono tutte per la partenza. Spero di ricevere una tua lettera che mi dica ciò che hai fatto in questi giorni e se hai parlato col Ministro.[69] Ti scriverò nuovamente domani: intanto ricevi una cordiale stretta di mano

dall'aff.mo Amico
C. Clericetti

**9.**

Carissimo Amico

Lucino 18 8bre 1873.

finalmente oggi potei parlare con Brioschi: fu un colloquio lungo, calmo, spassionato, che ha però lasciato in me una assai triste impressione. Non potrei ripeterti tutto quanto egli ebbe a dirmi: fatto sta che ebbe a mostrare un egoismo meraviglioso, ed un sentire così alto di se così esagerato da farmi stupire. Egli venne a Milano, da Varese con un piano tutto preparato: questo consiste nel lasciar le cose come sono, nel non promettere nulla a nessuno, lasciando a tutti la libertà di optare: v'era un dispetto nel suo accento malamente velato da cortesi parole. Per quanto spetta a me disse che non era punto offeso meco, che qualunque fosse la mia decisione rimarremmo amici: ma quanto alla mia nomina a Professore Ordinario, essere una cosa attualmente impossibile perché ha già presentato da un mese il prospetto della riconferma annuale: che nel suo pensiero il mio valore intellettuale era sempre progredito di anno in anno, che ero stato molto utile allo sviluppo dell'Istituto, che dovevo certo diventare professore titolare, che non mi aveva mai dimenticato, che aveva anzi sul tavolo a Varese i miei ultimi lavori e specialmente quello sulle Volte[70] da esaminare etc ma che questo non era il tempo per passare a tale proposta e che ad ogni modo <u>sotto una pressione non l'avrebbe giammai</u> fatta. Disse che i Professori titolari nominati per suo mezzo valgono assai più di quelli fatti nominare dagli altri, che lui aveva in testa tutto un Sistema circa le Scuole degli Ingegneri e circa le nomine dei Titolari e che tale sistema non l'avrebbe abbandonato se non colla morte.[71] Aggiunse che ad un sistema di idee egli sacrifica tutto, anche gli amici: Allora gli dissi che doveva considerare come una vera disgrazia per me l'essermi imbattuto con lui, anzi una vera disgrazia per tutti il dover dipendere da lui. Rispose essere la cosa verissima: che si pentiva d'aver fatto nominare titolare il Colombo[72] perché era alieno affatto dagli studii. Insomma circa questo lato della quistione, mi espose tutto un mondo di idee sue, l'una più sublime e antipratica dell'altra, ripetendo sempre che <u>sotto una pressione non avrebbe fatto nulla</u>. Si vedeva l'uomo indispettito: disse che era alquanto offeso con te perché quando

---

[68] Si riferisce ironicamente a Brioschi; si veda la lettera successiva.
[69] Antonio Scialoja.
[70] Potrebbe trattarsi della Memoria: "Il principio della cerniera nelle volte", *Il Politecnico - Giornale dell'ingegnere architetto civile ed industriale*, 1873, v. 5, pp. 482-492; 531-544; 577-590. Nel 1875 l'ampliamento di questo lavoro (presentato con il titolo "Ut tensio sic vis") sarà premiato dal Regio Istituto Lombardo con £ 1200 nell'ambito del concorso ordinario dell'Istituto (Classe di scienze matematiche e naturali); la commissione sarà composta da Brioschi, Casorati, Codazza, Pasi e Tatti. Nel 1876 tale contributo sarà pubblicato con il titolo: "La cerniera ed il principio dei punti di egual pressione nelle volte", *ibidem*, v. 8, pp. 100-114; 167-179; 235-241; 329-339; 470-478; 535-542; 612-626; 711-720; 745-758.
[71] Clericetti sarà nominato ordinario l'anno seguente; si veda la lettera 15.
[72] Giuseppe Colombo, ordinario dal 1865.

eri a Roma e lui pure vi si trovava, non gli confidasti nulla. Disse che era suo pensiero farti prima Direttore della scuola preparatoria in Milano: poi piano, piano, Vice Direttore dell'Istituto: la quale dichiarazione è in perfetta contraddizione col fatto che egli respinse la proposta fatta da Casorati, come sai, che egli dovesse pensare ad incaricarti della sua supplenza all'Istituto. Mi disse e ripeté parecchie volte che Roma non è città adatta a fondarvi una scuola di Ingegneri: disse e ripeté che egli l'avrebbe osteggiata con ogni suo potere e che avrebbe cercato di impedire che venisse aperta come ha impedito che si aprisse quella di Ferrara pure decretata dal Parlamento.[73] Mettiti bene in mente queste parole, perché ebbe a ripetermele parecchie volte: egli cercherà di impedire che la scuola di Roma venga attuata e deve appunto recarsi costì fra pochi giorni. Disse che se anche la scuola potesse avere un cominciamento, non potrebbe mai avere molti scolari e siccome il credito dei professori dipende dal numero della scolaresca, quella di Roma non potrebbe mai aspirare che ad una misera esistenza. Disse che quanto al nominare il Paladini a Professore Straordinario era uno sproposito perché troppo novello. Aggiunse che quanto a me invece la posizione era netta: che poiché m'era stata offerta la nomina a professore titolare, era evidente che potessi decidermi ad accettarla e trasferirmi a Roma: ma che non credeva probabile che la cosa potesse effettuarsi senza dipendere dal Consiglio Superiore dell'Istruzione. Disse che la scuola di Milano aveva credito non solo pel suo nome e per la valentia dei professori, ma altresì perché lui, il Direttore per la sua posizione nella Banca,[74] avviava direttamente i giovani nella carriera tecnica, assumendoli ad impieghi stipendiati: che questo non avrebbe certo potuto fare la Scuola di Roma. Insomma caro amico capisci bene che dopo questo, non mi rimane che di venire a Roma. Aspetto dunque per decidermi che tu mi mandi delle assicurazioni positive sulla mia nomina e sullo stipendio: le aspetto ansiosamente e non vorrei che le minaccie del Brioschi arrivassero a capo da intorbidare le cose. Brioschi è poi deciso, vista la piega delle cose, ad attendere all'Istituto nel venturo anno, un po' più di quanto abbia fatto nel decorso biennio: insomma è indispettito. Ma egli non ha tempo di attendervi e la scuola rimarrà tuttavia in balìa del segretario. Devi prepararti a lottare d'influenza col Brioschi che cercherà indubbiamente di mandare a monte la scuola e le nomine. Gli domandai se sapeva che anche il Beltrami andava a Roma. Rispose asciuttamente di sì e non aggiunse verbo su questo particolare. Ho tempo fino alla fine del mese per decidermi: scrivimi dunque qualche cosa di positivo sulla nomina, sullo stipendio, sulla indennità etc. Mia moglie è pure assai ansiosa di conoscere in quali acque nuotiamo: le passerà ogni cruccio di lasciare Milano quando la cosa sarà veramente decisa. Figurati che domandai a Brioschi se, pure abbandonando il pensiero della nomina immediata, voleva promettermi che l'avrebbe attuata entro l'anno scolastico. Mi disse che al momento, cioè adesso, come stanno le cose, non voleva e non poteva promettere nulla. Se mi verranno in mente altre [circostanze] del [notevolissimo] colloquio avuto col Brioschi te le scriverò successivamente: intanto ti mando le prime impressioni. Intanto debbo occuparmi di me: aspetto dunque una tua lettera, ma torno a raccomandarti di stare in guardia perché Brioschi cercherà di impedire che la scuola venga aperta. Addio: in attesa di una tua lettera ti mando una cordialissima stretta di mano

aff.mo Amico
C. Clericetti.

"Camerlata per Lucino"
NB: Paladini è assai indeciso ma la mia opinione dopo averlo udito parlare a lungo, si è che egli resterà a Milano.
Aggiungo una parola: Brioschi non crede affatto che il Municipio Romano voglia dare una somma per l'impianto della scuola. Affretta la cosa per darle un altro punto d'appoggio

---

[73] Si veda ad esempio: A. Fiocca, L. Pepe, "L'Università e le scuole per gli Ingegneri a Ferrara", Ann. Univ. Ferrara, sez. VII, Sc. Mat., v. XXXII, 1986, pp, 125-166.
[74] Banca di costruzioni; si veda la lettera precedente.

**10.**

Carissimo Amico

Milano 26 Ottobre 1873.

Ho ricevuto la tua lettera che rispondeva a quella che ti inviai il giorno prima del mio colloquio con Brioschi: Speravo vivamente di riceverne una in risposta alla Seconda: desideravo conoscere le impressioni che ponno aver fatto su di te le dichiarazioni di quella nuova Sfinge. Sarei anche molto desideroso di sapere se tu l'hai visto in questi giorni poiché Brioschi trovasi costì e se sia a tua cognizione ch'egli abbia già tentato presso il Ministro di demolire la nuova scuola come ha ripetuto di voler fare. Quanto a me, è una quistione assai complessa e devi trovare assai naturale, come la trovano tutti, la mia indecisione. Sono ben lontano dal voler fare una pressione su di te: tu mi offri £ 2000 per la Direzione del Gabinetto, le quali insieme al soldo di professore titolare formano £ 7000. Ora qui, il mio stipendio complessivo è di £ 7200, e bisogna aggiungere che sulle £ 2200 pagatemi dalla Società d'Incoragg$^{to}$ la Società stessa paga la Tassa di Ricchezza mobile. Ecco il perché io misi la condizione che il mio stipendio costì fosse dalle 7 alle 8000 lire. Ora con ciò io non intendo punto che tu cerchi di aumentare la cifra di £ 2000 stabilita pel gabinetto: anzi non lo voglio affatto perché non desidero di entrare in impegni eccessivi col Governo che pagando sempre male, pretende molto. Poi non mi dici nulla delle trattenute e sulla riduzione di soldo che dovrei sopportare nel primo anno: la mia è ad ogni modo per questo lato che pure è tanto importante, una quistione di Aritmetica: Aggiungi che ho preso quest'anno un appartamento nuovo nel quale ho fatto spese e debbo farne assolutamente altre per poterlo abitare foss'anche per un anno solo. Tutto questo fa sì che se io mi decido a venire adesso a Roma dovrei sopportare nel periodo di un anno una spesa fortissima. Dall'altro lato, dopo che il Brioschi m'ebbe a scoprire così nudamente le sue batterie, dopo ch'ebbe a rivelarmi intero il suo carattere profondamente egoista, io resterò certo a disagio in questo Istituto, tanto più che la persona del Segretario[75] m'è proprio pesante, indigeribile. Insomma io ti confido tutto, anche per scritto, perché spero vivamente che la nostra amicizia vorrà durare. Ma ti è proprio necessario un Professore di ponti <u>per questo prossimo anno scolastico</u>? La scuola di Ponti deve appartenere al 3° anno e benché un'ambizione di scuola di Ingegneria vi fosse già costì, dubito che ci siano Allievi pel 3° Anno. Invece tutta la parte delle costruzioni che comprende la Resistenza dei materiali, i solai, i tetti, i muri, le volte etc devono formare il 1° Corso della scuola di costruzioni ed è a questa, mi pare che potrebbe aspirare il Cerradini.[76] Ma sai che sarebbe veramente curioso che un mio scolaro dovesse essere professore titolare, ed io straordinario! Tu capirai benissimo a cosa tendo con tali riflessioni: il professore di Ponti sarà forse inutile quest'anno, mentre sarà necessario il venturo: con ciò avrei un anno di tempo e non sprecherei la £ 2400 che mi costa l'affitto di questo appartamento e le spese pel medesimo che verranno a superare certo il migliaio di lire. Io ho perso ogni affetto nel nostro Istituto perché non ho più fede nel suo Direttore, il quale m'ha impedito anche di assumere una posizione, or son due anni, che già a quest'ora, per altre combinazioni accadute in questo intervallo di tempo mi avrebbe messo al sicuro da ogni evento: di te ho ben altro concetto e ben altro affetto, perché sei uomo di cuore, circostanza a cui io do moltissimo peso. Scrivimi dunque se credi possibile la combinazione che ti ho esposta: ma in ogni caso, guai se il Brioschi ne trapelasse qualche cosa! Sarei rovinato! Credo che tu debba in questi giorni recarti a Firenze colla Commissione d'Inchiesta. Ti prego di telegrafarmi a Lucino scrivendo "<u>Camerlata espresso per Lucino</u>, indicandomi il giorno in cui sarai a Firenze perché farei colà una corsa per parlarti della combinazione proposta e per sentire da te se la cosa è possibile e se tu credi conveniente di farlo. A Brioschi mi basta rispondere per la fine del mese. Tu mi domandi se il Saviotti può assumere un Corso di ferrovie: non lo credo.[77] Ma anche le ferrovie formano un Corso pel 3° Anno e mi pare improbabile che tu possa impiantare addirittura per intero la Scuola. Aspetto dunque una tua lettera ed un tuo telegramma. La Cecilia è venuta a

---

[75] Ing. Giuseppe Giovannini.

[76] Si veda la lettera 8.

[77] Carlo Saviotti prende servizio come "assistente alle cattedre di Statica grafica [corso tenuto da Cremona] e Meccanica industriale [corso tenuto da Giulio Pitocchi] e incaricato delle lezioni sulle strade ferrate". Si veda: *Annuario della Istruzione Pubblica del Regno d'Italia del 1873-74*, Regia Tipografia, Roma, 1874, p. 142.

Milano per condur fuori per pochi giorni il Guido che ha fatto gli esami: essa ti manda tanti saluti. In attesa, ricevi una caldissima stretta di mano

dall'aff.mo Amico
C. Clericetti

NB: Sono chiamato a Torino, per la metà di Novembre circa, da una lettera del Sindaco Rignon,[78] insieme all'Ing.re Tatti per decidere in una quistione tecnica

**11.**

Milano 6 Decembre 1873.

Carissimo Amico

Eccomi ad inviarti un saluto e insieme le mie più calde insistenze perché tu non abbia a rovinarti la salute con eccessive occupazioni. Ho imparato ad apprezzare il valore pratico altissimo del proverbio "Val più un asino vivo d'un dottore morto" e me lo ripeto di sovente a modo di giaculatoria. Bada che ti torneranno i tuoi mali di capo e allora bisognerà pure ad ogni costo che lavori meno: al dolore fisico aggiungerai in allora un cruccio morale insopportabile. Noi siamo assai dolenti della vostra partenza:[79] la signora Elisa ha già imballato e spedito una gran parte della mobiglia *[sic!]* e dopo domani lascierà *[sic!]* anch'essa Milano; l'Itala sembra contentissima di tutto questo trambusto: è una novità ed è naturale che le vada a genio. Sono anch'io molto occupato, non però come te: ha le lezioni all'Istituto, ed è il meno: poi quelle alla società d'Incoragg.to dove ho incominciato un Nuovo Corso di lezioni di Meccanica pratica per Capomastri. Aggiungi che mi sono capitate varie incombenze, non lucrose certo, ma che esigono tempo e pensieri a disimpegnarle. Ho fatto una corsa a Torino con Tatti per quella faccenda del Tempio Israelitico,[80] ed una gita a Calolzio cogli allievi del 3° Anno per assistere alle prove di resistenza del Ponte sull'Adda. Polli m'ha incaricato di sviluppare in una Tavola, i disegni del metodo ch'egli propone per la Cremazione dei cadaveri, dietro richiesta venuta da Zurigo, dove pare che l'incenerazione verrà riattivata prima che da noi: poi desiderò di avere anche i disegni per l'Edicola Crematoria che andrebbe costruita nel Cimitero.[81] Ti assicuro che non ho perduto il mio tempo, il quale è anche più limitato d'una volta, perché di sera non studio più perché non vi reggo. Per altro non ho dimenticato il velocipede: quando posso, scappo via alle 4 pom. dalla Scuola di disegno e corro al Veloce-Club, di dove, sul Biciclo, percorro in ogni senso la piazza Castello e mi spingo anche sino al Rondò, fuori di Porta-Sempione: il sudore mi sgocciola dai capelli sulla faccia, nonostante l'aria fredda rigida che mi avvolge: poi scappo a casa per le 5.30 a sedermi al pranzo. Vittorio è stato assai gentile e cortese verso di me; mi ha scritto una cara letterina nel mio giorno natalizio ed una alla Cecilia nel suo giorno onomastico. La signora Elisa ci ha permesso di levarlo qualche volta dal Collegio insieme al Guido quando voi altri sarete tutti costì. Il Guido nostro continua bene ed è sempre contento del collegio: non è pervenuto a nostra notizia alcun lamento di sorta sulla sua condotta. Il Prof.e Riva continuerà a dire come ha detto a me, che quel ragazzo finge (!). Il Brioschi non si vede mai: mi basti il dire che dall'epoca in cui avemmo quel lungo colloquio assieme in Ottobre, io non l'ho più visto che un solo istante in strada. So però che egli pretende di occuparsi dell'Istituto e che attende (almeno ebbe a dirlo) a dare un carattere di Stabilità al Corso Preparatorio, dove il Jung ha preso il tuo posto, abbandonando il caro tono

---

[78] Felice Rignon. Per la natura della questione tecnica si veda la lettera successiva.
[79] La famiglia stava raggiungendo Cremona a Roma. Si veda la lettera 8.
[80] Clericetti e Tatti vennero chiamati dalla Giunta municipale di Torino per aggregarsi alla locale Commissione di Ornato al fine di studiare la stabilità della Mole Antonelliana che era in costruzione dal 1863. Il 23 dicembre del 1873 essi presentarono una relazione nella quale giudicavano insufficiente la stabilità della cupola: ne proponevano una parziale demolizione per venire sostituita con una in ferro più bassa. L'architetto Antonelli, che aveva allora 75 anni, fu irremovibile e la cupola rimase pressoché fedele al progetto originale.
[81] Nel 1874 Clericetti, insieme a Giovanni Polli, costruì per incarico degli Eredi di Alberto Keller il primo forno crematorio nel Cimitero Monumentale di Milano; nel 1876 ricevettero la medaglia d'oro all'Esposizione di Igiene a Bruxelles per tale apparecchio.

del Parini.[82] L'Istituto nostro per altro cammina bene: abbiamo un terzo anno assai numeroso, poiché conta più di 60 Allievi e lavorano. Rilevo dai giornali che costì piove: noi abbiamo invece una Superba stagione: un'atmosfera limpida purissima, un sole di primavera durante il giorno ed una splendida luna di notte che fa sognare al lago di Como. L'altro giorno a Calolzio, il panorama dell'Adda e delle montagne di Lecco era qualche cosa di divinamente bello. Se posso fare qualche cosa per te scrivimi che lo farò volentieri: e del resto lo sai. Addio: nella speranza di vederti presto a Milano foss'anco per una sola giornata, ricevi un caldissimo saluto

dall'aff.mo Amico
C. Clericetti

**12.**

Carissimo Amico

Milano 28 Decembre 1873.

Dacché hai lasciato Milano, non ebbi più tue notizie se non dagli altri. Né io ti faccio carico di questo, ma il fatto mi lascia supporre che tu sia troppo gravemente occupato per poter trovare un quarto d'ora perduto: il che potrebbe condurre a gravi inconvenienti per la tua Salute. La signora Elisa e le ragazze saranno a quest'ora già installate nella nuova abitazione e spero che vi si troveranno bene: forse voi altri non vi ricordate più Milano se non perché vi avete lasciato il Vittorio.[83] Ma io penso assai volte a te, alle nostre passeggiate in comune, alle nostre corse sul Velocipede e vorrei che tutto quanto accadde dal Settembre in poi non fosse avvenuto. Ho dovuto recarmi parecchie volte a Torino coll'Ing.re Tatti per quella faccenda del Tempio Israelitico. In ultimo vi rimanemmo quattro giorni per stendere il Rapporto e leggerlo alla Commissione di Ornato municipale, e completare così il nostro mandato.[84] Al mio ritorno, la vigilia di Natale, speravo di ricevere una risposta alla lettera che ti scrissi tempo fa, ma non ho trovato nulla. Allora pensai di scriverti nuovamente, per augurare a te e alla tua cara famiglia la miglior salute ed ogni prosperità pel nuovo Anno. La Cecilia fu obbligata a letto per parecchi giorni in causa di febbre reumatica. Ora sta meglio. Essa pure vi invia i suoi migliori augurii. Ho saputo che il Brioschi ha tentato un colpo di testa contro la tua Scuola in unione al Menabrea e che la somma stabilita per l'impianto non passò se non dietro minaccia di ritirarsi dal Gabinetto da parte dello Scialoja.[85] Si vede che Brioschi se l'è proprio legata al dito e che farà tutto il male possibile alla nuova scuola. Sarà tua cura principale adesso di indurre il Municipio di costì a stanziare senz'altro la rilevante somma promessa per la fondazione della Scuola. Se si arriva questo, la scuola sarà assisa incrollabilmente [sic!], né le opposizioni del Brioschi o di altri potranno più distruggerla. Io non ho parlato ancora col Brioschi dacché sono a Milano. E siccome non ho relazione intima di famiglia con nessuno dei colleghi sento molto la tua perdita. A Torino ho visto il Codazza, il quale ebbe a dirmi che ritornerebbe volentieri ad abitare Milano. Non so precisamente a cosa possa egli mirare, pel presente o se non fosse che un semplice voto. L'Istituto nostro va bene e conta tuttavia 200 allievi: nel novero degli Allievi Industriali del 1° Anno, trovo inscritto anche il Guido Perelli e non capisco perché vi si trovi tuttora, trattandosi di un Elenco stampato solamente alcuni giorni sono. Il Gioda[86] mi domanda parecchie volte tue notizie e non posso accontentarlo: spero però che vorrai scrivermi presto. Se vedessi Rodriguez, salutalo tanto da parte nostra. In attesa di tue notizie e pregandoti di riverire per me la tua Signora, ricevi un'affettuosa stretta di mano

dall'aff.mo Amico
C. Clericetti

---

[82] Giuseppe Jung fu assistente di Cremona al Politecnico dal 1867 al 1873, anno in cui Cremona si trasferì a Roma. Successe allora a Eugenio Bertini (che aveva seguito Cremona) sulla cattedra di Matematica al Liceo Parini.
[83] Si veda la lettera 6.
[84] Si veda la lettera precedente. La relazione e le Appendici corredate da tavole e da disegni sono contenute nel volume *Nuovo tempio israelitico in Torino*, Tipografia C. Favale e Compagnia, Torino, 1874, pp. 59-155.
[85] Sugli Indici dei *Rendiconti del Parlamento Italiano* si legge che Brioschi "Ragiona sul progetto di legge per l'impianto di una scuola d'applicazione per gl'ingegneri in Roma"; Menabrea era membro della Commissione di finanze.
[86] Probabilmente si tratta di Carlo Gioda che nel 1873 era Regio Provveditore agli Studi di Milano.

**13.**

Carissimo Amico

Milano 27 Gennaio 1874.

Dirai che sono poltrone per aver tardato tanto a rispondere alla tua lettera, ma il fatto che sta che *[sic!]* la salute mi serve poco, che le lezioni diurne e serali mi stancano per cui alla sera, quando ho finito il lavoro della giornata, non faccio altro. Vi sono però giornata in cui mi sento bene allora lavoro indefessamente. Ma l'inverno è stato crudo e tali giornate rare. Nella prima tua lettera mi domandavi informazioni sopra Scuole, Istituti ed opifici esteri dove un giovane Ingegnere potesse, sia completare la propria educazione, sia darsi a qualche ramo speciale. A tale domanda non rispondo perché ho saputo che venne anche indirizzata al Loria in più tempo che a me e che egli ti ha già risposto. Del resto la domanda quanto ad Istituti di educazione teorica non la comprendeva, parendomi che se un giovane è già ingegnere, perché ha già percorsa una scuola si applicazione italiana, non deve assolutamente aver bisogno di altro insegnamento teorico
*[illeggibile]*
Infatti da quanto ho saputo i giovani già nostri scolari che subirono l'esame erano proprio dei peggiori. I nostri ragazzi Guido e Vittorio sono in buona salute: siamo però assai dolenti tanto io quanto la Cecilia che alcune questioni insorte fra di loro, ci abbian tolto il piacere di avervi per qualche giornata assieme a Natale come era vivo nostro desiderio. Ed avremmo anche *[illeggibile]* alla cosa nella speranza che avendoli assieme , ogni diverbio sarebbe scomparso, se non fosse che il Vittorio scrisse una lettera gentilissima per altro. Alla Cecilia rifiutandosi di venire a casa nostra.
*[illeggibile]*
Tu *[illeggibile]* dell'allievo Gustavo Padoa, che è appunto uno dei miei scolari *[illeggibile]* Non posso dirne bene. Manca quasi sempre alla scuola o per *[illeggibile]* indisposizione fisica o per poltroneria *[illeggibile]*. Così è certo che verrà respinto agli esami. Ha una certa attitudine al disegno ed aveva incominciato benino, poi mancò a lungo e ora è da un bel pezzo che non lo vedo. Gli altri professori sono pure malcontenti del Padoa
*[illeggibile]*
Per esempio il Prof. Poli accompagnava gli allievi a visitare una fattoria con tutte le dipendenze ed una *[illeggibile]* da riso ed un Mulino: *[illeggibile]* un Ponte di ferro uno di muratura *[illeggibile]*. Di un manufatto importante *[illeggibile]* fuori dall'ordinario era *[illeggibile]*che gli allievi potessero ispezionare . Il Loria li impegnava a visitare stazioni ed opifici ferroviarii per mostrare *[illeggibile]* di una strada e tutto il materiale *[illeggibile]*. Una ferrovia e il Colombo li *[illeggibile]* a visitare officine meccaniche speciali etc.
*[illeggibile]*

**14.**

Milano 7 Agosto 1874.

Carissimo Amico

Eccomi ad inviarti un cordiale saluto ed una feroce stretta di mano all'inglese. Avrei dovuto scriverti assai prima, ma tu mi vorrai perdonare chiamandomi anche poltrone, se vuoi. Spero che la tua Signora sia in buona salute e che vi siate accomodati con soddisfazione nella vostra nuova dimora. Ho terminato l'altro ieri i miei esami al nostro Istituto ed ora incomincio a riposare davvero. Il Brioschi è immerso più che mai nella nuova corrente di idee, di occupazioni e non vorrei che finisse per affogarvisi. Il Banco di costruzioni[87] cammina zoppicando poiché ha perduto il credito degli uomini d'affari: molti e molti si sono rovinati per avere dovuto vendere con gravissima perdita le azioni di quella povera impresa in cui sognavano un Eldorado: ma ciò che cammina ancor peggio, a quanto pare, è l'affare dell'Isola d'Elba che il Brioschi, come sai, ha assunto in proprio, e pel quale ora non si trovano i Capitali.[88] È probabile che l'affare dovesse essere abbandonato colla perdita del deposito assai rilevante. E l'Istituto? L'Istituto fa da se come l'Italia del 48; ma se vi fosse un po' d'energia nei professori, si potrebbe richiamare al dovere chi ha sempre preteso d'insegnarlo agli altri. Io mi confido pienamente in te e come non lo farei con alcun altro, specialmente in iscritto, perché conosco la squisita delicatezza del tuo carattere. Il Brioschi pensa all'Istituto e alla sorte dei suoi professori com'io penso alla China: noi siamo senza appoggio alcuno presso il Governo. Io credo che sia ormai tempo, e tutti me lo dicono, di poter essere nominato in pianta stabile: ma quando penso che non ho appoggio alcuno per la negligenza del Direttore, né alcuna voce amica che s'interessi a me che ho sempre fatto il mio dovere e giovato agli interessi della Scuola, assai più di altri che vennero nominati professori ordinarii appena usciti dai banchi della scuola; quando penso che molti miei scolari, addetti agli Istituti tecnici sono in posizione stabile e molti in condizioni finanziarie assai migliori delle mie, quando penso a tutto questo e a tante altre cose, mi viene proprio la bile. Se il Brioschi si rovinasse colle sue nuove imprese, com'è probabile, se lo meriterebbe pel suo egoismo. Io mi congratulo vivamente teco per la nuova onorificenza che ti venne conferita a Berlino, della quale hanno discorso con compiacenza tutti i nostri giornali.[89] Mia moglie coi ragazzi si trovano ora a Cornigliano ligure ai bagni di mare: domattina anzi li raggiungo per ricondurli a Milano. Il Guido ha superato gli esami del corso ed ora desidero che subisca gli esami di concorso ad un posto semi-gratuito nel Collegio Longone,[90] perché la spesa annua è assai rilevante e grave per me. Ma i concorrenti sono almeno il triplo del numero dei posti e però ho ben poca fiducia che si voglia favorire me. Noi andremo a Lucino anche quest'anno a passare la vacanza ma credo per l'ultima volta perché quella casa è ormai passata in altre mani. In settembre spero che avremo fuori con noi per un po' di tempo il tuo Vittorio: così potessi avere la tua compagnia. Desidero ovviamente tue notizie e per mia parte prometto che durante le vacanze ti scriverò parecchie volte. Rammentati *[sic!]* alla tua degna Signora, alla Signora Elena e alla cara Itala e ricevi un cordiale saluto

dell'aff.mo Amico
C. Clericetti

---

[87] Si veda la lettera 8.
[88] Dal verbale della tornata del 24 febbraio del 1873 alla Camera dei Deputati: "È approvata la convenzione stipulata nel dì 20 marzo 1873 tra le finanze dello Stato, l'amministrazione cointeressata delle regie miniere e fonderie del ferro in Toscana ed il signor commendatore Francesco Brioschi per l'accollo a quest'ultimo della escavazione delle miniere Terranera e Calamita nell'isola d'Elba e la vendita del materiale escavato." Tale progetto, pur essendo stato ratificato dal Parlamento con Legge 3 giugno 1874 n. 2083 non poté essere attuato in quanto - come anche riportato da Clericetti in questa lettera - Brioschi non riuscì a trovare i capitali necessari per metterlo in atto.
[89] Cremona ricevette il premio Steiner dell'Accademia di Berlino per i suoi lavori di Geometria pura.
[90] Convitto nazionale Pietro Longone, istituzione ad oggi ancora attiva a Milano.

**15.**

Lucino 9 Ottobre 1874.

Carissimo Amico

Di ritorno da una corsa montana, trovai oggi il tuo graditissimo Telegramma. Accetta i più caldi ringraziamenti miei e della Cecilia poiché sa bene che la mia nomina è dovuta principalmente all'affettuosa premura che ti sei dato per raggiungere l'intento. Sarei anzi curioso di sapere quale condotta abbia tenuto il Brioschi in questa occasione e non sarei meravigliato che avesse, non contrariato, perché in questa occasione la cosa sarebbe stata mostruosa, ma per lo meno cercato di protrarla ancora.[91] Sento che sono obbligato di scrivergli per ringraziarlo, ma è un atto che mi pesa e non saprò trovare le parole: ad ogni modo aspetterò la sua comunicazione ufficiale. Non credevo che la cosa dovesse essere trattata così prestamente ma quando ho letto nei Giornali l'annuncio della riunione del Consiglio Superiore, pensai che fra le altre, avrebbero pure discussa la quistione della mia nomina ed ero altresì certo che tu m'avresti inviato un Telegramma per annunciarmi l'esito poiché conosco per prova il tuo buon cuore. Non so se avrai ricevuto da Milano 3 fotografie comprendenti il Progetto della Cupola in ferro pel Tempio Israelitico di Torino: il progetto è fatto, ma non si potrà passare alla costruzione se la Comunità ebraica non trova i denari per la spesa. E vuolsi una somma forte poiché l'edificio fu piantato a mezzo dall'Antonelli e con una cupola che non può durare lungo tempo.[92] Sono stato a visitare i lavori della ferrovia da Chiasso a Lugano, da Bellinzona per la linea del Gottardo. Se sapessi in che modo sono tratti [sic!] dalla Società tedesca, i poveri imprenditori italiani. La maggior parte sono antichi Cottimisti, cioè Sub-Appaltatori che nella speranza di un lento guadagno, si assunsero in proprio la costruzione di un tronco della linea ed hanno firmato contratti terribili senza conoscerne la portata. Ed ora vi perdono, per le angherie, i soprusi, e la malafede della Società, tutto quanto hanno prima guadagnato in 10, 15 anni di lavoro. Anche qualche ingegnere novello, è stato preso al laccio di quei contratti leonini ed ora si dibatte invano. Uno anzi, l'Ing.$^{re}$ Tassi, mio povero Amico, s'è suicidato questa primavera per disperazione. Intanto lavorano i Tribunali Svizzeri ed italiani con una furia di cause. Non so come il Brioschi abbia accolta la nomina del Bonghi a Ministro, ma son d'avviso che gli suoni amara: le azioni del suo Banco sono così in ribasso che non può ricavarne motivo alcuno di consolazione e tanto meno di gloria, lui che l'ama tanto.[93] Ho visto Rodriguez a Tremezzo dalla signora Kramer: saprai che prende moglie sposando in Dicembre una signora di Amburgo giovane e bella ed altresì agiata. Auguro di cuore al mio Amico una vita felice... A proposito della signora Kramer; non ricordo come, ebbi a parlarle di te. Essa soggiunse "Cremona non mi ama: lo so di certo". Mi credei in dovere di farle per tuo conto una dichiarazione contraria ed allora soggiunse "Almeno so che una volta Cremona non mi voleva bene". Devi dunque averne fatta una grossa.
Addio: rinnovandoti i più vivi e sinceri ringraziamenti, ti prego di riverirci tanto la tua degna Signora. Con una calda stretta di mano credimi

aff.$^{mo}$ ed oblig$^{mo}$ Amico
C. Clericetti

---

[91] Si veda la lettera 9.
[92] Si vedano le lettere 11 e 12.
[93] Si veda la lettera 8.

**16.**

Carissimo Amico.
                                                                                                        Milano 19 Nov. 74.

L'ingranaggio delle Scuole è nuovamente attivato su tutta la linea: m'intendo all'Istituto e alla Società d'Incoragg.to. L'Istituto ha quest'anno meno scolari del solito il che mostra che il tuo assorbe molto, e sarà presto il più numeroso, se pure non lo è già. Sono già tre settimane che sono a Milano con la famiglia: il soggiorno di Lucino m'ha annoiato quest'anno e sentivo proprio il bisogno di rientrare nella Società. Poi già in campagna finisco a fare quasi nulla: lo spettacolo d'una natura così tranquilla, così uniforme mi rende poltrone e finisco a vegetare come una pianta. Sono certo che a Portici voi altri vi sarete divertiti assai più di noi: avevate almeno lo spettacolo del mare, sempre grandioso ed imponente: e che non annoia mai. La mia Cattedra è messa fuori a Concorso in forza di una Legge del 68. Brioschi però mi disse che trattasi d'una pura formalità ed altrettanto ebbe a ripetermi il Bonghi, con cui ebbi il piacere di trovarmi una sera, durante la sua permanenza qui. Sembra che il Ministro abbia trovato un accomodamento per accontentare i pavesi.[94] Nella riunione dei professori e di varii Consigli Amministratori del Ghislieri[95] etc che tenne a Pavia, il Bonghi dopo aver detto che le proposte dell'Università, non gli sembravano abbastanza concrete, li invitò ad occuparsi della istituzione di una Scuola Normale per le Scienze Naturali, lasciando che l'Accademia Milanese svolga il suo compito naturale, di formare una Scuola Normale per le lingue e le lettere: li invitò pure a presentare un Progetto di completamento della facoltà Matematica. Pare che tali proposte siano state accolte favorevolmente: non so poi se si pensi ad attuarle immediatamente. Il Professore Cavallini, nostro collega, come sai, presentò al Direttore una formale domanda per essere nominato professore Ordinario minacciando di dimettersi se la sua domanda non venisse soddisfatta. Ciò ebbe luogo qualche tempo prima che si pensasse a me: il Consiglio Direttivo non credette di presentare la proposta al Consiglio Superiore: si decise di scrivere al Cavallini dimostrandogli l'impossibilità della cosa, stante la secondaria importanza del suo insegnamento e la lucrosa posizione in cui trovasi dall'altra parte il professore attuale: si aggiunse preghiera di voler desistere dal minacciato proposito di dimissione e gli si assegnò un aumento di £ 500 sullo stipendio. Il Prof. Cavallini si è accontentato ed ha fatto bene: è un signore ed ha comperato di recente una Villa a Moltrasio sul Lago di Como (la vecchia casa Passalacqua) lavora molto come ingegnere e non ha perciò bisogno di una posizione stabile migliore, tanto più che non ha famiglia. E come vanno i tuoi mali di capo? Io li soffro tuttavia di frequente, ora che non vado mai al Veloce-Club: m'è mancata la tua compagnia e non ci trovo più il conto. Saprai che Rodriguez, prende in moglie una signora Amburghese: l'ho incontrato dalla signora Kramer queste vacanze. Nel mese venturo Rodriguez andrà a prendere ed a impalmare la sua promessa in Germania: poi si fermerà qualche giorno qui. Ha però bisogno anche lui di essere più franco in salute. La Cecilia sta benino ed altrettanto i due birbanti: spera lo stesso della tua cara famiglia. Tanti rispetti alla signora Elisa ed una calda stretta di mano dal

                                                                                                        tuo aff.mo Amico
                                                                                                        C. Clericetti

---

[94] Probabilmente si riferisce allo scontento dei docenti dell'Università di Pavia per il forte calo di iscrizioni dovuto alla fondazione del Politecnico di Milano.
[95] Regio Collegio Ghislieri.

**17.**

Carissimo Amico.

Milano 20 Marzo 1875.

Colgo con piacere l'occasione della venuta costì del sig Italo Maganzini per inviarti un cordialissimo saluto. Il Maganzini, come sai, è una dei Concorrenti ai posti del Genio Civile ed ha specialmente in vista di ottenere la pensione per l'Estero. È veramente un giovane distinto e certo uno dei migliori Allievi che abbia avuto l'Istituto nostro: molto intelligente ed operosissimo sarebbe sempre in ogni caso un utile acquisto pel Corpo del Genio Civile ed io faccio voti perché egli ottenga quello a cui aspira. Spero che la signora Elisa sia in buona salute e ti prego di presentarle i miei rispetti: la Cecilia sta bene ed altrettanto i due birbanti che mi fanno diventar vecchio maledettamente! Ti ringrazio della Memoria che mi inviasti sopra un lavoro del Ceruti.[96] All'Istituto abbiamo un nuovo professore nel signor Alberto Errera, distinto Economista, il quale dà lezioni di Economia Industriale, specialmente per gli Allievi Meccanici. Il Brioschi ha dato le dimissioni dalla Direzione del Banco-Costruz[i], ma fu pregato di rimanere in carica fin dopo l'Adunanza Generale che dovrà tenersi in Aprile. È assai abbattuto ed a ragione perché il Banco è in uno stato miserrimo e deve liquidare al più presto.[97] Del resto qui nulla di nuovo. Sai tu quando ho ricevuto la mia Nomina a Professore Ordinario? L'altro ieri 18.[98] Oh bacato regno d'Italia... che Dio ti conservi i tuoi linfatici Amministratori che hanno preso per simbolo la lumaca e per motto "Chi va piano va sano". Abbiamo avuto un inverno assai rigido, moltissima neve ed una straordinaria Mortalità nella popolazione: Fino a 49 in una giornata. Contro il tuo probabile desiderio, io faccio sempre voti per il tuo ritorno a Milano. E Rodriguez? non ha ancora preso moglie? Mi dicono che aspetti la calda stagione. Non voglio fare Epigrammi ma ne avrei uno qui sulla punta della lingua. Addio: tanti rispetti alla signora, tanti saluti alla signora Elena e all'Itala. Se non mi scrivi ti strapazzerò. Con un shake hands formidabile e britannico credimi

Affezz.mo Amico
C. Clericetti

---

[96] Potrebbe trattarsi di Valentino Cerruti.
[97] Si veda la lettera 8. In effetti la Banca chiuderà nel 1876.
[98] Si veda la lettera 15.

**18.**

Milano 21 Luglio 1875.

Carissimo Amico

Ricevo ogni tanto notizie indirette di te e della tua famiglia da comuni Amici che vengono di costì e tutti mi dicono che lavori troppo, che sei eccessivamente occupato e che ne soffri in salute. Oh beato quell'anno in cui facevamo le nostre corse mattutine sul velocipede a respirare l'aura pura di bastioni e fuori di Porta e ritornavamo a casa rinvigoriti e pronti al lavoro! Così fossimo ancora assieme... Ho incontrato l'altro ieri mattina il Ceradini che arrivava allora da Roma e doveva recarsi a Cremona: mi parve in eccellente salute. Ho dato ieri l'ultima Lezione all'Istituto: è stata una tirata lunga anche quest'anno: d'estate poi è una doppia fatica il far lezione. In quest'anno però noi non la sentiamo l'estate e non mi ricordo d'aver mai passata una stagione così curiosa: pioggia continua e freddo d'Aprile. Mio fratello Pietro, affranto dal lavoro e rovinato in salute dal rigidissimo inverno scorso, ha dovuto partire per cercare nel riposo, nelle acque di Recoaro, nell'aria di campagna un ristoro alle sue forze ed un rimedio ad una triste malattia che lo minaccia davvicino. Adesso fa la cura del latte in un paesino di montagna sul Lago di Lugano. Ho ricevuto il Rapporto che mi hai scritto sulle Esperienze del Conti sulla Pietra Serena:[99] pare insomma che tutti i suoi studii e le somme spese dal Governo abbiano fruttato ben poco: mi dicono che anche il Richelmi[100] ha stampato una Critica delle esperienze del Conti stesso sull'attrito.[101] Gli esami all'Istituto sono incominciati: darò i miei dal 2 al 6 di Agosto, dopo di che si ha in Progetto una gita piuttosto lunga cogli Allievi del 3° Anno.[102] Dapprima erasi deciso di recarci al Gottardo, ma essendoci stati di fresco quelli di Padova, si temette di recare troppi disturbi ai lavori in corso. Allora il Brioschi espose il progetto di visitare i lavori della ferrovia di Cosenza (appaltati dalla Banca di Costruzione,[103] e ora sospesi per una causa pendente fra essa ed il Governo): di là a Napoli e quindi a Palermo, dove ci troveremmo nei gñi[104] stessi del Congresso Scientifico e dove il Direttore radunerebbe contemporaneamente la Commissione, di cui faccio parte, incaricata di dare gli Esami per la Cattedra di Costruzioni di quella Università. Mi sembra difficile che il progetto possa attuarsi perché gli Allievi che anelano tanto di far gite durante l'anno scolastico, non la pensano ugualmente dopo gli Esami: vi saranno i rimandati che naturalmente si rifiuteranno, poi tutti quelli che non avendo compiuti i lavori alla scuola di disegno, vogliono anzi approfittare della partenza degli altri per terminarli. I Professori stessi non amano stornare una parte del tempo destinato alla campagna: insomma mi par difficile che il viaggio abbia effetto. Ma d'altra parte siccome il Brioschi colla Commiss$^e$ esaminatrice deve pur recarsi a Palermo è probabile che egli voglia proprio effettuarlo. Ho saputo ieri che presto farai una corsa a Milano per recarti in Germania. Io spero e desidero che ciò avvenga quando noi saremo tuttora qui per trovarci un poco assieme. Noi passeremo quest'anno la vacanza a Tremezzo sul Lago di Como. Tanti rispetti alla tua Signora anche da parte della mia Cecilia. Inviandoti una cordialissima stretta di mano, col desiderio vivissimo di salutarti presto personalmente credimi

Aff.$^{mo}$ Amico
C. Clericetti

---

[99] "Relazione intorno alla Memoria del sig. Colonnello Pietro Conti, avente per titolo : "Sulla resistenza alla flessione della pietra serena" della Commissione composta dei soci Betocchi, Cremona, Beltrami (relatore). Letta nella sessione del 7 marzo. 1875", *Atti della R. Accademia dei Lincei*, 1874-1875, v. 2, pp. 408-416.
[100] Si tratta di Prospero Richelmy.
[101] P. Richelmy, "Impressioni prodotte dall'esame della Memoria del colonnello Conti intorno all'attrito", *Atti della Reale Accademia delle scienze di Torino*, 1875, v. 10, pp. 773-805. Richelmy scrisse un altro contributo l'anno successivo sul medesimo argomento "Nuovi appunti alle osservazioni presentate dal sig. colonnello Conti in difesa della sua Memoria sull'attrito", *Atti della Reale Accademia delle scienze di Torino*, 1876, v. 11, pp. 663-673.
[102] "Effemeridi dell'Istituto Tecnico Superiore nell'anno 1874-1875", in: *Programma del Regio Istituto Tecnico Superiore di Milano. Anno 1875-1876*, 1875, pp. 13-16.
[103] Si veda la lettera 8.
[104] Giorni.

**19.**

Milano 20 Luglio 1877.

Carissimo Amico

Sono felice di avere finalmente avuto tue notizie ma in pari tempo ne sono desolato. Figurati! Tu hai lasciato Roma ieri, 19 ed io ci sarò domani sera 21 e contavo tanto sull'improvvisata che volevo farti! Parto stasera per Catania e l'unica mia fermata è appunto Roma, dove il primo obbiettivo era di stringerti la mano e passare qualche ora piacevolmente assieme. Rilevo però dalla tua lettera che la tua Signora e il Vittorio sono tuttavia colà e mi riserbo di far loro una visita.
Ti portavo i saluti del Culmann da Zurigo che mi ha parlato con tanta ammirazione dei tuoi lavori. Vi fui cogli Allievi dell'Istituto, con Stoppani Colombo e Martelli e fummo accolti festosamente dai Professori e dagli Allievi del grande Politecnico. Gli allievi ci hanno dedicato una Kneipe[105] in una Birreria dove passammo la serata in circa 200! Fra i brindisi e le tazze di birra. Ora vorrai sapere cosa diavolo io vada a fare a Catania: ecco. Vi sono invitato dal Sindaco di quella città,[106] per esaminare le condizioni di stabilità di un teatro di recente costruzione e riferirne in proposito a quella Giunta Municipale. Il compito è serio ed ho aspettato ad accettare, ma infine è nella cerchia principale dei miei Studii. Salirò l'Etna, un mese dopo aver salito lo Spluga ed il Gottardo. Dunque alla tua Scuola pare che tutto sia già finito: noi abbiamo appena incominciato gli Esami e parecchi Colleghi daranno altre Lezioni nei giorni che separano gli Esami Speciali da quelli generali. I miei esami speciali li terrò al mio ritorno dalla Sicilia. Costì fralle *[sic!]* vergini bellezze delle Alpi, passerai certo giorni calmi sereni e felici ed io li auguro di cuore a te che hai tanto bisogno di riposo ed alle tue ragazze Elena ed Itala che vorrai salutare per me. La Cecilia ed i ragazzi sono in buona salute e si apprestano a partire per Monticello dove anche quest'anno, passeranno la vacanza e dove verranno raggiunti dal Guido appena avrà compiuti i suoi esami di ammissione al Liceo. Addio: ti auguro la più lieta vacanza durante la quale riceverai certo mie notizie e fors'anche da Catania se appena ne avrò tempo.
Ricevi un cordialissimo saluto

dall'aff.mo Amico
C. Clericetti

---

[105] Dal tedesco, significa *pub*, *birreria*. Probabilmente Clericetti ha usato un termine errato.
[106] Potrebbe trattarsi di Francesco Tenerelli.

**20.**

Monticello 17 7bre 1877.

Carissimo Amico.

Sono lietissimo di aver ricevuto tue nuove e di saperti in eccellente salute. La tua rapida descrizione delle bellezze naturali delle regioni Alpine che hai visitate, mi ha messo l'acquolina alla bocca e il desiderio vivissimo che io fossi stato in tua compagnia. Ho bensì valicato a piedi nello scorso Giugno il S. Gottardo, lo Spluga e il *[Briga]*, ma ho potuto appena intravederne le magiche bellezze perché distratto dalla numerosa compagnia in cui eravamo e dalla rapidità della nostra gita che non ci consentiva riposo, né il tempo bastevole a godere l'incanto di quei romiti valloni, di quelle alpestri cime. Ho veduto come sai la tua Signora a Roma ed era impaziente di raggiungerti per togliersi al caldo infocato della capitale. La gita che ho fatto in Sicilia, fu rapida come un Express-train: vidi e ammirai la gentile Catania dove rimasi cinque giorni a studiare la quistione di quel teatro e quando mi ebbi assunto tutti i dati che m'occorrevano e istituite quelle Esperienze di resistenza che trovai opportuno a formarmi un chiaro concetto delle sue condizioni di stabilità, me ne partii. Potei così togliere molti dubbii insorti ed aquetare *[sic!]* l'animo di parecchi che avevano timore: a suo tempo mi occuperò di stendere la Relazione per inviarla colà. Vi trovai parecchi antichi Allievi del nostro Istituto che m'accolsero tutti assai cortesemente e con vera espansione d'affetto. Mi accompagnarono a visitare le antichità di Siracusa, che per altro sono allo stato rudimentale, come la più gran parte di quelle di Roma: sono ruderi, sostruzioni e bott-lì. Trovai però interessante il fenomeno acustico detto l'orecchio di Dionigi nelle vetuste cave di pietra, dove sudarono e morirono negletti tanti poveri schiavi e prigionieri. Feci per mare anche il tragitto da Messina e Palermo, dove mi fece gli onori di casa l'Architetto Basile che vi sta costruendo un grandioso e magnifico teatro di stile greco, un edificio veramente monumentale.[107] Visitai con lui le poche antichità Arabe e Normanne e lo stupendo Duomo ed il Chiostro di Monreale che sono fra i più preziosi gioielli del Medio-Evo italiano. Ma ciò che ammirai sopra tutto in Sicilia è la incantevole bellezza della costa marina nei dintorni di Palermo e quel tratto da Taormina a Siracusa, fu per me come una rivelazione e mi sentii crescere a mille doppi l'ammirazione per quell'eletto popolo greco e per le antichissime popolazioni italiche, che ebbero un sentimento così potente e così elevato del bello da scegliere a luogo di dimora i più stupendi punti del globo. Visitai finalmente Girgenti e là potei finalmente ammirare un tempio greco in piedi, quello della Concordia, più le rovine grandiose di molti altri. E qui pure, come non ammirare la squisitezza del sentimento estetico di quegli antichi che sapevano approfittare delle bellezze e degli accidenti del terreno, per dar risalto alle loro opere Architettoniche, coll'elevare i loro templi sopra promontori e rupi in cospetto dell'immensità del mare! Com'è stupendo l'effetto delle tinte infocate di quel cielo tropicale e di quell'ardente sole sul colore giallastro delle vetuste colonne scanalate, sui frontoni, sugli architravi di quei colossali edifici! La luce in Sicilia ha qualcosa di così ridente che si direbbe viva: sono particelle animate che colpiscono l'occhio, tanto è il brio del suo splendore e il fascino che esercita! Ma basta della Sicilia. Ora sono a Monticello colla famiglia: cioè ci sono e non ci sono perché ad ogni 3 o 4 giorni mi tocca di recarmi a Milano per svariate faccende. Il Ponte di Lodi sull'Adda ebbe a subire guasti durante l'ultima piena. Fu dunque nominata una Commissione e cui pure appartengo per decidere il da farsi: e ad ogni tratto bisogna che ci riuniamo o a Lodi o a Milano per discutere. Poi si sta costruendo un ponte sull'Oglio (ferrovia Treviglio-Rovato) in un'arcata di muro di $42.^{m}00$ e sono incaricato di fare ogni settimana delle esperienze sulla resistenza dei mattoni da impiegare pel medesimo. Così ad ogni infornata di mattoni che escono dalla fornace di Cremona bisogna che ne sperimenti parecchi e riferisca in proposito. Ma a quest'ora ti avrò stancato di notizie: ne aggiungerò una sola ed è penosa pel nostro Istituto. Il Prof Stoppani, malcontento di Milano, ci abbandona e si trasporta a Firenze a quell'Istituto di Studii Superiori. Tanti rispetti alla tua Signora anche da parte della Cecilia e saluti da noi tutti per la signora Elena e pel Vittorio. Ricevi una calda stretta di mano dal tuo affezionatissimo amico

C. Clericetti

Scrivimi nuovamente: te lo raccomando

---

[107] Giovan Battista Basile, che progettò il Teatro Massimo di Palermo.

**21.**

Milano 13 Giugno 1878.

Caro Amico.

La morte di tuo fratello[108] ha colpito tutti dolorosamente. Artista distinto e originale aveva creato una scuola propria e aveva molti ammiratori. La solita provvidenza lo ha colpito nel più bello della sua carriera e nei giorni migliori della sua vita! In pochi anni è questa la seconda sventura[109] che colpisce la tua famiglia e noi dividiamo profondamente la tua afflizione che il tempo solo può lenire.
Ho saputo che eri a Milano in questi giorni e sono passato all'Albergo del Biscione nella speranza che tu vi fossi alloggiato e che potessi almeno stringerti la mano: ma tu non eri alloggiato colà. Comprendo che in questi giorni dolorosi avrai preferito la solitudine, ma avrei desiderato tanto di vederti e con me mia moglie.
Sarà una lieve consolazione per te, ma sarà pure un conforto nel tuo domestico lutto il vedere come la città intera ha condiviso il tuo dolore. La Cecilia desidera di essere rammentata a te, alla tua signora e a tutta la famiglia. Accogli un cordiale saluto ed una stretta di mano

dell'aff.mo Amico
C. Clericetti

**22.**

Carissimo Amico.

Milano 20 Dicembre 1878.

Mando un cordiale saluto a te e i miei doveri alla tua Signora, insieme ai più vivi ringraziamenti per le squisite e affettuose cortesie che ho ricevute da voi altri durante il mio soggiorno costì. Sono prove di amicizia che non dimenticherò mai perché le ore passate colla tua famiglia furono le migliori per me: colgo pure l'occasione della vicinanza delle feste Natalizie e del Capo d'anno, per inviarvi i migliori augurii, i voti più sinceri pel vostro benessere e la vostra prosperità.
Accludo alla presente una mia sembianza che troverai alquanto antica rispetto a me, o viceversa più giovane di me. Ma io non ci ho colpa se i ritratti non invecchiano e se, per rispettosa deferenza lasciano all'originale questo misero vanto. Giunto a Milano trovai un freddo rigidissimo e tale che passerei volentieri tutto l'inverno a letto: e invece ho dovuto riprendere immediatamente la croce delle mie Lezioni all'Istituto e alla società. Anche il povero Lombardini fu "Assunto ai misteri della seconda vita" come disse il Tecchio in Senato annunciando la morte del Re:[110] aveva però 84 anni e mi pare che a tale età si possa scendere calmi e rassegnati nella tomba. Il Tatti, commosso fino alle lagrime, disse bellissime cose sulla fossa del suo vecchio Amico e anche l'amico Tamanini di Trento[111] è immerso nel lutto per la repentina morte di suo padre che lo lascia solo nella deserta casa: finirà a trovare una buona volta la pianta a cui appiccicarsi che Bertoldo non trovava mai: voglio dire finirà ad ammogliarsi. Essendomi sfuggita tale poetica figura, ti prego di non mostrare questa lettera alla signora Elisa ché quanto alla mia Cecilia non la vedrà di certo.
Oggi giungeranno da Treviso i nostri sposi[112] per passare il Natale con noi: la morte ci va mietendo attorno tanti amici e parenti che bisogna ben stringere i nodi fra i superstiti. Ho pensato molte volte, che quando s'arriva ad una certa età converrebbe cambiare ambiente per non assistere da vicino alla scomparsa di tanti Amici. Volevo essere allegro con te perché è la fine dell'anno, ma appunto precisamente per questo, mi passa dinnanzi come una funebre rassegna di vittime del 78 che mi forzano il pensiero e la mano. Ieri c'è stata la distribuzione dei premi al Ginnasio Beccaria e l'Emilietto ha portato a casa due premi, l'uno per la

---

[108] Tranquillo, morto il 10 giugno.
[109] Potrebbe riferirsi alla morte del fratello Giuseppe (Peppino), avvenuta nel 1876.
[110] Vittorio Emanuele II.
[111] Francesco Saverio, architetto-ingegnere, figlio di Michele, pure ingegnere.
[112] Potrebbe trattarsi del fratello Emilio e della moglie.

Scuola, l'altro disegno. È un ragazzo di molto puntiglio e che ho fiducia sarà per riuscire bene. Del Guido spero bene nel prossimo anno. Ti rinnovo tanti ringraziamenti a te e alla Signora Elisa, anche da parte della Cecilia. Tanti saluti alla signora Elena, al Vittorio, all'Itala per me e tutti i miei: i miei doveri alla tua signora e per te una stretta affettuosissima di mano

<div style="text-align:right">C. Clericetti</div>

### 23.

Carissimo Amico

<div style="text-align:right">Milano 12 Maggio 1879</div>

Sono lietissimo che tu venga fra noi entro il mese e spero ed insisto fin d'ora perché trovi modo di venire una giornata a desinare da noi, e ne faccio conto. Dunque verrai cogli Allievi, perché Nazzani mi scrisse che sulla fine del mese farà la gita d'istruzione venendo nell'Alta Italia. Alla Scuola nostra non s'è ancora pensato alla gita: già la si fa sempre così tardi che v'è tempo. È però probabile che si visitino i lavori della Pontebba[113] e forse il Brennero.

Ti do la notizia che il cielo è sereno e che sussiste la possibilità ed anzi un certo grado di probabilità, che abbia a durar tale per almeno 20 minuti primi! Almeno a tali conclusioni conducono gli Studii metereologici fatti in questi ultimi sette mesi in cui il <u>Dinamismo Endogeno</u> del Padre Cecchi[114] ne ha fatte d'ogni sorta. Un'altra notizia che forse però è già a tua cognizione si è che anche il nostro Collega Jung ha trovato finalmente quella tale <u>famosa pianta</u>[115] sotto le forme di una gentile signorina Cantoni, sorella di un nostro Allievo del 3° Anno.[116] E così la nostra Scuola continua a dare per parte dei professori un alto esempio di moralità e di osservanza scrupolosa del precetto Biblico "Crescete et etc"

Noi siamo tutti in buona salute ad onta della perversa stagione e ne ringraziamo... noi stessi perché abbiamo saputo preservarci dall'influsso del famoso Dinamismo. Tanti doveri alla signora Elisa e a tutta la tua gentile famiglia. Accogli un cordiale saluto ed una stretta affettuosa di mano

<div style="text-align:right">dell'aff.mo Amico<br>C. Clericetti</div>

### 24.

<div style="text-align:right">20 May 1879</div>

My dear Clericetti,

I received your letter with great pleasure. Accept my best thanks for the invitation to dinner you have the kindness to give me, taking it for granted that I should come to Milan with our scholars. But the matter is otherwise. It is settled that only the scholars of our highest class will go on a little journey almost exclusively for the sake of instruction in hydraulics. On this account and further for economical reasons, only two professors (Nazzani of hydraulics and Favero of bridges and railroads) with their Assistants will lead the little company. Ought to give the bad example of increasing the expense by the presence of useless person? Perhaps you think that I as a Senator, can run on railways, free of cost, but also on this point, you would be mistaken. I am yet a Senator <u>in partibus infidelium</u>; the "august Assembly" has not yet vouchsafed to occupy itself about us and throw its doors wide open to give us admittance.

---

[113] Si riferisce alla ferrovia pontebbana, inaugurata il 30/10/79.
[114] Padre Filippo Cecchi, al secolo Giulio Isdegerde. La meteorologia endogena (o dinamismo endogeno) si occupava di studiare un metodo che permettesse di prevedere i movimenti tellurici anche mediante l'osservazione dei fenomeni metereologici.
[115] Si veda la lettera precedente.
[116] Bice Cantoni sposò Giuseppe Jung il 12 agosto; il fratello Vittorio si laureò al Politecnico il 9 settembre dello stesso anno.

Moreover, I cannot suspend my lessons on higher geometry; they must be continued till the middle of June. Our travelers will spend in Milan the 30$^{th}$ and 31$^{th}$ instant; I am sure, it is unnecessary that I should recommend them to you and to my other old colleagues.

You give me the comforting information that the sky above your head is clear, with the probability of its remaining for ten minutes such as it is. On my side, I am sorry to be obliged to tell you that our sky is as cloudy and threatening as it has been these seven months. Last night we had a dreadful music of innumerable deafening claps of thunder which prevented my sleeping in such a manner that at present I am not free from headache.

The progress of meteorological Science is very consoling.[117]

It has been made known to me that also our friend Mr. Jung has found his "tree"; but you are not right at all in slandering this botanical species. Mrs Cecilia is a pleasant and fruitful tree, on which everybody would be willing to be hanged. I wish Mr. Jung a similar fate.

And now, good bye, my dear friend.

Please, inform me about your plans for the holydays. Till when will you remain in Milan and where are you going afterwards? Probably I shall be obliged to remain here till the middle of July; then I will go on the Alps. My family returns your hearting compliments; give my remembrances to Mrs Cecilia and my love to your dear boys.

With the hope of seeing you in Milan, in the month of July,

Your faithful
L. Cremona

**25.**

5$^{th}$ June 1879

My dear Clericetti,

I return you my best thanks for your kind letter of 3$^{rd}$ instant. I am very glad of your proposal to keep hence forward our correspondence in English, I accept it with the greatest pleasure; it will be very good practice for me, though I must spend much more time on it than if I write in another language. I here renew the expression of my gratitude to you and your colleagues and scholars for the hearty reception of our little company. Nazzani wrote me three letters from Milan, informing me of your kindness, and especially he told me of the toasts at the banquet at the Public Gardens. It is a great satisfaction to me to see how the young Roman school has obtained the esteem of its elder sisters of Milan and Turin, and how promptly currents of sympathy spring up among the students of such distant cities. As for me I feel much consoled hearing from you that my old friends and the scholars of the I. T. S. have not forgotten me. But I cannot keep myself from the doubt of which I spoke to you in my last letter. That is to say, I fear, that discretion and modesty have not been as strictly maintained as could have been wished; am I wrong, my dear friend? I should be very satisfied if I'll mistaken. Our young travelers are at present at Spezia; they will arrive in Rome tomorrow evening.

I learn from your letter that the bad weather has not yet left off troubling you; we are less unfortunate. These last eight days all rain has escaped and the sun is splendid, even too hot. However, I fear, it is too late to remedy seven months of inclement weather and the crops will be poor everywhere. How nice if now it were really possible sensibly to diminish the taxes! Misfortune surrounds us on all sides; river floods in the North, volcanic eruptions in the South; poor Italy seems to be threatened by God's anger, like Egypt in Moses's time. If so it be even the prospect of the increase of insects horrifies me. Speaking of volcanoes, I eagerly long for a visit to Etna; and probably I shall at length decide upon running down there.

And you? From association of ideas, my thoughts fly to the Alps; no sooner will my holydays begin than I intend undertaking a series of excursions in the Val d'Aosta and lateral valleys; and you have made some plan, have you not?

I am very glad of the good news of your boys, but I wish Mrs Clericetti would imitate their example.

---

[117] Si veda la lettera precedente.

All the members of my family are in good health and give me the charge of sending their love to your lady and sons.
And now, good bye. I always am, my dear Clericetti

Your affectionate
L. C.

**26.**

15 June 1879

My dear Clericetti,

Your letter has somewhat calmed my apprehensions as regards the passage of our students through Milan. If I say that I already perfectly knew Mr. N…i's[118] character perhaps you will wonder that I trusted him with the care of leading the young Company. But at first the matter was settled otherwise. Professors F.[119] and N. were both designed for accompanying the students, and the former should have been the Commander of the Company. But on the day of departure, F. was not quite well and I said he would join the travelers after one or two days. In this way the matter was sent off some days, it always being thought that F. was leaving on the morrow; and when it was evident that F.'s health did not get better it was too late to give another Prof. the charge of overtaking N. No doubt that N. is a very clever man in his branch of science; but my wish would have been to avoid any parade and any sounding of trumpets. Meanwhile I was informed that another inconvenience had taken place; N. seems to have thought that his dignity of Professor did not allow him to introduce the Engineers C. and S.[120] who, being two assistants of our School, were of the Company. C. & S. are two very good-natured intelligent young men; the former, as an old scholar of Turin, in that city made his own way and on some occasions introduced his friend also; but Sciolette for instance in Milan, was completely overlooked thrown into the shade.  At the dinner of the Public Gardens, Sciolette was rudely reproved by N. because he arrived to the appointment ten minutes after him, but not even then was he introduced to the persons with whom he was to dine; and then S. would have been quite alone. if some Milanese students, taking him for one of the Roman scholars, had not had the kindness to keep him with them: poor S. felt so mortified that he had scarcely re-entered Rome than he declared his wish to give over his office at our School. Very small trifles, my dear friend; but I avail myself of them in order to prolong my English exercise…
Our scholars were greatly satisfied with the hearty reception they met with among their Milanese colleagues.
Prof. Favero has not yet got over his illness; and N. is sorrowful on account of the dangerous condition of his daughters; one of whom has been seized with diphtheria and the other has caught a nasty disease by intercourse with woman servant![121]
We wish Mrs. Clericetti a speedy recovery; please, give her our remembrances. Returning you a British "shake hands", I am, my dear Clericetti,

Your affectionate
LC

---

[118] Nazzani's.
[119] Favero.
[120] Carlo Sciolette, ingegnere architetto; non si identifica C., l'altro assistente
[121] Una delle figlie morirà nei giorni successivi: si veda la lettera seguente e la 28.

**27.**

Milan June the 24<sup>th</sup> 1879

My dear friend

I beg excuse for my having not answered before to your last letter, but I have been very much engaged these last days. Few days ago I received from Nazzani the announcement of the death of his poor little girl from the diphtheria and feel very sorry for him and his wife: it was really a bad conclusion of his gay journey with the Students. Now we are on the eve of starting with our scholars, but in a very different direction. We go to Turin there to remain a couple of days: then we go through the Cenis to Bellegarde to inspect the Works on a railway line in construction and from thence to Geneva. We shall stay two days also in that charming town and then we shall proceed to visit the Hydraulic works of rectification of the Aar and thence to Friburgo and perhaps to Berna so that we shall return by the way of the Simplon crossing the mountain on foot, from Briggs (I believe) to Domodossola. It is rather a long and fatiguing journey, to which I would personally have prefered *[sic!]* a visit to the railway works of the Pontebba-line which had, I must confess, the disadvantage of not presenting any Warehouse or Mill for the instruction of the Mechanical-Scholars of our 3<sup>d</sup> year. You will receive a letter from me either from Ginevra or from Friburgo. We have had nearly a month of splendid weather: our sky is as beautiful and calm and blue as the most charming Lombard sky mentioned by Manzoni and all the sad pronostics[122] *[sic!]*, deduced from the observations of the Stars, and the Historical remembrances of other unhappy years or periods of years, brought to light by Archeologues *[sic!]* have been forgotten. So the world is not going to a smash as it was feared by many, and we shall be deprived of the sight of a splendid Catastroph[123] *[sic!]*. M<sup>r</sup> Tamanini was lately in Milan for few days: he came to see us and left many salutes[124] *[sic!]* for you. As I had foreseen since his father's death, he is at last going to marry a nice young girl from Trento of whom he showed to us the portrait. The marriage will take place I believe next month or in August, and they intend to visit together Rome and Naples. Tamanini is sure to call at your house, because he wishes to present his bride to your family: but very probably at that period you will all be absent from the Eternal City. My brother Emilio the officer is going to leave Treviso with his wife in order to follow his Regiment to Naples, where they shall remain certainly for a couple of years. Every one *[sic!]* is now on the eve of depart,[125] either for bathing places or Sea shores or the mountains or the country, so that in a short time, Milan will be a desert, as it always is during the Autumn. Adieu my dear friend. My kind remembrances to Signora Elisa and all your family in which my wife hearthly *[sic!]* joins me.
Believe me ever

affectionately yours
C. Clericetti

---

[122] Corretto da Cremona in *prognostics*.
[123] Corretto da Cremona in *Catastrophe*.
[124] Corretto da Cremona in *salutations*.
[125] Corretto da Cremona in *departure*.

**28.**

5 July 1879

My dear friend,

You have not occasion to make an apology for your delay in answering me; it is now my turn to thank you for your very interesting letter first of all for the kind wishes which you and Mrs Clericetti sent me on the occasion of my Saint's day.

Nazzani's wife and children left Rome immediately after the death of the eldest girl; but three days later Nazzani himself was obliged to overtake them (at Zibello near Parma) because the second one also had caught the diphteria. They passed some days in the greatest desolation, aggravated by the alarm which the name of the dreadful disease had thrown among the people of that house and country. Now the news is a little better, the danger seems to (have ceased) (to have passed away).

The plan of the trip you are on the eve of undertaking with your scholars is not wanting in interest. It has some parts in common with my route for Alpine excursion. But the time is quite different; in fact, I cannot start from Rome before the middle of July, as I already have informed you. I wish you a happy journey, free, as much as possible, from the discomfort of the hot season.

When in Geneva, please give my kind remembrances to Prof. Schiff[126].

I am sorry I shall not be able to see Tamanini and his bride when they come to Rome, as at that time I shall be a great way off from this city and above the sea-level. I cannot understand what attraction Rome and Naples under the August sun can have for nearly married persons, who are able to choose their route. Mergellina and Posillipo are a true Eden but their charms cannot be fully enjoyed when the season is too hot. Rather than in such a burning paradise, I would prefer to make my nest among the cool forests of Switzerland, near an imposing glacier, or on the shores of the lakes of Sweden, according to Prof. Mantegazza's advice, who is a very competent man on the scientific laws regarding love. When you write to Tamanini, please communicate my wishes for his happiness and my regret for the circumstances which will prevent my meeting him and his bride on their wedding tour. I hope to call on them, in Trent, when I visit the Tyrol again.

On the contrary I am sure I shall have the pleasure sometimes of seeing your brother Emilio and of being introduced to his lady, as they are to live a long while in South Italy.[127]

My wife has made up her mind to pass the holydays with the children at Perugia, an artistic and agreeably situated town, much more wholesome and less warm than Rome. They will leave at the end of July as soon as Vittorio will have got through his examination of "licenza liceale". At the end of August an artistic and agricultural exhibition will take place at which the inhabitants of Perugia expect a great concourse of visitors.

Give my and my householder's remembrances to Mrs. Clericetti and your boys. I shall look forward to your letter from Switzerland with great longing. Good bye, my dear friend

Yours sincerely
L.C.

---

[126] Non identificato.
[127] Si veda la lettera precedente.

**29.**

My dear Friend.

<div style="text-align:right">Milan July the 30<sup>th</sup> 1879.</div>

As you have not made yours appearance amongst us these last days, I must infer that you are still detained in Rome, either by the Senate alone or by some other engagement. The weather keeps beautiful and very probably we shall have a splendid autumn to enjoy, but we have not found yet a Country place for ourselves . In fact I have few days and had no time to enquire about villas or Cottages to be let, and less again to wander about in search of them. Next week also I shall be much occupied by my Exhaminations *[sic!]* at the Polytechnic and also on the first days of the next month, as it this year my turn to attend to the Esami Generali. In a few days you will I suppose leave behind our sunny shores and wander among the Alps to regain strength of body & mind amidst the splendid Scenery of a virgin nature. We have crossed the Simplon last month from Brigue to Iselle, stopping the night at the Convent. How I wished to be able to stay for a week at least, among those few Friars in those hospitable walls in that atmosphere of peace! There was a portrait of the first Napoleon hanging from one of the walls of the reception Salas: it was a very good oil painting of the first Consul dressed in the military dress of the time, when he crossed the Simplon and ordered the Convent to be erected. That portrait made a singular impression on every one of us, I believe: perhaps because it recalled to our memory the raising glory of the family, just a few days after the event that closed, it appears forever, the expectations of the descendant of a great man! How fate is unmerciful on the head of poor humanity!

Poor M$^r$ Favre, the eminent Engineer of the 1$^{th}$ Gothard tunnel, was not far from the day, when the two ends of his Gallery would meet, and fate strikes him to death, just in the bosom of the Mountain that he was perforating: the source of his glory was also the cause of his death![128] So it was with poor Mengoni: on the very Eve of pulling down all the scaffoldings that shaded the front of the Building he had just finished and to which all his glory is connected, he takes a false step and is plunged into the void space and crushed to death at the feet of his monument![129] Centuries pass on but the fate of humanity is always the same: to be ballotted *[sic!]* by a misterious *[sic!]* force like the waves by the wind, unconscious of past & future! The day before yesterday I went with my wife & children to spend the day on the Adda, at Vaprio, at Monastirolo at Trezzo. At Vaprio we admired the remnants, very rapidly faiding, of the Madonna by Leonardo which he painted in the house of his friend Francesco Melzi[130] at the same time when he was searching the laws of fluid motion by gazing on the rapid current of the river, minding at the feet of the palace. We spent an hour at Monastirolo an admirable example of the Italian Villa! Of the last two Centuries, that speaks of richness and family-greatness which also are now past for ever *[sic!]*. Then at Trezzo we wandered about those ruins of an historical Castle and saw the remaining Abutments of the famous Arch thrown across the river by Azzone Visconti, which was called by the contemporaries the 8$^{th}$ Marvel of the world. In fact it was the greatest Stone-arch which has ever been built up to this day, and it only lasted 70 years having been destroyed by the soldiers of the unhappy Count Carmagnola when he besieged the Castle.[131] Adieu my dear friend: hoping to see you in a few days, before your ramble on the Alps, present my duties to M$^{rs}$ Cremona & family if they are still in Rome and believe me ever

<div style="text-align:right">Your affectionately friend<br>C. Clericetti</div>

---

[128] Louis Favre, morto il 19 luglio 1879 per un infarto durante un sopralluogo al cantiere del Gottardo, sette mesi prima che le due parti si congiungessero.

[129] Giuseppe Mengoni morì a Milano il 30 dicembre 1877 precipitando dalla cupola centrale della Galleria Vittorio Emanuele II pochi giorni prima dell'inaugurazione.

[130] In realtà si trattava di Gerolamo Melzi, che ospitò Leonardo da Vinci tra il 1511 e il 1513; Giovanni Francesco è il figlio di Gerolamo.

[131] Il castello di Trezzo d'Adda e il ponte furono costruiti da Bernabò Visconti; il ponte venne distrutto nel 1416 durante un assedio da parte di Francesco Busone, detto "Il Conte di Carmagnola". Tale personaggio diede lo spunto al Manzoni per la tragedia omonima.

**30.**

Milano 22 Gennaio 1880

Caro Amico

Poiché insisti, eccoti i conti dei tre barattoli di Solfato di Chinino: per riguardo agli ultimi due, vedrai dalla fattura che furono da mio fratello commessi fin dalla fine di Dicembre: e furono spediti subito. Ma giorni sono capitò un avviso dall'Ufficio postale, qualmente all'Ufficio postale di Roma giaceva un pacchetto diretto a persona di ignota dimora nella Capitale. La persona di ignota dimora eri tu e nota che al tuo nome erano aggiunti i tuoi titoli di Commendatore, Senatore del Regno etc.!!!!! Ecco la causa per cui avrai ricevuto in ritardo il Chinino. Sono lietissimo di sentire da te che la tua signora si sia pressoché ristabilita, veramente lietissimo e la è con me la mia Cecilia, la quale pur troppo gode pochissima salute e dacché sono tornato di costì, non ha potuto uscire di casa più di 3 volte.
Per riguardo alla tua pensione all'Istituto Lombardo, ti dirò che me ne sono occupato immediatamente. Andai prima di tutto dal Poli, venerabile ottantenne dall'aspetto florido, dalla mente limpida, invidiabile. La sua signora era seduta al suo fianco sopra una poltrona, da dove non può muoversi dal 71 in poi per una paralisi alle gambe. Il Poli è tutto a tuo favore, ma soffre ai piedi; e il rigidissimo inverno lo tengono in casa. Andai dal Ceriani: era un po' incerto a chi dare il proprio voto per la pensione, ma dopo un lungo colloquio mi disse che ove ti fermassi dovessi salutarti a suo nome, il che equivale a dire che il suo voto sarebbe stato per te. Andai dal Biffi: mi dichiarò senza reticenze che egli avrebbe votato per te. Visitai il Sacchi: è malato pur troppo e nella impossibilità di presenziare la seduta, come il Poli. Però mi dichiarò che di pieno diritto la pensione spettava a te. Visitai l'Ascoli con cui pure sono in buone relazioni: il suo voto è assicurato per te. Andai da Polli ma non lo potei vedere: è una lampada che si spegne lentamente, fumigando. Lo viddi *[sic!]* quindici giorni fa e m'ha fatto senso: la sua mente è confusa: vi ha tuttavia nel suo intelletto un vano lucido che vorrebbe nascondere con ogni potere, la rovina del resto. Ieri era a letto e temo assai. Si tratta di una paralisi progressiva. Viddi *[sic!]* finalmente il Biondelli, altro mio amico: il suo voto è assicurato. Oggi poi fui all'adunanza all'Istituto, dove presentai quel lavoro sull'opera del Green (l'Area-momento).[132] La votazione deve essere stata fatta oggi stesso ed io non ho alcun dubbio dell'esito a tuo favore.[133] Non credo nemmeno che il Sangalli abbia tentato di scavalcarti, almeno per quanto è a mia cognizione.[134] Per riguardo poi alla mia nomina a Membro effettivo, ti dirò che son rimasto soccombente per 2 o tre voti al più nelle nomine fatte nell'ultima Adunanza. Gli amici mi hanno sostenuto, ma i Medici hanno portato il Taramelli che ha vinto. Vi ha, come saprai, un po' di attrito fra i matematici e i Naturalisti al nostro Istituto e le parti non cedono nelle votazioni. Spero però di riuscire ora che la morte del Frisiani lascia libero un posto da Membro effettivo.[135] Casorati me n'ebbe a parlare oggi stesso e così pure il Colombo. Voglio sperare che i Medici non vorranno questa volta rifiutarmi il loro voto. A quest'ora la votazione per la pensione è fatta ma non ho dubbio alcuno dell'esito. Se domattina, dopo la lezione potrò conoscere l'esito in modo indubbio te ne manderò notizia immediatamente. È ormai tempo di riprendere la nostra corrispondenza inglese: faccio voti perché il Polli possa guarire ma non lo spero e il suo stato mi dà una tristezza indicibile.
Addio. Tanti doveri alla tua Signora anche da parte della Cecilia e mille augurii perché si ristabilisca completamente assai presto. Credimi presto

aff.mo Amico
C. Clericetti

---

[132] Titolo della comunicazione di Clericetti: "Il metodo dell'area-momento nella determinazione delle condizioni di resistenza delle travi elastiche rettilinee" (*Il Politecnico-Giornale dell'ingegnere architetto civile ed industriale*, v. 12, 1880, pp. 194-203; 310-327); l'opera commentata è: C.E. Green, *Trusses and Arches Analysed and discussed by graphical Methods, New-York*, 1879.
[133] Nella seduta del 22 gennaio l'assegnazione della pensione venne rimandata, ma nella seduta successiva del 5 febbraio venne effettivamente assegnata a Cremona.
[134] A Sangalli verrà corrisposta la pensione a partire dal 1° luglio 1880.
[135] Clericetti sarà eletto Membro effettivo dell'Istituto Lombardo nell'adunanza del 24 marzo 1881.

**31.**

26-1-80

My dear friend

Your letters of 22 and 23 reached me; and by them I heard not without emotion how much you worked for my service in the affair of the Lombard Institute.[136] Your conduct towards me is truly brotherly; and never will I forget the assistance you gave me. Several of our colleagues addressed to me some lines expressing their favourable purpose; even Polli made her daughter write me in the same sense. The ballot war put off, I think, by the request of those members who were here in Rome for the parliamentary debates.

I hope you will succeed in the next election for fellowship at the L. I.[137]; I should be very glad in forwarding it. Can you suggest anything to be done by me?

I send you here enclosed, together with many thanks, 39 frcs and 35 cts, I owe you for the three small bottles of sulphat [sic!] of quinine you had bought for me. I hope we shall no more want quinine, at least for a long time. Since 14$^{th}$ December, my wife has been always free from fever and has made uninterrupted thought slow progress in her convalescence. She charges me with returning you and Mrs. Clericetti her best thanks for the interest you took in her illness.

Good bye, my dear friend. Please, give my kind regards to your lady and my and all my household's love to your boys.

Yours sincerely
L.C.

**32.**

Milan February the 14$^{th}$ 1880.

My dear friend.

I owe you a letter some days since, but having been lately much engaged in many things, you will excuse my apparent forgetfulness. I have received the money you have send [sic!] me for the two bottles of Chinino: there was no reason at all for such a hurry. I must congratulate you for the Pension, which has been conferred upon you by our Lombard Institute: one of the best sort of honours [sic!] which can be wished and sought for, and very much deserved by yourself. My wife is very far from being in a good state of health. She has passed a very bad winter having been obliged to keep her bed the most part of the time, without having suffered any real illness. Her constitution is very weak and faintly: sometime she feels well enough, but an hour of walk fatigues her to such an extent, that she is obliged to go to bed. I do hope that the coming spring may restore her strengths: the state in which she is now, and has been long since, without encluding [sic!] any danger, spreads a shade of melancholy aver all the family. M$^{rs}$ Cremona is well at last and we are very glad for it: poor Lady she has a long fight against the fever. I hope you have passed merrily the last days of the Roman Carnival: as for us, there is no more such a thing as the Carnevalone. Great efforts era made every year to revive the spirit of the people but without any success. There are now so many new sources of expenses for the families, and so different kind of amusements, that no money is left for the Carnival. Then, there is no denying that our population have become more serious and that the thought of future is not, as it was once, forgotten, in order to enjoy the present. We have got an indigestion of <u>Beneficenza</u>, under every shape. We must eat and drink, no matter if beyond necessity, for the sake of <u>beneficienza</u> we must dance and frequent theaters, no matter if beyond the private Budget, for the sake of <u>beneficienza</u>. We must work in order that the poor may enjoy life better, and get to misery ourselves, in order to be beneficiated in our turn. What a splendid circonlucution [sic!]! In a short time, healthy and wise people, shall be oblige to retire in order to give place to the <u>rachitici</u>, <u>the</u> <u>scrofolosi</u>, the <u>tignosi</u>, the fools of every description, the deaf and dumb, the crooked, the liberati dal Carcere etc etc.

What a magnificent result it will be to change the world into an immense hospital! Adieu my dear friend. Accept an english shake of hands and present my respects to all your family

You dearest friend
C. Clericetti.

---

[136] Si riferisce all'assegnazione della pensione; si veda la lettera precedente.
[137] Lombard Institute.

**33.**

Milano 27 febbraio 1880.

Caro Amico.

Dai colleghi Colombo ed Hajech, ho saputo or ora che sono fra i proposti per la Nomina a Socio effettivo dell'Istituto Lombardo. È questa la terza volta almeno che fui proposto. Si tratta ora di sostituire il defunto Frisiani. Le proposte furono fatte nell'Adunanza di Giovedì scorso giorno 19 e devono essere votate Giovedì prossimo 4 Marzo. I candidati sono quattro: la proposta firmata da maggior numero di membri fu quella del Prof Pietro Pavesi (Zoologo) che non conosco: poi vengo io con <u>due</u> firmatari di meno. Questa volta confesso che sarei proprio dolentissimo se non avessi a riuscire e poiché me l'hai promesso, domando il tuo efficace soccorso. I matematici hanno proposto me (Casorati Colombo, Hajech, Ferrini Celoria, Schiaparelli, Beltrami non era presente, ma non credo di dover dubitare del suo voto). Ma i professori medici di Pavia e gli altri cioè i Naturalisti proposero uno degli altri tre. Non so se il Cornalia abbia firmato la mia proposta, ma se stessi a quanto mi disse Colombo, pare che sarebbe in mio favore, ma non ne sono certo.
I due Cantoni (fratelli) il nuovo eletto Taramelli sono pel Pavesi. Insomma ora domando il tuo aiuto e te ne sarò molto grato. Se tu scrivessi al Prof Cantoni, il fisico, al Verga (Il Biffi no perché appartiene all'altra Classe) al Cantoni Gaetano e a qualche altro, credo che riusciresti a fare in modo che nella prossima votazione il Pavesi non raccolga i ⅔ dei voti dei presenti che si richiedono.[138] Perdonami il disturbo, ma aiutami se credi di poterlo fare e non ne dubito perché l'hai promesso. Addio: tanti doveri alla famiglia. Tronco per impostare subito

aff.mo amico
C. Clericetti

*[Segue elenco scritto da Cremona:*
Verga
Garovaglio
Polli
Cornalia
Brioschi
Hajech
Stoppani
Schiaparelli
Mantegazza
Cantoni Giovanni
Cremona
Sangalli
Casorati
Colombo
Ferrini
Corradi
Cantoni Gaetano
Celoria
Beltrami
Maggi*]*

---

[138] Durante la seduta del 4 marzo, nessun candidato proposto raggiunse il numero di suffragi richiesti. Si rimandò quindi l'elezione alla prima seduta ordinaria di giugno. Durante quella adunanza – 17 giugno – ancora non si raggiunse il quorum. Si rimandò al 29 luglio 1880, quando venne eletto Guglielmo Körner. Pietro Pavesi sarà nominato Membro effettivo il 22 febbraio 1883.

**34.**

3-3-80

My dear Clericetti.

Your last letter had scarcely reached me, than I wrote to Cantoni Giovanni, Verga and Cornalia, in order to second your election as a M.E.[139] of the L. Institute. I did not wrote to Cantoni Gaetano, because I (have not enough intimacy) am not intimate enough with him but I begged the Physicist to speak with his brother about the affair. Likewise I suggested to Cornalia and Verga to engage their friends in forwarding your election. It is a pity that poor Polli is so ill and perhaps at the (end of his life) (point of death); and that Brioschi, Stoppani, Mantegazza and myself are away from Milan. I am very sorry that in this moment I am not able (to absent myself from Rome ) (to take a run away from Rome.)
But I hope that, in spite of these unfavourable circumstances, you will succeed, <u>quod est in votis</u>.
(I neglected telling you that, having occasion to write to Hajech before I knew he was one of the subscribers for you, I suggested to him to promote your election).
Perhaps I have not yet thanked you for the friendly and affective interest you took about the affair of my pension.
I assure you, our friendship always grows more (soldered) riveted and therefore I eagerly wish that the Institute may recognize your services in the cause of science and confer upon you the deserved prize.
My wife's health continues good enough; the little remnants of her illness will vanish only by time, the warm season and charge of air.
We all are very sorry for your lady's bad state of health; can you not think of carrying her on the high Alps, ex. gr.[140] in the fine and green Engadin? It seems to me, that pure, delicious air would recall the dead to life!
Please give her our kind regards and always believe me

Yours sincerely
<u>LC</u>.

*[nota indirizzata a Cremona:*
Dear Professor,
here is your letter.
Thanks for Wednesday . If not this week perhaps some other Wednesday.

Yours sincerely
Mary Isabella Taddeucci*]*

---

[139] Membro Effettivo.
[140] Exempli gratia.

**35.**

Rome, April the 2$^d$ 1880

My dear Clericetti

I was very surprised to receive your letter from England; and after having read it, I saw once more that you possess the energetic travelling spirit of a true Englishman. I am very sorry for Mrs. Clericetti; she must have been much afflicted in leaving her poor sister. Have you received further tidings of her?
As to the affair of the Lombard Institute, I hope you will oppose your English spirits to the temporary unsuccessfulness. You must behave like a true Philosopher and see in the fact only a consequence of Italian idleness. As you know, I did not neglect writing to my friends; but also your rivals did not remain inactive. One of them wrote even to me. This is sufficient to explain the division of the votes which permitted none of the candidates to succeed in the election.
Many thanks for your purpose of sending me a copy of the English translation of your paper on American suspension bridges.[141] Do you not think that the suffrage of the English Engineers is worthier of being desired than that of the Milanese physicians and Pavese Naturalists?
You were surprised to find the Quaternions taught in Owen's College at Manchester. Quaternions are speciality of the British Mathematicians; besides Hamilton's original and masterly works, Great Britain possesses many, more or less elementary, treatises or introductory books on that subject. But continental mathematicians see in Quaternions only a method, useful in treating some arguments, which however are also accessible by other analytical proceedings; while one could not make one's self master of the Quaternions without great trouble and long waste of time. And after all, though in England there are several eminent Mathematicians, the level of public teaching in higher departments of Mathematics are, generally speaking, not very lofty.
And now good bye, my dear friend. Give my kind regards to Mrs. Clericetti and my love to your boys.
I am ever

Yours sincerely
LC.

**36.**

1880, June, 4

My dear Clericetti

I have just received your English paper on the Theory of modern American suspension bridges.[142] Intending to read it thoroughly as soon as I shall be a little less busy. I must not put off the hearty thanks I owe to you. I hope your lady and your boys are in quite good health. It is not so with my wife, who since many days is forced to keep her bed. As soon as she will be able to get up, I will take her to the Bagni di Lucca, where my family will pass all the holydays.
Good bye, my dear friend
Your sincerely

---

[141] C. Clericetti, "The Theory of Modern American Suspension Bridges", *Minutes of Proceedings of the Institution of Civil Engineers*, v. 60, 1880, pp. 338-359; ristampato in: *Van Nostrand's Engineering Magazine*, New York, v. XXIII, 1880, pp. 111-122.
[142] Si veda lettera precedente.

**37.**

Milan June the 8${}^{th}$ 1880.

My dear friend.

I am sorry to understand that M${}^{rs}$ Cremona is not well yet, but am sure that the Sea air and water will do more than medicaments to restore her health. My wife has returned from England a fortnight ago, leaving her sister in a state better than could be hoped. She cannot recover completely but she may live on. In the meanwhile, for more than two months I was left alone with our children, a thing which would have been impossible if they were babies; but they are quite young men and taller than myself.

I read with pleasure your first Discorso to the Senate,[143] of which you kindly send a copy to me, and I quite agree with your ideas on the Superior Council of Public Instruction, and about Universities. The nomination of an effective member at our Lombard Institute, is to take place on the next meeting that is on the 17${}^{th}$ and the proposals must have been made on Thursday last (the 3${}^{d}$). I was at the meeting, but there were very few members. Casorati Beltrami, Celoria were absent, perhaps on account of the funeral of one of the Professor at Pavia; of our friends, only Colombo, Cornalia, Ferrini, Schiaparelli were present at the meeting, and also Hajech, but I know nothing of the result. Polli will never again be in a state to assist to any meeting, and I fear, not even to go out from the house. The finishing up of a long life is such a sad thing! I saw him last night: he complains a weakening of his constitution and his mind is not steady. He seems to have received a shock from M${}^{rs}$ Kramer's death: he speaks of her as if she was still alive and as if he had to attend to her every day. The weather is good at last and it is to be hoped that it will continue so. Nothing new in our city: the wealthy people are occupied with the silk worms though it is a source of revenue nearly exhausted for Lombardy. Then the summer begins and Milan is left a desert. I have not found yet a house for the next Autumn, so I do not know what one shall do; but I will look about again. Cecilia sends many salutes to M${}^{rs}$ Cremona and yourself and family: she came back from England in a better state of health than she was for many a month. The journey, the visit to her country and relations did her good. She knows that she will not see again her sister, who can not live long, but she has seen her and nursed her for two months, and she feels content.
Believe me ever

Your devoted friend.
C. Clericetti

---

[143] Si tratta di: "Riforma del Consiglio Superiore della Pubblica Istruzione: discorso del senatore Cremona pronunziato al senato nella tornata del 9 aprile 1880", Forzani e C., tipografi del senato, Roma, 1880.

**38.**

Milan June the 18[th] 1880

My dear friend.

I suppose you are Aware that yesterday the Lombard Institute had to vote for the nomination of an effective Member in the Class of Mathematical Sciences; but very probably you do not know yet, that the result has been the same as usual, from two years since. Nobody has been elected and the votation *[sic!]* is again postponed to November. They had three successive votations *[sic!]* and in every one of them I had 7 votes given to me and 9 against me. I really feel provoked to see this constant war against myself by Doctors who go to the Institute to read such trash as they do.

To be sure I will some day succeed, I do not want to be received as for Charity sake. Now that D[r] Polli is dead there is another place vacant, but I do not want to take his place, but the place wich *[sic!]* was only vacant untill *[sic!]* these last days, I mean to say the place of Frisiani. If you and Brioschi and Mantegazza had been there, I would have succeeded, and <u>I must really entreat you for friendship's sake to come to Milan for the next votation</u> *[sic!]*, <u>as I beg you to induce also Mantegazza to do the same</u>. Colombo wrote to me last night all the particular's adding "Io non ci capisco più nulla: nessuno vuol cedere malgrado le più grandi proteste di accordarsi su un nome".[144]

<u>You must really do me the favour to come to Milan next November</u>. It is more than two years that things are going on at the Lombard Institute in a way that is not the best for its prosperity: this antagonism must cease either by mutual concessions or in some other way. Professor Polli has ceased to suffer. Galanti[145] and myself read something at the Cemetery before his Coffin. Not one Medical man said a word in praise of a man whose life had been all dedicated to the improvement of Science! Only later when the body was cremated with Gorini's apparel,[146] D[r] De Cristoforis said something in memory of the deceased, but only relating to his part in the diffusion of Cremation: only as a réclame to Cremation! I was not present because I have completely withdrawn myself from the Antipatic *[sic!]* theme, sorry that my name has been mingled with some others which I despise very much. Adieu my dearest friend: my best regards to M[rs] Cremona and family.

Your devoted
C. Clericetti

---

[144] Né Cremona né Mantegazza saranno presenti alle successive adunanze, del 1/7, 15/7 e 29/7, quando viene nominato Guglielmo Körner; saranno invece presenti in quella del 24 marzo 1881, quando finalmente viene nominato Clericetti. Membri effettivi presenti: Cremona, Cantoni Carlo, Belgiojoso, Verga, Biffi, Biondelli, Cossa Luigi, Hajech, Beltrami, Celoria, Poli, Maggi, Schiaparelli, Casorati, Körner, Corradi, Mantegazza, Colombo, Sangalli, Ceruti, Ferrini, Strambio; Cornalia è assente per malattia.

[145] Si tratta probabilmente di Antonio Galanti, che scriverà una commemorazione di Polli insieme a Cornalia.

[146] Quando, nel 1874 Clericetti e Giovanni Polli, avevano costruito il primo forno crematorio a gas nel Cimitero Monumentale di Milano (si veda la lettera 11), Paolo Gorini - che già stava studiando il processo di cremazione - considerò che l'impianto di Polli-Clericetti avrebbe potuto essere migliorato e nel 1877 inaugurò presso il cimitero di Riolo (LO) un nuovo sistema di forno crematorio che utilizzava la combustione di fascine di pioppo ed era quindi più semplice ed economico. Tale forno crematorio venne installato in tutto il mondo e ha funzionato fino agli anni '80 del Novecento. Si può dunque comprendere il tono piuttosto risentito di questa lettera.

**39.**

17 Nov. 1880

My dear Clericetti

You are perfectly right, but I am not guilty on account of the Biblioteca V. E.[147] and of a general inquiry on public libraries, I am so very busy that I can no longer dispose even of a minute for myself and my friends. On this account, your letters still lie on my writing table. For the last five months, I have not been able to write a single English line! But, I beg you as warmly as I can, not to discontinue sending me your kind English letters, which give me greater pleasure than you perhaps may think.
But I do not overlook my duties as your friend. As soon I received your letter of Nov. 5th, I sent to the Secretary of the Lombard Institute my proposal with your name. And I repeat now my old promise: when it will be the matter of the election, I will run down to Milan, notwithstanding my oppressive business.
I hope to be able another day to write you a longer letter. We are in good health; I desire the same to you and yours. Give my respects to Mrs. Clericetti; and good bye!

Yours sincerely
LC.

**40.**

My dear Clericetti

Your letter of the 3d Dec.[148] caused me great astonishment and vexation. Indeed you have not good luck! I was preparing everything in order to repair to the station and come to Milan, when your telegram reached me!
But, cheer up! You must be calm and wait for a better future. Without doubt, I will be at Milan, in March. Let us hope that your other friends will not fail in that final[149] trial.
I send you, for the new year, my best wishes for you and your family: be all happy and healthy for ever! And if you love me, you must wish me to be delivered from my chains.
Good bye. I always am

Yours sincerely
LC

---

[147] Biblioteca Vittorio Emanuele. Dal giugno 1880 F. De Sanctis, Ministro della Pubblica Istruzione, aveva nominato Luigi Cremona Regio Commissario per il riordinamento della Biblioteca Vittorio Emanuele di Roma che si trovava da tempo al centro di accese polemiche. Cremona, dopo circa due anni di accurate inchieste, provò molte delle accuse formulate.
[148] Questa lettera non riporta la data, però potrebbe inserirsi qui.
[149] Cremona aveva usato il termine "utmost". La persona che gli correggeva le lettere lo sostituì con la parola "final" e aggiunse a margine la seguente spiegazione: *"utmost trial" sarebbe quel eccesso di dolore e di prova che si può sopportare. Non vorrà significare ciò, m'immagino.*

**41.**

Milano 16 Marzo 1881.

Caro Amico.

Aspetto un tuo cenno per sapere se debbo fare qualche passo personalmente presso qualcuno, relativamente alla prossima elezione all'Istituto Lombardo.[150] Il Brioschi prima di partire ha scritto al Cantoni Gaetano e al Körner e mi ha dato incarico di scrivere al Casorati e al Beltrami a Pavia. Non l'ho fatto perché credo che tu abbia già parlato con loro ed aspetto in proposito istruzioni da te. Com'ebbi a dirti, il Cantoni Giovanni mi promise il suo voto senza reticenze, e a parte scienziato, mi pare uomo sincero. Del Cornalia piuttosto dubito ad onta che mi abbia promesso e parecchie volte, di votare per me. Nell'ultima Adunanza fu tirata in scena anche la proposta del Golgi, distinto professore, a quanto sento, ma più giovane di me e nominato da poco Socio Corrispondente.

Il Brioschi mi disse di aver parlato anche con Sangalli, promettendo di appoggiare il Golgi nella prossima futura occasione, se in questa occasione il Sangalli votasse per me. Il Sangalli però avrebbe solo risposto "Chissà". La tua presenza farà molto e non trovo modo di dirtelo appieno e di ringraziarti quanto meriti. Desidero però che mi si presenti una occasione se non altro per farti comprendere quanto ti sono grato. Posso io contare su Casorati e su Beltrami, o debbo scrivere, o debbo fare una corsa a Pavia?

La prima volta che tu avrai occasione di recarti al Ministero della Pubblica Istruzione, mi farai il favore di informarti riguardo ad un Mandato di pagamento per £ 225.64 che sono dovute fino dalla fine di Novembre. Nel Novembre (25-28) fui costì per l'Esame dei titoli e la proposta di nomina di un Professore per la Cattedra di Scienze delle Costruzioni a Napoli, alla quale fu nominato il Bonolis.[151] Il 13 successivo Dicembre una lettera firmata dal Tenerelli, allora tuttavia Segretario Generale n° di posiz$^e$ 66. N° del Prot. Gen. 4995. N° di partenza 8233 mi avvertiva che con Decreto del medesimo giorno era stato provveduto al pagamento di tale indennità e aggiungeva "Ella riceverà avviso quando il relativo Mandato sia ammesso a pagamento presso codesta Tesoreria". Ma l'avviso non è mai venuto. Del resto questo con tutta tua comodità.

Nessuna notizia relativamente alle elezioni del Consiglio Superiore: riguardo poi al tuo birbone della Perseveranza[152] a tuo riguardo, mi pare (e non credo che sia troppo tardi) che tu dovresti smentire la cosa: il tuo soggiorno nella Capitale e quindi lontano dalla sede di quel giornale può essere un eccellente motivo da allegare per aver ignorato la notizia pubblicata dalla Perseveranza, fino a questi ultimi giorni. Dunque a rivederci presto. Tanti doveri alla signora Elisa ed augurii alla promessa sposina.[153] Accogli una stretta di mano

dell'aff$^{mo}$ Amico
C. Clericetti

---

[150] Si vedano le lettere precedenti.
[151] Potrebbe trattarsi di Alfonso Bonolis, traduttore con Francesco Mazza del *Manuale di formole, tavole e notizie di uso frequente agli ingegneri, architetti...* (di J. Claudel), Pellerano ed., Napoli, 1879.
[152] *La Perseveranza: giornale del mattino* era un quotidiano milanese conservatore che iniziò le pubblicazioni il 20 novembre 1859. Era il giornale dei grandi proprietari terrieri lombardi; diretto sino al 1866 da Pacifico Valussi, poi da Ruggero Bonghi fino al 1874. Sotto la direzione di Carlo Landriani iniziata dal 1880 si avviò versò un lungo declino e chiuse il 20 maggio 1922 mentre era direttore Tomaso Borelli.
[153] Si tratta di Elena. Si veda anche la lettera successiva.

**42.**

Milano 28 Marzo 1881.

Caro Amico.

Spero che il viaggio di 40 ore a cui ti sei sobbarcato per rendere un servizio ad un amico, per dimostrare una volta di più che non sei solamente un illustre scienziato, ma un uomo di cuore, non abbia recato nocumento alla tua salute. L'Istituto Lombardo mi ha comunicato con lettera del Presidente Cornalia, la nomina, e vi ho risposto oggi in brevi parole. Il Brioschi che mancava all'Adunanza di Giovedì, volle però essere presente e presiedere la riunione dei professori dell'Istituto Tecnico Superiore, tenuta Sabato sera, per prendere accordi sulla Nomina del Consiglio Superiore. Quando si tratta del suo interesse personale, diretto, non manca mai. Egli comunicò varie liste di nomi che erano state proposte da diverse Università[154] e fra altre, anche quella di Napoli che portava anche il nome tuo. Si discusse alquanto su tali liste. Pavia accettava Brioschi purché noi proponessimo Cantoni.[155] Pisa voleva Betti etc. Si discusse alquanto, ma si chiacchierò poco, poco assai, com'è nostro costume. Avrei amato che per l'Ingegneria fosse scelto il Curioni e per mio conto respingevo Cantoni, come pure il Cannizzaro non si voleva da molti. Però parecchi colleghi mi dimostrarono l'utilità di non offendere Pavia, con cui già da anni esistono attriti che sarebbe meglio sopprimere. Allora accettai il Cantoni. Ieri mattina si fece la votazione. Posso assicurarti che il tuo nome figura in quasi tutte le nostre schede: dico quasi perché so che in una non figurava e quella scheda non portava di comune colle altre cioè con quelle convenute, che il nome di Brioschi. Tu hai dunque ottenuto di certo un buon suffragio dalla nostra Scuola, come lo avrai da Napoli certo e da altre Università.
So che a te importa poco perché sei sicuro ugualmente di entrare nel Nuovo Consiglio,[156] ma non potrai a meno d'essere lieto del voto dei tuoi vecchi Colleghi. Adesso il Vecchio Consiglio Superiore dovrà suicidarsi coll'aprire esso stesso le 700 e più schede che creeranno il Nuovo. Faccio voti perché faccia presto. L'Elenco dei professori universitari, pubblicato dal Ministero, sembra stato fatto con poca diligenza perché sonvi parecchi errori. Io intanto sono indicato come Professore di Scienza delle Costruzioni e il Saino come Professore di Ponti.[157] Pare impossibile che si facciano le cose con tanta leggerezza. Abbiamo la città molto animata; assai più dell'ordinario. Fa un curioso senso l'incontrare sul corso Venezia, dal Ponte al Dazio, dunque in piena città, una vera locomotiva che trascina un vagone carico di prodotti per l'Esposizione.[158] L'occhio, come tutto, è schiavo delle abitudini. Abituato a vedere tali macchine in aperta campagna o sotto le tettoie delle stazioni, l'incontrare una vera locomotiva da ferrovia ordinaria in una strada di città, fa un senso curioso almeno a me. Mi par più grande più colossale del solito. Gli è come quando mi accade di vedere una persona passeggiare nella parte piana del nostro giardino pubblico, dove le piante sono basse e i cespugli minuscoli. La persona mi par più grande del naturale, perché non vi trovo il rapporto ordinario tra l'altezza di un uomo e quella di una pianta. Ma queste sono chiacchiere e te ne domando perdono. Noi partecipiamo vivamente alla Gioia della tua famiglia all'avvicinarsi dell'epoca del matrimonio della sig.na

---

[154] Sulla G.U. del 18/03/1881 si trovano le "Norme per l'esecuzione della legge 17 febbraio 1881 sul Consiglio superiore di pubblica istruzione". In particolare l'Articolo 1 recita: "I professori ordinari e straordinari delle Regie Università, delle Scuole di applicazione per gl'ingegneri dell'Istituto tecnico superiore di Milano, dell'Istituto di studi superiori di Firenze, dell'Accademia scientifico-letteraria di Milano e delle Scuole superiori di medicina veterinaria, nel giorno indicato da apposita circolare ministeriale, si adunano per designare, mediante votazione con schede segrete, sedici persone che, in conformità dell'articolo 2 della legge 17 febbraio 1881, saranno dal Ministro proposte alla nomina Regia per far parte del Consiglio superiore."
[155] Giovanni Cantoni.
[156] In effetti, Cremona verrà nominato membro su proposta ministeriale il 12 maggio 1881, ma si dimetterà all'inizio di ottobre (si veda la lettera del 25 luglio a Clericetti). Sarà presente nel Consiglio superiore più volte negli anni e ne diverrà vicepresidente (01/05/1888-31/05/1889; 01/06/1890- 31/05/1894).
[157] Anche Antonio Sayno era docente di Scienza delle costruzioni.
[158] Esposizione Nazionale, tenuta a Milano dal 5 maggio al 1° novembre 1881. Clericetti tenne una conferenza sui grandi manufatti eseguiti in Italia e fece parte della Giuria dell'Esposizione come Segretario e Relatore per l'Ingegneria ed i Lavori pubblici. Si veda *Esposizione nazionale del 1881 in Milano - Relazione generale*, Tip. Bernardoni, Milano, 1883.

Elena.[159] È un doloroso distacco sempre, per un padre per una madre e mi sovviene l'Amaury del Dumas,[160] ma sono lagrime che somigliano a quelle provocate dal riso: è un allegro pianto. Tanti doveri alla signora Elisa e a tutta la famiglia. Ricevi ancora un ringraziamento ed una stretta di mano.

aff.mo
C. Clericetti

**43.**

Milano 25 Maggio 1881.

Caro Amico.

L'Esposizione[161] è aperta da tre settimane e noi aspettiamo di giorno in giorno che tu mantenga la promessa fattaci di venire: la tua camera è pronta da un pezzo e non ci manca che la tua desiderata presenza.

Non so se i tuoi sposi[162] siano venuti qui, o ci siano attualmente. Vorrei che così non fosse perché abbiamo sempre contato sopra una loro visita. L'Esposizione è riuscita assai interessante ed ha superato l'aspettativa generale, non tanto per la quantità dei prodotti esposti, quanto pel buon gusto nella disposizione e distribuzione d'ogni parte e per quell'impronta artistica che caratterizza ogni prodotto delle nostre industrie, dai più ricchi ai più modesti ed usuali.

Questa mi pare la vera nota dominante della Mostra: la scelta dei giardini fu infelice perché la città ora manca di un ritrovo geniale pel dopopranzo. Ma d'altra parte ne ha avvantaggiato molto l'effetto generale e tanto più per l'aggiunta del Giardino Pubblico della Reale Villa che è una vera oasi. La crisi ministeriale che perdura da tempo e alla quale finora non vedesi una fine è stata certo una disgrazia pei Milanesi che hanno fatto sforzi e spese grandi, nella speranza di un concorso grandissimo di tutti gli italiani: si contava naturalmente, benché meno, anche sugli stranieri ma la faccenda improvvisa di Tunisi[163] e la irritazione nazionale contro i francesi che vi è sorta e che anzi va crescendo, fornì altra delusione. Erano in vista delle corse di piacere da Parigi a qui: ora non bisogna pensarci più e contare solo sui nazionali. Però sono certo che l'Esposizione piacerà e ti ripeto che noi aspettiamo una tua lettera che ci avverta della tua venuta. Ma bisogna venire per restarci una quindicina di giorni e però scegli bene il tuo tempo perché non ci va meno. Ho visto dai giornali che anche tu fosti ufficiato dal Sella ad assumere il Ministero dell'Istruzione, ma pare che tu abbia rifiutato.[164] Pel Sella il non essere riuscito a comporre un Ministero nei centri, fu cosa proprio deplorevole e che dispiace tanto a chi lo ama e lo stima.

---

[159] Elena andrà a vivere a Macerata con il marito Adolfo Perozzi.
[160] Novella di Alexandre Dumas padre, composta nel 1843.
[161] Si veda la lettera precedente.
[162] Si veda la lettera precedente.
[163] Nel 1881 il governo della Terza Repubblica francese stabilì con un'azione di forza il protettorato sulla Tunisia, che era obiettivo dichiarato dei propositi coloniali del Regno d'Italia. Si venne così a creare una situazione di crisi politica tra i due stati; questo episodio venne descritto dalla stampa come "lo schiaffo di Tunisi". Alla data di stesura della presente lettera, il protettorato francese sulla Tunisia era appena stato ufficializzato con il Trattato del Bardo (il giorno 12 maggio).
[164] Cremona rifiutò l'invito di Quintino Sella con questa lettera del 19 maggio: "Carissimo Amico, Tu mi facesti un'offerta, che mi resterà come uno dei più preziosi ricordi della mia vita. L'essere da te stimato capace di venirti in aiuto nella difficile impresa a cui ti sei sobbarcato è per me altamente onorevole e lusinghiero. Ma potrei io darti un aiuto efficace? aggiungerei io forza al tuo Ministero? Non lo credo. Ad ogni modo le mie opinioni e i miei precedenti politici mi vietano di pormi contro i caduti, associandomi ai successori. Perdona illustre amico, se ti mando per iscritto una risposta diversa da quella che tu, per tua benevolenza mostrasti desiderare, che io ti portassi a voce. Se ricuso l'alto onore non è per cagion tua, ma della situazione, la quale m'impone dei doveri che tu certamente vorrai apprezzare, pur non approvandoli. Confido che mi resterà intera la tua amicizia, come io sarò sempre Tutto tuo L. Cremona". La minuta di questa lettera si trova alla segnatura 051-11846 (8533) ed è stata pubblicata in A. Brigaglia, S. Di Sieno, "L'opera politica di Luigi Cremona attraverso la sua corrispondenza, Seconda Parte. Il crollo delle speranze e il lavoro organizzativo", La Matematica nella Società e nella Cultura, *Rivista della Unione Matematica Italiana*, (I), 2010, pp. 137-179.

La continua persistenza poi del Farini a non volersi sobbarcare al compito di comporre un Ministero o ad entrare in una combinazione, è inesplicabile con criteri ordinari perché, non raggiunge così l'intento, se pure è tale il suo pensiero, di non voler essere demolito. La sua astensione ostinata lo demolisce e a spiegarlo non mi rimane che di sospettare che sia vero quanto mi venne detto del Farini. Che cioè egli appartenga segretamente , volente o meno, alle secrete sette della Romagna,[165] e che tale legame gli impedisca di accettare un Ministero. Lo si dice anche costì?

Troverai Milano molto animata. Esposizione, Ind Teatri, Panorama, Ippodromi, Tramway, ferrovia elettrica, Esposizione Musicale, Esposizione Artistica, Pallone areostatico (non compiuto) etc etc. Non si studia non si lavora più: i milanesi si moltiplicano e vanno innanzi e indietro dall'Esposizione, occupati a curiosare. C'è un brulicame tutto il giorno là che pare un grande Alveare. Le signore non vanno ai bagni ed hanno anche sospese le solite visite per dedicarsi all'Esposizione: e là tengono i loro crocchi quotidiani. Dunque avvertimi quando potrai venire. Tanti doveri alla signora Elisa e saluti ai figli. Accogli nell'attesa una cordialissima stretta di mano e credimi sempre

Devot$^{mo}$ Amico
C. Clericetti

**44.**

Milan July the 25$^{th}$ 1881.

My dear friend.

We have been expecting you week after week, these last two months, as you had promised to come again and visit our Exhibition.[166] But you did not come, nor send word you would: we just saw for few moments your darling sposi,[167] who seemed also in a hurry to get away from Milan. However we hope to see you in September when the hottest part of the Summer will have passed and the Exhibition will be more frequented than it is lately. On Saturday (the day before yesterday) we had a very agreeable visit at our Institute from your Scholars, accompanied by Prof$^{rs}$ Guj, Sinigaglia and Ceruti.[168] Brioschi was absent from town, but the Professors were all present to receive them. Your Scholars left a very good, very agreeable impression and I am very glad to relate it to you. My Guido has passed well the Exhaminations [sic!] which he had to repeat for the Licenza liceale, and has also passed well enough the Exhamination [sic!] of Geometria proiettiva under Jung: so things are progressing in a good direction with him. As for Emilio he passed as usual without Exhaminations [sic!]. I have been with Colombo and Martelli to Genova to Spezia, Livorno and Florence with our Scholars of the last year.[169] The works in the harbor of Genova are very well conducted as usual by our Genio Civile. Maganzini is there as active and intelligent as ever: he lately married and brought home at Genova a very nice and handsome Sposina from Verona. We visited the Dandolo[170] nearly finished in the Arsenal of Spezia and the Lepanto[171] in construction by the Brothers Orlando at Livorno, two great iron monsters like the Duilio,[172] and even bigger, but very probably as useless: at least we can say that they are <u>incognéte</u>.

We were all very sorry to learn of your resignation from the Superior Counsel of the Public Instruction. Your real motives are not known to us, but we suppose you resigned as a demonstration against some of the

---

[165] Forse si riferisce al fatto che Domenico Farini - seguendo le orme di suo padre, Luigi Carlo - fosse massone.
[166] Si vedano le lettere precedenti.
[167] Si vedano le lettere precedenti.
[168] Enrico Guj, Francesco Sinigaglia e Valentino Cerruti. Verrà pubblicata la "Relazione del prof. Francesco Sinigaglia intorno al viaggio fatto dagli Allievi della R. Scuola d'applicazione per gl'Ingegneri di Roma", Tip. Barbera, Roma 1881.
[169] Si veda: O. Garuti, P. Orlando, "Relazione del viaggio d'istruzione fatto dagli allievi ingegneri civili del R. Istituto Tecnico Superiore di Milano nell'anno 1880-1881", *Il Politecnico - Giornale dell'ingegnere architetto civile ed industriale*, v. 14, 1882, pp. 53-65; 196-218.
[170] Corazzata costruita su progetto di Benedetto Brin, Direttore del Genio Navale; fu completata nel 1882.
[171] Corazzata costruita su progetto di Benedetto Brin, Direttore del Genio Navale; fu varata nel 1883.
[172] Corazzata gemella della Dandolo. Venne costruita nel Cantiere navale di Castellammare di Stabia e varata nel 1876.

Decrees of the present Minister[173] relating the studies of Mathematics in the Licei, which have been reduced almost to nothing. I hope and wish that you will not be long out the Counsel. The season is very hot and dry: one fancies and dreams, as always in this season, of the Green valleys of the Alps, of the running streams of the thick woods. But we cannot go into the country: my Exhaminations *[sic!]* at the Institute will not take place until next month, at the end of which I must be in Milan as one of the Giuria of the Exhibition.[174] Somewhere we shall go for sometime but cannot say where or when. M$^r$ Betocchi is everlasting backwards and forwards from Rome to Milan: he does so like to explain in a sort of Conference which lasts even two hours, the Exhibition made by the Ministry of Public Works; every time he learns that one of the Schools of Engineers is coming, he starts, comes here remains a day and runs back. Adieu my dear friend: our united compliments to M$^{rs}$ Cremona. Hoping to see you, receive an English shake of hands and believe me ever

<div style="text-align:right">Your affectionate friend<br>C. Clericetti.</div>

**45.**

<div style="text-align:right">Milano 29 Ottobre 1881.</div>

Caro Amico.

Quello che ho potuto fare io per te è ben poco a paragone di quanto hai fatto per me e per cui ti sarò grato tutta la vita, come ti è gratissima tutta la mia famiglia. Tu sei uno scienziato di cuore: lo dissi tante volte e i fatti mi hanno sempre dato ragione. Ed è proprio vero quello che dici che quando s'è oltrepassato "il mezzo del cammin di nostra vita", quando le frivole passioni non hanno più presa in noi, il pensiero si raccoglie volentieri intorno agli amici più cari e il loro esempio serve di sprone a proseguire nel cammino della vita: e lo spirito nostro attinge coraggio dal coraggio degli amici.

L'Esposizione[175] sta per essere chiusa: il tempo è triste e una pioggia minuta ha scemato il freddo degli scorsi giorni, ma l'umidore penetra fino nelle ossa. E non di meno la città è assai affollata: sono le migliaia di persone che hanno aspettato gli ultimissimi giorni per visitare la mostra: I Cafè, i restaurant sono rigurgitanti e percorrendo la sera il Corso, giù fino alla Piazza del Duomo e nella Galleria è uno spettacolo curioso il vedere tante centinaia di bocche in movimento: ogni fatica va a finire in appetito. Presto però la Città riprenderà la sua calma abituale: ognuno rientrerà nella propria orbita o in quella che s'è acquistata nuovamente con le Onorificenze e il ricordo dell'Esposizione sarà pei milanesi come quello di una brillante meteora. Però qualche cosa rimane ed è una maggiore fiducia degli italiani in se medesimi e un'aspirazione meno vaga al progresso. Se tu poi fossi venuto prima e con maggior tempo disponibile e ci fossimo venuti insieme mediante i molti tram a vapore, qua e là nella campagna milanese e nella Comasca, avresti riscontrato i sintomi di un vero rinnovellamento generale: dappertutto sorgono opifici e alti camini industriali e nuove case di campagna e villette circondate da giardini. Dappertutto anche nei più modesti comuni sorgono locali per gli uffici municipali e per le scuole e le vecchie case coloniche e i Cascinali si vanno imbiancando come si preparassero ad una festa: insomma è uno spettacolo consolante questo della patria nostra che finalmente si è scossa dal lungo letargo e vuol progredire e rifarsi a nuovo migliorando tutto. La Cecilia ha ricevuto stamane una carissima lettera della signora Elena ed è lietissima che un dono così modesto, così piccino, le sia piaciuto. Ella risponderà, e intanto invia tanti saluti a te alla signora Elisa agli sposi, all'Itala e al Vittorio. Attendo di conoscere i tuoi desiderii relativamente all'abito da £ 38. e a qualsiasi altra cosa o informazione che ti bisognasse. Io ho passato due giornate a Lucino presso i nostri vecchi amici, gli Olginti di Como: ho visitato, come a riposo della mente, una vecchia chiesucola lombarda,[176] dimenticata nella solitudine della campagna, nascosta da alte piante secolari e dove una vasca che ha servito per tanti secoli per l'acqua lustrale, fu nei tempi romani una bellissima urna sepolcrale che in

---

[173] Guido Baccelli, Ministro della Pubblica Istruzione dal 2 gennaio 1881 al 30 marzo 1884.
[174] Si vedano le lettere precedenti.
[175] Si vedano le lettere precedenti.
[176] Si tratta di S. Bartolomeo al Bosco ad Appiano Gentile (CO).

cinque diverse cavità, dovette conservare i resti cremati di un'intera famiglia. Poi nel X° secolo un Suddiacono Valperto se ne valse per seppellirvi i resti dei suoi parenti e lasciandovi le iscrizioni romane, ne aggiunse scioccamente altre, per dirsi lui stesso autore dell'urna![177] Sono mistificazioni di cui il Medio Evo ha lasciato tanti esempi. Ma alcune di quelle iscrizioni sono assai singolari e accennano a quel riposo assoluto che tanto speriamo al di là della vita. *Bona nocte Vade dormitum!* E mi pareva che nel silenzio della campagna l'aria mi andasse ripetendo quelle care parole! Addio carissimo Amico: tanti doveri alla signora Elisa e a tutta la famiglia. Ti accludo la lettera del Ghezzi con cui accetta la tua proposta: attendo ora di conoscere le tue intenzioni

aff.mo amico
C. Clericetti

**46.**

Milano 21 Nov.e 1881.

Caro Amico.

Perdonami se non t'ho scritto prima per ringraziarti della Magnifica Monografia su Roma che hai proprio voluto inviarmi. Sono stato assai occupato in tutti questi scorsi giorni, occupatissimo, non solo per la riapertura dei Corsi all'Istituto e alla Società d'Incoraggiamento,[178] ma per altre cose che mi hanno assorbito assai tempo. Quella Monografia è un lavoro assai pregevole per la ricchezza di notizie, di raffronti, di dati che contiene ed ho cominciato a leggerla. Ti ringrazio vivamente del preziosissimo dono: poi ti dirò che ho anche aspettato la giornata di oggi a scriverti per darti notizia dell'esito della Estrazione della Lotteria. L'esito fu completamente negativo e per te e per noi. Non abbiamo né il numero di ieri 2797 né quello d'oggi 2387. Già non si faceva gran conto sulla Lotteria: ma insomma sopra 1000 premii non azzeccharne *[sic!]* neppure uno! Speravo che oggi sarei stato un signoretto: rimango invece al verde e senza trovarmi in campagna, il che però è già un fatto singolare. Tutto qui è tornato in calma e della brillante Esposizione[179] non rimangono che le spolpate baracche. Il nostro Guido ha superato lodevolmente gli esami che gli mancavano ed è entrato regolarmente nel 2° Corso Preparatorio: questa è certo una consolazione assai poveretto! Ieri mattina dei 101 colpi di cannone sparati per la Regina,[180] uno me lo sono dedicato a me stesso perché era anche il mio Compleanno: purtroppo sono 46! D'ora in poi ne conterò uno ogni due anni, per brevità di calcolo. Non avendo ricevuto da te alcuna comunicazione relativamente all'Assistente pel Ceradini, debbo ritenere che il Guidi non sia stato eletto a Torino: il suo nome infatti non compariva nei Giornali che diedero l'elenco dei professori nominati il mese scorso nelle parecchie Università: ma forse gli daranno l'incarico semplice senza nominarlo Professore Straordinario. Abbiamo avuto tre settimane di cielo così sereno, di sole così splendido e di temperatura così mite che non avevamo nulla da invidiare a Roma. Quanto volentieri sarei corso in campagna: e difatti moltissime famiglia non sono rientrate ancora. Se torno un'altra volta al mondo, voglia che ci sia pronta una Casa di campagna, altrimenti torno indietro!
Spero che la signora Elisa sia tuttora in buona salute e le invio tanti doveri: la Cecilia è raffreddata da qualche giorno ed obbligata in casa: una cosa leggiera. Addio: ricevi un affettuoso saluto, un ringraziamento ed una stretta di mano.

aff.mo amico
C. Clericetti

---

[177] Alla nota 1 dell'articolo di A. Garovaglio "Il culto di Mitra – Il battesimo ed i battisteri" (*Giornale della società storica lombarda*, 1889, Vol. 6, Fasc. 1, p. 161) si legge: "Quest'urna dall'epoca che ve la fece trasportare Valperto, trovavasi nella chiesuola di S. Bartolomeo al Bosco presso Tradate, provincia di Como fino al dicembre dell'anno 1884. Avvisata dell'esistenza dell'interessante monumento dal solerte ed intelligente Prof. Ing. Celeste Clericetti, ora defunto, la Commissione Archeologica conservatrice dei Monumenti della provincia di Milano, questa aperse tosto trattative col proprietario del suddetto oratorio Giovanni Grazioli, che assenziente il R. Ministero della Pubblica Istruzione, poté farne l'acquisto per la somma di L. 400; ed ora è uno dei monumenti più importanti del nostro Museo Patrio d'Archeologia, o che maggiormente attira l'attenzione degli interessati in simili studi."
[178] Si veda la lettera n. 1.
[179] Si vedano le lettere precedenti.
[180] Margherita Maria Teresa Giovanna di Savoia, nata il 20 novembre 1851.

**47.**

Milano 26 Decembre 1881.

Caro Amico

Mille augurii i più sinceri e cordiali, a te e a tutta la tua degna famiglia: mille voti perché la vostra salute perduri eccellente. Buone feste, buon capo d'anno. Io speravo di rivederti entro il mese: era il mio principale desiderio nella eventualità che fossi chiamato a Roma pel Genio Civile. Invece da quel che sento, sembra che il Nuovo Regolamento approvato dalla Camera, abbia mutate le disposizioni per l'ammissione degli Allievi: sembra che non si faranno più esami. Ma noi da lontano, non siamo mai bene informati e però non lo so di certo. Ma se è così, se la scelta degli allievi dovrà farsi dietro le Classificazioni ottenute nei Corsi delle scuole, allora il Genio Civile sarà riempito di napoletani, dove i 10 sono più numerosi dei nostri 8. È un pezzo che il mezzo dì trionfa: ma ne parlavi anche tu e già due anni sono, ed è un serio guaio.
In famiglia siamo tutti in buona salute e la Cecilia ha superato se stessa, poiché da un mese circa non ha avuto alcun ritorno della sua periodica e violenta emicrania. Anche la stagione quest'anno è singolarmente mite: non ricordo un altro inverno così sereno e così tiepido. Tutto dunque invita all'allegria, se non ci fossero tanti pensieri, tante continue preoccupazioni e nessuna distrazione. Non so più nulla del Ghezzi: dopo il tuo telegramma in risposta al mio non l'ho più visto: gli ho detto allora tutto quello che doveva fare e credo che la sua partenza per costì sia imminente. Ma certo tu sai meglio di me: quello di cui posso assicurarti, si è che hai fatto un buon acquisto per la scuola. Ho fiducia che non avrai che a lodarti di lui e tanto più buona è la scelta, perché a differenza di molti altri assistenti, egli ha una decisa preferenza per la Carriera dell'insegnamento: e questa è la ragione principale per cui ha accettato per la tua scuola e rifiutato l'altro posto a Conegliano dove avrebbe dovuto insegnare cose estranee alle sue inclinazioni.
È assai grazioso il duetto Baccelli-Sbarbaro[181] e degno di figurare in qualche teatro diurno. Ma un altro duetto un po' troppo ripetuto in questi ultimi giorni è il Bismark-Leone:[182] Oh se il Papa volesse andarsene una buona volta davvero e per sempre: questi signori oltremontani provino dunque ad avere in casa quella delizia: noi pagheremo anche le spese ben volontieri. L'immenso disastro di Vienna[183] ha destato anche qui una sensazione immensa ed è certo che i teatri ne risentiranno nella entrante stagione di carnevale. È però singolare il fatto che nel vero paese dei teatri, l'Italia, e dove la trascuratezza e il lasciar andare le cose senza troppo occuparsene, sono note caratteristiche, non sia mai accaduto nulla di simile. Speriamo dunque ancora che non accada nulla: ma quando si pensa che la distruzione pel fuoco è la morte naturale dei teatri, quando si pensa alla Leggi della probabilità, si deve arguire che una volta o l'altra deve pur accadere un incendio. Alla Scala si sono prese in questi giorni molte precauzioni, ma quel vecchio edificio è così difettoso e mancante delle condizioni necessarie alla sicurezza e all'igiene! Il Brioschi è a Milano, ma dovrebbe esserci assai più di frequente o nominare un sostituto durante le sue assenze: perché con quel nostro Segretario[184] che vuol fare tutto lui, con quelle sue vedute d'una spanna e quel suo carattere intollerabile, o si è in continue liti o si sopporta ciò che non si dovrebbe. Addio. Tanti doveri alla signora Elisa e i migliori auguri di me della Cecilia e dei nostri figli. Accogli una stretta di mano e credimi sempre

aff.mo amico
C. Clericetti

---

[181] Guido Baccelli, ministro dell'Istruzione (dal 1880 al 1884) venne aspramente criticato da Pietro Sbarbaro - giornalista, scrittore, filosofo, professore universitario di diritto e politico - per aver espulso dall'Ateneo di Cagliari due studenti di fede repubblicana. Si veda, ad esempio: *Alberico Gentili, Atti del Convegno (S. Ginesio, 11-12-13 settembre 2008)*, vol. 2, Giuffrè ed., 2010, pp. 206-207.
[182] Si riferisce ai rapporti politici tra Bismarck e il papa Leone XIII.
[183] Il 23 marzo 1881, nel rogo che divampò nel Teatro dell'Opera di Vienna, persero la vita novantadue persone.
[184] Ing. Giuseppe Giovannini.

**48.**

Milano 12 Genn° 82.

Carissimo Amico.

Ieri mattina ho visitato il Ghezzi che è tuttora a letto e dimagrato in modo singolare. Ha sofferto una infiammazione intestinale, che fu trascurata a lungo: la malattia acuta è vinta ma gli rimane una prostrazione grande di forza e una tale debolezza di stomaco che non può digerire il cibo. Egli si mostra più che altro preoccupato del pensiero che tu voglia lasciarlo in libertà se pel 1° di febbraio non si troverà costì, mentre già da due mesi ha preso stanza in codesta città. Gli promisi di scriverti, ma prima ho voluto parlare anche col suo medico curante D$^r$ Belloni: ci andai, non era in casa. Ieri sera quando tornai dall'avergli parlato, ricevetti la cara tua. Il D$^r$ Belloni mi disse che attualmente il Ghezzi è in stato di convalescenza, che non ha nulla al piloro né ai polmoni, e che, altro non intervenendo a complicare le cose, crede che pel 15 Febbraio potrà partire ristabilito. Non so che volesse alludere colla frase "altro non intervenendo" ed è per questo che lo richiesi se il Ghezzi fosse minacciato di tisi polmonare o soffrisse al piloro, il che egli ebbe ad escludere. Le cose sono dunque in questo stato che il Ghezzi è dichiarato in convalescenza ma non può reggersi in piedi e che non si nutre che di brodo e di tapioca. Nonostante, è giovane, e i giovani fanno presto a rimettersi e di qui alla metà di Febbraio c'è ancora un mese. Ma io ho voluto dirti precisamente come stanno le cose: il Dottore è informato della destinazione del Ghezzi e del cruccio del giovane di dover tardare, e della tema sua che tu non possa soprassedere. Io vedrò ancora il Ghezzi fra pochi giorni e ti scriverò di nuovo: intanto ti prego di non voler prendere alcuna determinazione. Per parte mia ho assicurato il Ghezzi che tu lo avresti aspettato almeno fino alla metà di febbraio: povero giovane! Trovo bene di aggiungere avermi il Ghezzi assicurato di non aver mai sofferto prima malattia alcuna: dunque non è fiacco di costituzione o malaticcio di natura: ma il male lo ha così abbattuto che fa senso. Ripeto che entro una settimana ti scriverò nuovamente e voglio sperare che le cose si volgano al bene. Ho passato molti giorni di malessere io pure, senza però sospendere le mie occupazioni, ma con dolori al capo assai forti e molesti. La verità è che per star bene ho bisogno ogni tanto di cambiare per una giornata il genere delle occupazioni, di avere per esempio qualche incarico che mi obbliga a star fuori di città una giornata. Se l'opportunità si presenta, bastano quelle poche ore a guarirmi, altrimenti i mali si sommano. È per questo che avrei accettato volentieri di entrare nella Commissione Provinciale per la conservazione delle antichità patrie; tale Commissione si dimise lo scorso anno per disaccordi fra i membri, due quali [sic!] sono nominati dal nostro Consiglio Comunale, uno dal Ministero. Parecchi della Giunta mi proposero insieme ad un altro, ma il Consiglio volle riconfermare uno dei rinuncianti (il Boito) e nominare un pittore (il Bertini) per l'altro, escludendo solo l'altro rinunciante, il Massarani. Il Ministero deve ora nominare il proprio membro e se tu non fossi in rotta col Ministro, avrei pregato te di propormi: non si tratta di posto che dia lucro alcuno: anzi il contrario perché non pagano nemmeno le spese, ma che mi conviene un po' per simpatia all'argomento, poi perché appunto obbliga ogni tanto ad una visita in provincia, ad una scampagnata di poche ore, ossia ad una diversione.[185] Ti dirò che gli Assessori che mi proposero (A mia insaputa) furono: Cusani,[186] Tagliasacchi,[187] Vimercati,[188] nonché il Sindaco[189]: ma chi comanda qui è un certo Circolo di Casa d'Adda[190]. Spero che la tua signora siasi completamente ristabilita: la Cecilia è sempre fiacca ma insomma tira innanzi senza medici. I ragazzi poi o giovinotti stanno sempre bene. Io avrei proprio caro che tu leggessi quelle poche pagine sulle Esperienze di Wöhler per dirmi se non ho ragione di spiegarle nel modo che ho fatto e di concludere che gli sforzi repentini, ripetuti sono più pericolosi perché inducono allungamenti ed accorciamenti massimi doppi di quelli che si producono per un egual carico che cominci da zero e vada

---

[185] Clericetti fu dapprima membro e in seguito vice presidente della Commissione Conservatrice dei Monumenti ed Antichità della Provincia di Milano dal 1882 al 1886.
[186] Ing. Luigi Cusani.
[187] Ing, Gioachimo Tagliasacchi, assessore alla Commissione edilizia.
[188] Vimercati Gaetano, assessore a varie commissioni amministratrici del prestito civico.
[189] Giulio Belinzaghi.
[190] Potrebbe trattarsi della casa di Carlo D'Adda, senatore e membro del Consiglio comunale di Milano dal 1871 al 1886.

crescendo fino al valore finale.[191] Tanti doveri a tutta la tua famiglia. Ti ringrazio vivamente di tutte le informazioni che mi hai dato nella tua lettera del 28 Dicembre riguardo al Nuovo Ordinamento del Genio Civile: ne parlai poi al Brioschi che naturalmente è favorevole agli Esami e mi disse che il Perazzi[192] lo era pure e che avrebbe indotto i suoi amici a propugnare la medesima soluzione: tanto più la cosa è adottabile se si accetta il provvedimento che mi dici suggerito dallo stesso Baccarini[193] di pagare le spese ai giovani poveri. Addio a settimana ventura. Ti stringo con affetto la mano

Devot.mo
C. Clericetti

**49.**

Milano 17 Genn° 1882.

Carissimo Amico

Il giovane Ghezzi è tutto racconsolato per la cortesissima lettere che avesti la bontà di scrivergli e ti è sommamente grato del tuo buon cuore. Ma... ma io temo. Per ora ti dirò solo che non v'è alcun miglioramento: il suo stomaco rifiuta il cibo: passa cattive notti. Da quel poco che so di un consulto tenuto l'altro ieri, v'è effettivamente qualche cosa ad un polmone come appunto dubitavo io. Ci metteranno una pezza... forse: ma ho paura. Permettimi pel momento di non dirne altro: te ne scriverò ancora la prossima settimana. Non vorrei danneggiare la sua posizione presso di te, ma lo stato del Ghezzi mi rattrista profondamente ed io debbo essere sincero teco. A settimana ventura: ed intanto ti scriverò di altre cose. Io t'ho pregato perché volessi leggere quelle poche mie pagine sulle Esperienze di Wöhler ed ora te ne dico la ragione. In quel breve studio ho dato una nuova interpretazione ai noti numeri del Wöhler che mostrano l'influenza sulla rottura del ferro e dell'acciaio dovuta alla ripetizione degli sforzi.[194] E vi sostengo che l'influenza che dirò deleteria degli sforzi ripetuti, rispetto a quelle di sforzi crescenti per gradi infinitesimi da zero al valore finale dipende dalla legge d'elasticità rappresentata da una breve formula di cui parla il Poncelet, il Lamé ed altri, per la quale un carico repentino produce oscillazioni nella barra attorno alla posizione di equilibrio, e la più ampia rappresenta un allungamento doppio del definitivo ossia di quello corrispondente ad uno sforzo eguale, ma cresciuto per gradi infinitesimi. Avrei dunque dato una specie di dimostrazione teorica dei risultati di Wöhler e le formule che se ne deducono dovrebbero essere più attendibili delle altre finora proposte per tener conto della ripetizione degli sforzi. Ora la stessa formula conduce alla legge $S=1/2(t+S_o)$ vale a dire se lo sforzo permanente o minimo ($S_o$) della barra è nullo, si ha $S=1/2t$. Cioè <u>quel più piccolo sforzo specifico</u> (S) <u>che può produrre la rottura, quando sia ripetuta un numero illimitato di volte, è la metà di quello sforzo specifico</u> (t) <u>che può produrre immediatamente, applicato una sola volta</u>. Ora voglio esporti un singolare ravvicinamento di risultati che comincia a rischiarare per me un dubbio che ho da molto tempo. Alcune esperienze, specialmente parmi di Fairbairn hanno mostrato che il cimentare una barra di ferro fino alla rottura, ne aumenta il coefficiente di rottura.[195] Il Tresca poi ha messo in chiaro l'altro fatto che il cimentare i materiali con sforzi superiori al limite di elasticità, questo limite si sposta e può spingersi più presso la rottura. Ebbene non v'ha una analogia fra questi fenomeni della materia inorganica con quelli che presentano i nostri muscoli e i nervi che coll'esercizio si rinforzano e divengono più elastici e più capaci a sopportare sforzi maggiori? E le facoltà della nostra intelligenza la Memoria etc non partecipano delle stesse proprietà? Qui mi pare di vederti ridere, ma bada che queste cose non le ho stampate e sono pensieri che espongo a te. Non arriveremo ad

---

[191] C. Clericetti, "Sulla determinazione dei coefficienti di sforzo specifico, dietro la esperienze di Wöhler", *Il Politecnico - Giornale dell'ingegnere architetto civile ed industriale*, v. 13, ottobre-novembre 1881, pp. 548-563. Sullo stesso argomento si veda la lettera seguente.
[192] Potrebbe trattarsi di Costantino Perazzi, ingegnere e politico italiano.
[193] Alfredo Baccarini, già direttore dell'ufficio provinciale del Genio Civile a Grosseto, nel 1882 è Ministro dei Lavori Pubblici.
[194] Si veda la lettera precedente.
[195] Sir William Fairbairn scrive su tale argomento ad esempio in: *Useful information for engineers*, Longman, London, 1860.

ammettere che le leggi dell'elasticità debbano essere le medesime tanto in un organismo vitale quanto in uno bruto? Oppure non dovremo ammettere che vi sia una specie di vita primordiale anche nei materiali? Non ridere e lasciami finire. Pochi giorni sono lessi di alcune interessanti Esperienze istituite dal D$^r$ Paul Bert a Parigi che suppongo sia lo stesso Attuale Ministro dell'Istruzione in Francia. Al D$^r$ Bert premeva di determinare sperimentalmente i limiti delle dosi sotto le quali si possono somministrare le sostanze Anestetiche: cioè le minime dosi, producenti i primi sintomi dell'Anestesia e la massima producente la morte subitanea: le dosi comprese fra questi due limiti formano ciò che il Bert chiama "<u>Zone maniable des Anestetiques</u> [*sic!*]". Lo scopo pratico della ricerca è dunque evidente. Ecco che ti trascrivo i risultati di tali esperienze[196]

|  | Cane | | Sorcio | | Merlo | |
| --- | --- | --- | --- | --- | --- | --- |
|  | 1° limite | mortale | 1° limite | mortale | 1° limite | mortale |
| Cloroformio................... | 9. | 19 | 6. | 12 | 9. | 18 |
| Bromuro d'Ethyle........... | 22. | 45 | 7. | 14 | 15. | 30 |
| Etere............................. | 37. | 74 | 12. | 25 | 18. | 40. |

Ne conclude, naturalmente, che <u>la dose minima, quella che produce i primi sintomi, che anche ripetuta non conduce ad effetto letale, è la metà di quelle che produce la morte immediata: La dose che può produrre</u>[197] <u>la morte quando sia ripetuta illimitatamente, è la metà di quella che la produce a un tratto con una sola somministrazione</u>. Non è singolare questo ravvicinamento? Dovremo concludere che le Leggi della statica dovranno in avvenire entrare anche nel campo della Farmacopea? Perdonami se t'ho stancato ma dimmi cosa ne pensi. Tanti doveri alla tua gentile signora anche da parte di mia moglie. A settimana prossima. Accogli intanto una stretta affettuosa di mano

aff$^{mo}$ amico
C. Clericetti

*[Scritto a lato:]*
NB. Grazie dei 2 numeri del Diritto[198] che ricevo ora e che leggerò domani.

---

[196] Tale risultato viene riportato, per esempio, in E. Bérillon, *L'œuvre scientifique de Paul Bert*, Paris, Picard-Bernheim, 1887, p. 100.
[197] Probabilmente manca un "non": "La dose che <u>non</u> può produrre la morte...".
[198] Potrebbe trattarsi della rivista *Il Diritto*, pubblicata da Enrico Lai.

**50.**

Milano 6 Luglio 1882.

Carissimo Amico.

È un pezzo che non ricevo notizie dirette sulla salute della tua Signora, ma pur troppo so che non accenna ancora a miglioramento.[199] Ed è una cosa che ci rattrista assai e non ti scrivo da parecchio tempo per non seccarti con Argomenti estranei, e che ti debbono recare noia più che altro. Del resto anche la Cecilia non è mai guarita della Corizza che soffre da quattro mesi e il cambiare Medico e Metodo di cura non ha valso ancora a nulla. Non è obbligata al letto, ma da un mese non esce di casa, cosicché anche noi non abbiamo fatto ancora progetto alcuno per la vacanza che pure è imminente. E tu che aspiravi tanto alle Alpi e che ora sei trattenuto da così dolorosa occupazione! quante volte ci penso! Mio fratello Emilio è stato finalmente promosso Maggiore e rimane nei Bersaglieri ma nel 7° e non più nel 6° Reggimento. Trovasi ora colla moglie a Firenze. È questa l'unica novità nella nostra famigliola se pure non ne è un'altra, quella che il nostro Guido s'è fatto maggiorenne, avendo ieri compiuto i 21 anni. Voglio sperare che superi bene tutti gli Esami del 2° Corso preparatorio, come ha superati quelli del 1°. Del tuo Vittorio è un pezzo che non ho notizie, ma egli deve essere innanzi un anno del Guido e spero che farà sempre bene. L'Emilio è stato ancora esentato dagli esami e perciò trovasi in piena vacanza. Ha un anno ancora di Liceo e pare che voglia dedicarsi alla Medicina. Al Politecnico nulla di nuovo ma parecchi professori sono sconfortati perché nessuno pensa a nominarli Ordinarii ed altri perché non sono neppure Straordinarii, ma semplici incaricati dopo 9 e dieci anni di eccellente insegnamento. Sarà in parte colpa di speciali circostanze, ma in parte è colpa del Direttore[200] che non si occupa o non vuole occuparsi del personale dell'Istituto. Io so che i tuoi allievi sono riusciti i primi anche quest'anno nel Concorso ai posti del Genio Civile e me ne congratulo vivamente con te che in pochi anni hai saputo creare una scuola così seria e così completa. Quanto ai nostri Allievi è però proprio un fatto che da parecchi anni nessuno dei distinti o dei buoni aspira a posti nel Genio Civile: si collocano meglio qui nelle Industrie che si vanno sviluppando e raffermando ogni giorno di più. Già saprai che il povero Cornalia è morto anche lui,[201] ma Colombo ne ha avvantaggiato molto, avendo sua moglie ereditata metà della sostanza del Professore ed essendogli anche immediatamente succeduto nella Pensione all'Istituto Lombardo. Il povero Cornalia non aveva che 57 anni ed il Brioschi ha letto di lui al Cimitero una Commemorazione funebre bella e commovente.[202] Gli eredi della Signora Kramer vanno rapidamente scomparendo: il Bossi,[203] la Venino, il Cornalia[204] se ne sono già andati. Il Vannucci,[205] il Caccia pure non tarderanno a seguirli e la sostanza destinata alla Fondazione Kramer va rapidamente crescendo.[206] Così vanno le cose. Addio. Ricevi tanti saluti anche dalla Cecilia col desiderio vivissimo di migliori notizie della signora Elisa.

aff.mo amico
C. Clericetti

---

[199] Infatti Elisa morirà 10 giorni dopo, il 16 luglio.
[200] Francesco Brioschi.
[201] Emilio Cornalia era morto l'8 giugno.
[202] F. Brioschi, "Emilio Cornalia", *Programma Regio Istituto Tecnico superiore di Milano*, Tip. Galli e Raimondi Milano, 1882-83, pp 71-73.
[203] Benigno Bossi.
[204] Emilio Cornalia era figlio della nobile Luigia Kramer, parente di Carlo, marito di Teresa Berra Kramer.
[205] Atto Vannucci.
[206] La "Pia Fondazione Edoardo Kramer", istituzione di beneficenza e di assistenza con sede a Milano venne fondata nel 1871 da Teresa Berra Kramer ed è attiva ancora oggi. In particolare la signora Kramer aveva destinato uno speciale premio biennale al migliore studente del Politecnico che permetteva al vincitore di completare e perfezionare i propri studi con un viaggio d'istruzione nei paesi europei più avanzati.

**51.**

Milano 24 Aprile 1883.

Carissimo Amico.

È un pezzo che non scrivo, ma la mia apparente negligenza è proprio dovuta al molto da fare. Credo fermamente che non avrai <u>mai</u> occasione di pentirti di aver tenuto presso di te il Vittorio e che ne avresti avuta molta se avessi agito diversamente. Il nostro Emilio puoi immaginarti con quale entusiasmo ha visitato Roma per la prima volta: non gli pareva vero. Sotto una scorza ruvida egli nasconde un animo sensibilissimo. È stato lieto anche di aver veduto il Vittorio e l'Itala in casa tua, ma quasi più non li riconosceva, tanto sono cresciuti. Io ti prego di non volere dare addosso agli Ambrosiani per la bomba dell'Esposizione mondiale.[207] Non si sa come sia nata la matta idea di passar sopra all'esplicito desiderio dei romani, ma è nata fra pochi, forse in un lieto Simposio e dopo aver bevuto molto: forse, ed è probabile, era un maneggio per far passare il Progetto di fabbricare la Piazza d'Armi e la Piazza Castello: progetto odioso a quasi tutti i milanesi, ed a ragione. I Giornali, meno la Perseveranza,[208] hanno mostrato poco senno, gonfiando il pallone: ma ora la situazione è affatto cambiata. Il Negri, che non si capisce come diavolo sia entrato in quel Comitato Provvisorio, se n'è dimesso e s'è dismesso anche dalla Presidenza della Costituzionale, tanto è la opposizione che l'idea ha incontrato qui, e tanto è il desiderio di mostrare che nessuno vuole soverchiare Roma. Però la mattata ha avuto un buon risultato; ed è quello di aver spinto la Capitale a svegliarsi e a dar corpo formalmente al suo progetto. Adesso Roma non può più esitare. Però in fondo a tutto questo v'è qualche cosa di serio e di grave: si dice che i risultati della Mostra dell'81 sono già sfruttati, che le Commissioni avute dagli Industriali sono esaurite e che fu fuoco di paglia. Si dice che parecchie Ditte debbano presto lasciare in libertà operai, e che ogni tanto è necessario qualche cosa di Straordinario, cioè un altro fuoco di paglia, per tener in piedi le nostre officine, per dar lavoro alle Maestranze. Viviamo dunque a scosse elettriche e il maggior risultato del progresso, sarebbe quello di ridurre tutto al provvisorio e di edificare sulla sabbia mobile, nella quale, secondo Byron, sono scritte le promesse delle donne. È ora di tornare alla età della pietra, alle case lacustri su palafitte, alla pesca all'amo. Rivive ora il Progetto di passare lo stretto di Messina con un manufatto stabile ed a proposito di ciò, ti pregherei di un favore. Qualche anno fa è uscito un opuscolo, od una Memoria in un Giornale tedesco, nel quale era esposto con un disegno illustrativo, il Progetto di passare la Manica con un ponte a molte <u>arcate di acciaio</u>, ciascuna della corda di 1000. metri.[209] Il Giambastiani si attiene a tale idea per lo stretto di Messina,[210] ma la dà per sua, il che non è. Ma io non mi ricordo affatto né l'anno, né il giornale, né l'opuscolo in cui era esposto tale Progetto. Se per caso tu te ne ricordassi, ti sarei immensamente obbligato se vorrai avvertirmene <u>subito</u> con una <u>Cartolina</u> postale, perché le mie ricerche finora sono state infruttuose e nessuno degli amici che ho interpellato si ricorda l'epoca o il titolo, benché ammettano di aver visto il Progetto. La Cecilia è ora a Novara presso la famiglia Zanetti e ci rimarrà per una settimana e voglio sperare che le faccia bene, perché insomma è sempre malaticcia. Sono dolente che non mi si presenti l'opportunità di venire a Roma: avrei amato visitare l'Esposizione, benché se ne parli così poco. Addio: tanti saluti all'Itala e al Vittorio. Ricevi una caldissima stretta di mano

dal vecchio Amico
C. Clericetti

---

[207] Si tratta dell'Esposizione Internazionale di Belle Arti, tenuta a Roma dal 21 gennaio al 1° luglio 1883.
[208] Si veda la lettera 41.
[209] Sembra strano che ancora in questi anni si parli in un articolo di un ponte di collegamento sulla Manica, perché pare che già Napoleone avesse pensato a un tunnel e nell'81 era già iniziata la trivellazione sul lato inglese. Clericetti aveva trattato di ponti nell'articolo "Sopra i ponti americani e sulle più recenti fondazioni tubolari", *Il Politecnico-Giornale dell'ingegnere architetto civile ed industriale*, v. 5, 1873, pp. 340-351; 437-446.
[210] Nel *Catalogo della esposizione collettiva del Ministero dei Lavori Pubblici alla Esposizione nazionale di Torino del 1884* si trova la citazione di una pubblicazione di Angelo Giambastiani: *Progetto di un ponte in acciaio per l'attraversamento dello Stretto di Messina, redatto dalla Direzione tecnica governativa delle ferrovie Novara-Pino e Genova-Asti*, Genova, 1884. Segue la spiegazione: "Questa pubblicazione serve d'illustrazione a due grandi quadri esposti nella sala della Esposizione collettiva del Ministero dei Lavori Pubblici, il primo dei quali rappresenta il prospetto del detto ponte composto di 3 arcate della corda di 1000 metri ciascuna; e di 2 arcate estreme, ognuna della corda di 500 metri. Il secondo rappresenta le basi della montatura dei grandi arconi."

**52.**

Nervi 1° del 1885.

Egregio Amico

Sono qui da tre settimane in cerca di riposo e di salute in questo angolo tranquillo del paese, fra questa flora orientale e questa splendida marina. Della burrasca passata conservo un capogiro insistente, un senso di vertigine per cui mi par sempre di dover cadere. Non cado mai però e dacché sono qui colla Cecilia, vado attorno solo: ho un semplice congedo dall'Istituto di Milano e ritorneremo colà per la metà del mese, perché la casa è là e naturalmente qui si è sempre sossopra. Il fratello Pietro che ora non ha impegni di sorta è qui pure e vi passerà buona parte dell'Inverno che è crudo a Milano e invece mite e benigno a Nervi. A Milano riprenderò in parte le mie antiche occupazioni, i miei molti impegni: ma già il migliore Augurio che ti posso mandare pel nuovo anno, è quello della salute che la Cecilia ed io ti desideriamo col cuore ottima e prospera. Conserva quella confidenza nella vita che è indispensabile ad ogni serio intento. Ti sei creato attorno una nuova famiglia col matrimonio di tua figlia[211] ed un'altra ti sorgerà accanto in un prossimo Avvenire:[212] puoi dunque, data la salute, guardare francamente in faccia al futuro. Quanto a noi non avemmo mai ragazze, e dei due figli, il minore, Emilio, è qui con noi a passare pochi giorni, poi ritornerà ai suoi studi a Torino, ed il Maggiore, Guido, non ha potuto passare il Natale con noi. Fra cinque mesi sarà agli esami: poi se tutto riuscirà bene, come pare, sarà Tenente, messo in un Reggimento ed inviato forse sul Napoletano. Dice che passerà con noi il futuro Natale! Grazie... ma la famiglia è spezzata per sempre. È un bene che la Cecilia o non pensa troppo all'avvenire od ha maggiore coraggio, il che è probabilissimo, di me: essa è rassegnata e sta bene in salute: anzi meglio di una volta. Ma già: sarà l'ozio e l'inazione che mi danno questi pensieri: ed io che non sono abituato a questo far niente, ne soffro le conseguenze. Perché vado pensando a tante cose, alle quale non avrei tempo da dedicare se fossi come sempre, molto occupato di mente: perché, dato pure che la mia attività sia inutile, è divenuta una necessità per me, un fuoco come quello delle Vestali antiche, che guai se si spegne: perché infine l'illusione è una necessità della vita. Qui la cura della salute è la suprema Legge, per tutti quelli che ci passano l'inverno: ed in gran parte sono malati di petto: Nervi è una vera Stazione climatica. La Flora Orientale domina nei giardini come saprai, e l'aria è imbalsamata da un purissimo effluvio: abbiamo dodici gradi di temperatura, mentre a Milano saranno allo zero. Tutto è splendido qui quando il cielo è sereno: ma quando piove, com'è successo in molti giorni scorsi, non si può quasi uscire di casa pel fango. E che fango! Sulla Via che corre lungo Nervi e la Marina, il fango raggiunge l'ideale: imbratta tutto, copre tutto e anche per poco, si rimane inzaccherati fino al colletto e se ne trasportano dei Chilogrammi con se. Addio dunque, caro Amico: rinnovo i più caldi Augurii perché ti sia felice il nuovo anno e lo sia pure per tutta la tua gentile famiglia e mandiamo insieme i più caldi saluti

Aff.mo Amico

C. Clericetti

---

[211] Elena, si veda la lettera 42.
[212] La figlia Itala si sposerà in effetti il 22 novembre 1890 con Vincenzo Cozzolino; si vedano le lettere 052-12133 e 12134 dell'archivio.

**53.**

Milano 14 marzo 1885.

Caro Cremona

Sono diventato poltrone anche riguardo allo scriver lettere, e l'antica attività mia, perché almeno mi vorrei concedere questo merito, è finita. Ho abbandonato Nervi con mia moglie nella seconda metà di Gennaio: non ero ristabilito, ma tanto da riprendere in parte le mie Lezioni: ma istanze continue e di tanti, hanno fatto sì che dovessi proprio recarmi a Verona per l'affare del Forte Masua.[213] Ci andai con mia moglie, senza però presentarmi mai al Tribunale, perché munito di una giustificazione Medica: restai sempre all'albergo, e là potevano venire ad interrogarmi sulla mia Relazione che avevo presentato mesi addietro. Tornai a casa e non ripresi altre Lezioni all'infuori di quelle al Politecnico: e cioè non ho date Lezioni alla Società di Incoraggim°, limitando la mia azione alla presenza personale. Mi pare di aver fatto tutto quello che potevo e anche di più, perché altri professori si sarebbero presi un più lungo riposo. Sto benino, ma conservo un senso di vertigine al capo, sempre meno intenso, ma pur sentito ancora. Conservo un senso di eccitabilità al cuore, che non avevo anzi prima della malattia; e ciò dipende dallo stato delle mia arterie, non più elastiche, ma ateromatose. Per questo poi non v'ha rimedio alcuno, e mi serva almeno di memento e di freno, ad ogni voglia di occuparmi troppo. Ho visto a Nervi il professor Rossetti[214] e la pur gentilissima signora: egli era a letto ed incomodato assai alla gola. Aveva voglia di fare tanto ma non poteva: mi parlò di te e dei colleghi Ferrini e Colombo di qui: faceva qualche passeggiata lungo quella stupenda marina, con un cappello di paglia a larghe tese. Figurati che là a Nervi si poteva far questo in Gennaio! Del resto se vuoi che definisca Nervi in poche parole, ti direi che è un borgo che sorge in mezzo al fango, tra il monte e il mare. In tutto questo tempo ho visto Brioschi una volta sola: fu al mio arrivo da Nervi, al ripresentarmi alla scuola. E quando mi disse che mi trovava pallido e dimagrato, gli osservai che lui pure non aveva la faccia d'un tempo. Rispose che per lui v'erano altre ragioni (e già tutti lo sanno). Ma lo trovai assai invecchiato, e la pelle del suo volto mi sembrava di cartapecora. L'affare del Chinino ha fatto senso anche a lui. Il Rag Maglione incaricato di liquidare la faccenda, chiamò a se tutto il Consiglio di Amministraz$^{ne}$ e quindi anche il Presid. Brioschi. Maglione voleva che si quotizzassero tutti di un tanto per acquetare gli affamati creditori: ma il Brioschi disse che lui non aveva più nulla, che soffriva personalmente, ma non poteva far altro.[215] Mi si dice poi che ora siasi dismesso dalla Presid$^a$ del Consiglio di Amministraz dei Cementi:[216] del resto io so che non l'avrebbero rieletto nella prossima adunanza. Ecco un uomo che pare finito: lo è certo almeno per gli affari e lo sarà forse per altre cose. Ai suoi Professori, Brioschi non ha mai pensato perché egoista fino al midollo. Ha avuto mille aspirazioni e soddisfatte, predicando sempre agli altri l'abnegazione: fortunato, se intorno al suo nome si fa ora silenzio. Tu sarai forse ancora occupato del tuo Progetto di riforma dell'Insegnamento Superiore che leggeremo nei resoconti del Senato:[217] ma sarai certo occupato d'altro e come sempre, perché la scienza nel caso tuo, non va dimenticata. Ho visto a Genova il Prof$^e$ Tardy che è anzi venuto a trovarci a Nervi, e che ti rammenta sempre. Ha una bella abitazione dove vive solo in mezzo a una bellissima raccolta di libri. S'è ritirato dall'insegnamento e perciò non ha occupazione alcuna ufficiale. Simpatico come sempre per la molta cultura, pel carattere; rimpiange ancora amaramente la moglie

---

[213] Il Forte Masua è una fortezza costruita dal Genio Militare del Regio Esercito Italiano tra il 1880 ed il 1885 e sorge in località Masua nel comune di Fumane (VR) sui Monti Lessini. Durante la costruzione avvenne un incidente: il 2 novembre 1883 morirono quattro muratori. Se ne trattò anche in Parlamento e probabilmente al Clericetti venne chiesta una relazione sulla sicurezza del cantiere.

[214] Potrebbe trattarsi di Francesco Rossetti.

[215] Si riferisce allo scandalo che aveva travolto la società anonima "Fabbrica lombarda di prodotti chimici", che controllava pressoché da sola la domanda di chinino, del cui Consiglio di Amministrazione Brioschi era presidente. Il rag. Giovanni Maglione era curatore provvisorio del fallimento.

[216] Si tratta della "Società Italiana dei cementi e delle calci idrauliche", costituitasi nel 1865 con il nome di "Società bergamasca" con sede a Bergamo; nel 1927 diventerà "Italcementi". Brioschi era presidente del suo consiglio di amministrazione dal 1876 e si dimette effettivamente nel 1885.

[217] Si vedano i *Discorsi del Senatore Luigi Cremona pronunziati in Senato nelle tornate del 30 novembre, 14, 15, 16, 17, 18 dicembre 1886 20, 21, 22, 24, 25 gennaio 1887*, Forzani e C., tipografi del senato, Roma, 1887. Un paio di bozze di lavoro relative a tale progetto di riforma si trovano alla segnatura 052-12002 (8689).

estinta. E a solo nominarla, gli vengono le lagrime agli occhi e si fa taciturno. Abbiamo visto anche il Prof$^e$ M$^{se}$ Piuma e la signora: ma tu lo conoscerai meglio di me. Pieno di denari, aveva bisogno di una informazione riguardo Milano, assai utile per lui e ho potuto dargliela efficace. Ma t'ho forse stancato con tante cose e tante parole: ti mando dunque un vivissimo saluto e una stretta cordiale di mano

aff$^{mo}$ Amico
C. Clericetti

**54.**

Civello 27 Agosto 1885.

Caro Amico

Prendo le cose in tempo per assicurare la vostra venuta qui, che è tanto desiderata dalla Cecilia e da me. Ti aspettiamo, insieme alla Signora Itala, col desiderio che restiate qui molti giorni con noi, che saremo lietissimi di avervi in nostra compagnia: e, bando alle reticenze: dovete venire e restare qui, ammenoché ben inteso, non vi piaccia il luogo o la poca compagnia dei tuoi più fidi Amici. E solamente vorrai avere la bontà di scrivermi qualche giorno prima (Maccio per Civello Prov Como) perché allora si troverà una Carrozzella a Portichetto o a Fino per prendervi. Si prende a Milano la linea Milano-Saronno-Como (Stazione in Piazza Castello). Ci sono molte partenze: alcuni treni s'arrestano a Portichetto e a Fino, come i treni delle 8.45 a 11.26 a – 6.32 – 8.57 (da Milano). Altri treni si arrestano solo a Fino, come 8.18 a – 11.9 a – 2.34 p – 5.10 p. La distanza da Fino a Portichetto è brevissima, ma desidero solo di conoscere il treno con cui partirai, per far trovare il biroccio nell'uno o nell'altro punto. Il 12 settembre che è l'anniversario del nostro Matrimonio, non è che un pretesto per avervi qui: ma ricordiamo la promessa, e ricordiamo che ne hai preso nota fin da anni fa. Dunque il pretesto o come vuoi, è eccellente per noi. E la promessa fatta da te, non può essere scritta sulla sabbia come il Byron voleva che fossero le promesse delle donne. (Mi sento come se le mie orecchie venissero tirate da mano femminile). In mezza ora al più si viene da Portichetto o da Fino a questo Civello, dove l'aria è eccellente, il panorama che si gode dalla casa, assai buono, benché ci troviamo a poca altezza: poi la ferrovia nuova per Varese e che viene da Como, ci passa al piede, e sarà aperta fra 15 giorni. Quanto alla casa, vorrete accontentarvi, spero. Stasera vado a Milano coll'Emilio, e domani egli avrà la visita militare come compreso nella leva di quest'anno. Ma già è esonerato dal servizio perché ha l'altro fratello militare. A Milano la giunta, ha fatto una povera figura coll'aver distribuito per premio delle scuole Comunali, un libro tradotto dalla Segur che dice corna dell'Italia. Si vede quanta trascuratezza domini qui da noi: tutti vogliono comandare, tutti essere al sommo, e nessuno si occupa dei dettagli. Ed è ben certo che l'Assessore della Pubblica Istruzione non ha neanche veduti i libri che faceva poi distribuire fra gli Allievi. L'ex-Sindaco Belinzaghi, deve proprio, in questa occasione, sfregarsi le mani: lui che si vantava quasi di essere indotto, perché la dottrina gli pareva così inutile (almeno pel Sindacato). Il piano regolatore dorme ancora, e si vede anche dai ciechi, quanto sia facile il dire e difficile il fare. Intanto continua la bazza dei costruttori: vi si fabbrica per fas e per nefas: vi si ammonticchiano i piani delle case, gli uni sugli altri e si trattano così gli inquilini come Sardelle da disporre a piani nelle botti. Il Milano Insubre, 2500 anni sono, non poteva essere molto diverso per gli scontorcimenti delle vie, la negazione delle Piazze, etc del Milano attuale fuori delle mura. Se non ci fossero gli alti Camini industriali. Bel monumento! Ma non voglio essere pessimista, e quando s'è in campagna è difficile l'esserlo. Ma anche qui però, questa incessante distruzione di boschi e di selve, tolgono la più bella attrattiva alla natura che va avvicinandosi all'aspetto d'una carta geografica: vuol dire che in avvenire si fabbricheranno delle carte molto grandi e vi si passeggierà *[sic!]* al di sopra (in città) e sarà lo stesso come il trovarsi in campagna. Mi viene in mente a questo proposito di utili sostituti, una certa scena, mi pare del Ferreol,[218] in cui un personaggio sviene, e una donna per farlo rinvenire gli va iniettando aria nelle nari col soffietto. E lui, ritornando in se dice "Quel bonheur! Il me semble d'être à la campagne![219] E qui ti lascio, anche per tema di uno scappellotto, e nella

---

[218] *Ferréol*, commedia in 4 atti di Victorien Sardou.
[219] Clericetti non chiude la citazione con le virgolette.

sicurezza di avervi qui entrambi almeno perché il Vittorio sarà forse a Roma, mando tanti saluti della Cecilia e insieme una stretta cordialissima del tuo

Vecchio Amico
C. Clericetti

**55.**

Milano 14 Luglio 1886.

Carissimo Amico

Non ho notizie della tua salute ma la spero eccellente tuttora. A questi giorni penserai forse alle ridenti Valli Alpine, ai Boschi di Abete dall'aria profumata etc, e avrai in progetto qualche gita interessante su quelle stupende alture. Spero dunque di rivederti presto, perché già in questo caso probabilissimo, vorrai pure passare da Milano e stringere la mano ai tuoi vecchi amici. Intanto la notizia più recente che posso darti è la morte di un vecchio Paladino dell'Istituto, cioè del Biondelli,[220] che abbiamo accompagnato al Cimitero: lottò molto colla morte, parve riaversi completamente, ma dovette cadere anche lui. Così, scomparsi Lombardini, Poli, Cornalia, Hajeck, Biondelli, etc l'ambiente della sala dell'Istituto, è cambiato in pochi anni. Ma anche all'infuori di quella ristretta cerchia di pensatori, anche l'ambiente sociale si va mutando: ieri mattina è morta a 72 anni anche la Contessa Maffei e quelle sale che per 50 anni furono il centro di una eletta società di letterati filosofi e politici, sono chiuse anch'esse per sempre. Il Poeta Maffei diede il nome, ma non il pensiero dominante a quella storica riunione, perché lui fu austriacante (e lo si dice ora), ma la Maffei fu audacemente italiana anche nei primordii della sua gentile dominazione. Sono notizie tristi. Quest'anno non si fa la gita cogli Allievi del 3° anno, perché troppo postergata, così da incontrare l'indifferenza degli allievi, già in preda al panico degli esami imminenti. Ma ci sono tante cagioni del fatto e queste riguardano il modo di condurre l'Istituto e vi è pure quella mancanza di entusiasmo che pur troppo, per l'età crescente, è andato scemando nei Professori più anziani, fra i quali vi sono anch'io. La Cecilia sta benino, ma ogni tanto, è a letto col suo male antico: quella fortissima emicrania, refrattaria a qualsiasi rimedio che l'abbatte così profondamente. L'Emilio tornerà oggi da Torino, ove ha compiuto assai lodevolmente il 3° anno di Medicina, mentre il Guido ha ancora un esame prima di compiere ed essere licenziato da quella scuola di Applicazione di Artiglieria e Genio. Mio fratello Pietro è pure tornato dopo una lunga peregrinazione. Ha visitato il Cairo e le Piramidi: poi le nostre colonie di Assab, di Massaua. Fu ad Adan e Suakim e rifacendo il mar Rosso, andò in Palestina, poi a Costantinopoli e quindi in Atene. Fu un lungo viaggio e ritornò in eccellente salute ricco di notizie e di emozioni. Mio fratello Emilio che è a Firenze colla moglie doveva pure recarsi a Massaua, ma non ha potuto andarvi, perché trattenuto invece a Firenze per un Processo giuridico nel quale era testimonio importante. Avrai visto ultimamente il Colombo, ma io non lo so, perché non l'ho visto che per pochi momenti all'adunanza per gli Esami e pel viaggio che andò a monte. Ma anche il Nazzani tuo dev'essere Deputato e davvero non so come, e lui e il Colombo possano utilmente attendere alla Cattedra: ma adesso si concilia tutto, meno lo spero, il Papa e il Re, o piuttosto la Chiesa collo stato. Noi andiamo ancora a Civello e desidero proprio che quest'anno almeno tu possa trovare il tempo di venire a trovarci colla tua figlia, e stare un po' di tempo con noi. In tale desiderio, ti mando un affettuosissimo saluto e insieme tanti doveri alla signora Itala, e saluti al tuo giovine Ingegnere.[221] Ti prego di non dimenticarci, mentre mi dichiaro

Il Vecchio Amico devotissimo
C. Clericetti

---

[220] Bernardino Biondelli.
[221] Vittorio.

**56.**

Civello 3 8^(bre) 86.

Carissimo Amico.

Anche per quest'anno la vacanza è sul finire: Agosto e settembre sono passati con una fretta che a me pare maggiore della solita.

Avremmo tanto amato una vostra visita, ma l'ultima tua lettera me ne tolse la speranza. Ma tu sarai certo preoccupato di quel lieto evento domestico che mi hai annunciato: cioè il prossimo matrimonio del tuo giovine Ingegnere.[222] Da quanto me ne dicesti, sembra un eccellente partito e la Cecilia ed io mandiamo i più lieti augurii ai promessi Sposi. A pensare che ho veduto i tuoi figli poco più di bambini, ed ora due hanno già formato una nuova famiglia e tu sei già nonno felicissimo, non mi par vero che tanto tempo sia scorso per chi pensa ora come pensava allora. Il Guido è ora Tenente di Artiglieria e destinato a Milano: così lo avremo in casa per tutto l'inverno prossimo, e l'Emilio ha compiuto il 3° Anno di Medicina. La Cecilia sta bene in salute, anzi meglio di anni sono. Quella incertezza che avevamo riguardo all'esito degli esami del Guido, è cessata dunque, ed è subentrata all'ansia negli animi nostri, una calma che non può se non giovare anche igienicamente. Ma io sento che non potrei più dare lezioni serali e tenere quelle lezioni di disegno alla Società di Incoragg^(mo) [223] alle quali attesi per 23 anni: ed ho data la mia dimissione. Ma pur troppo la società non dà pensione affatto. Così pure mi sono dimesso da tutti gli impegni che avevo riguardo la Scuola dei Capomastri,[224] la quale (nessuno certo se ne ricorda ora) fu pure una creatura mia. Ma da un pezzo ero stanco di tutto questo lavoro, in gran parte serale, dovetti ascoltare di più il mio stato di salute, e dimettermi. Rimango professore al politecnico esclusivamente ed ecco tutto. A proposito del politecnico, anzi delle Scuole di Applicaz^e in genere, so che si tratta di aggiungerle alle Università, credo anzi per tua iniziativa, e creare per esse una facoltà speciale. L'intento a me pare giusto perché non ho mai approvato la tendenza delle Scuole di Ingegneri, specialmente svizzere, di avvicinarsi troppo all'Officina, come del resto, si voleva anni sono fare anche a Milano secondo l'idea del Colombo. Ho saputo però che la Scuola di Milano sarebbe esclusa da questa misura. Vedo che tu lavori e studii sempre, e mi felicito con te. A giorni poi ci vedremo, se pure pel 10 sarai già di ritorno a Roma dalla tua villa di Porto Maurizio e mi spiegherai meglio le innovazioni che si vanno preparando per opera tua circa l'assetto delle Scuole di Applicazioni. In attesa dunque di stringerti la mano personalmente, fa tanti doveri e saluti per noi alla signora Itala e al Vittorio e credimi sempre

aff^(mo) Amico
C. Clericetti

---

[222] Vittorio.
[223] Si veda la lettera n. 1.
[224] Si veda la lettera n. 1.

**57.**

Milano 3 del 1887.

Caro Amico.

Sono convalescente d'una grave pleurite che mi obbligò al letto fin dal 4 Novembre: ma non sono guarito affatto né potrò ancora uscire per un pezzo. Tutto rimase sospeso e non mi resta che di poter dire, guarisco. Non è una consolazione alla mia età, non è un conforto al dolore alla noia patita, allo svincolo da tutto quello che formava il mio mondo intellettuale. Già ho colto ai primi di Ottobre a Trezzo il mio malanno ed era in corso quando venni costì alla metà di quel mese. Tutta la filosofia e il conforto sono ridotte per me alla parola <u>pazienza</u>. Tu pure devi aver sofferto un dolore profondo per lo scioglimento del matrimonio progettato di tuo figlio. Ma non ne conosco le ragioni e queste potrebbero anche esser tali da bilanciarli *[sic!]*.[225] Ti resta certo però il dispiacere della partenza di Vittorio e della sua assenza per due anni. Il motivo per tale decisione deve essere stato forte, per abbandonare un impiego e occupazioni avviate: ma è però vero, che considerato il fatto in se solo, un viaggio all'estero è divenuto un completamento necessario all'educazione tecnica o letteraria d'un giovane: così diceva anche la Cecilia al ricevere la lettera della signora Itala. Tutto il male non viene per nuocere (ma io personalmente non posso proprio dirlo). Un'altra ragione di conforto la troverai nella tua riconciliazione col Brioschi, che mi venne asserita da qualcuno al Politecnico, prima ancora che mi ammalassi. È bene che due uomini eminenti non si osteggino all'infinito: ma è certo che tu specialmente avrai dovuto coprire d'un velo di carattere politico, il movente di certe azioni. Ma io non voglio risollevare nulla e applaudo alla vostra riconciliazione, e faccio voti perché sia perenne. Il Brioschi è venuto a trovarmi nei giorni più pericolosi della mia malattia: ma già non era il caso di dire una parola del fatto. Il Loria mi supplisce all'Istituto, in attesa della mia guarigione completa, e il De Marchi continua da assistente. Desidero sapere le ragioni, se può dirle della rottura del matrimonio, perché era già in anticipazione, una gioia per te il pensarci. Sai la decisione presa ad enorme maggioranza dal nostro Consiglio Comunale, riguardo al monum° di Napol. 3°. Ma io dubito assai che quel monum° verrà collocato in posto: ad ogni modo si tratta dell'avvenire e gli avvenimenti non aspettano il piano regolatore. Ma il partito del <u>secolo</u> è furibondo, e i quattro Deputati radicali di Milano, eletti dalla maggioranza degli Elettori politici, aspettano la prima occasione per rovesciare il Negri ed il Consiglio abborrito. Sarebbe una commedia se non scapitasse almeno il decoro della Città. Ma con tanta smania di industrie di officine, di operai, l'elemento radicale, demagogo, cresce di continuo, è una marea a flusso continuo, non riflusso. Adesso aspetto tue notizie che speravo anzi di ricevere prima, e m'avrebbero consolato tanto nella depressione di spirito in cui mi trovo da un pezzo. La Cecilia invia tanti doveri a te e saluti alla signora Itala. Accetta una cordialissima stretta di mano dal tuo vecchio amico.

C. Clericetti

---

[225] Le ragioni si trovano nella lettera che Cremona scriverà a T.A. Hirst il 5 giugno 1887: "La rottura dell'impegno matrimoniale è derivata dalla scoperta spiacevolissima che le condizioni morali della famiglia della sposa non erano quali in buona fede Vittorio ed io avevamo creduto. L'unione delle due famiglie diveniva impossibile, senza grossi guaj in avvenire. Vittorio, messo sul bivio tra il distaccarsi dalla sua propria famiglia e il rinunciare alla sposa, ha adottato il secondo partito, sotto la condizione, da me accettata, di potersi allontanare dall'Europa per alcuni anni. Egli partì da Genova il 18 marzo e sbarcò a Buenos Aires il 18 aprile. È probabile ch'egli abbia là una buona e importante occupazione nella sistemazione idraulica d'una provincia occidentale, sotto le Ande." Si veda: L. Nurzia (a cura di), "La corrispondenza di Luigi Cremona (1830-1903)", v. IV, *Per l'archivio della corrispondenza dei matematici italiani*, Quaderni P.RI.ST.EM., Palermo, 1997-8, p. 212.

**58.**

Milano 26 Febbr⁰ 87.

Caro Amico.

Molte vicende sono trascorse dall'ultima volta che ti scrissi dal letto, tormentato da una lunga e dolorosa pleurite. Non mi è venuto però alcun cenno tuo di ricevuta. Abbiamo visto il Vittorio come sai, quasi sulle mosse per partire per l'America. Ma io non so a che punto sia la cosa. Quattro mesi continui di sofferenza: adesso il Prof Loria, funziona da mio supplente al Politecnico e continuerà tuttavia fino a che potrò dirmi guarito. Ma ciò che soffro attualmente non è più la pleurite, perché è seppellita da un paio di mesi: ma in conseguenza di una dilatazione dell'intestino, venuta anch'essa con tante sofferenze. Ma è ora *[...]* di sentire la tua voce perlomeno dallo scritto: adesso poi non *[devi]* preoccupazioni ministeriali. Ma è strano che non si possa finire a formare un Ministero e che il Robilant[226] si comporta da vera mummia e che il Farini al solito non accetti il che lascia supporre o almeno lo lascia a me qualche Società Segreta.[227] Che farne di Uomini come il Farini, così inutili dal lato Politico, così dannosi, perché si ritira sempre tutte le volte che la patria ha bisogno di lui. E che disgrazia che Minghetti sia morto adesso, adesso che non v'è più destra e sinistra. Non so quando, ma mi pare che la patria non si è trovata mai in circostanze così impossibili: la soluzione non si intravede ancora, e bisognerà purtroppo che si esponga la corona. Vedremo Crispi al Ministero, e ci sia per amor di Dio.[228] Ci sarai forse tu, ma veramente non so proprio dire come ti trovi con Crispi. Ci sarà il Saracco probabilmente, ma quel benedetto Luzzatti distinto economista, non è uomo politico, poi fu in predicato per l'istruzione: e allora lo potrebbe essere per gli interni per gli esteri e per qualunque altro Ministero.[229] La Cecilia ha passato bene la crisi di quest'inverno: un inverno così rigido e lungo: ha avuto da fare con me per quattro mesi di seguito. Il nostro Guido che è sempre qui è contento perché com'è proprio vero, i soldati hanno fatto il loro dovere a Dagoli *[sic!]*,[230] qualunque sia per essere la soluzione *[...]* della quistione. Addio caro Amico: ti ripeto che desiderai a lungo tue notizie: e ormai sono vecchie di parecchi mesi. Spero che hai saputo conservare la salute negli scorsi mesi. Dimmi se il Vittorio è ancora costì o sia già partito e dove si trovi.[231] Rimango sempre, col desiderio vivissimo di tue notizie, tuo vecchio Amico

C. Clericetti

---

[226] Nicolis Carlo Felice conte di Robilant.

[227] Si veda la lettera 43.

[228] In aprile Crispi diventa Ministro degli Interni.

[229] In effetti, nel tempo ricoprirà le cariche di ministro del tesoro, delle finanze, delle poste e telegrafi, dell'agricoltura, industria e commercio e dell'interno.

[230] Probabilmente si riferisce alla battaglia del 26 gennaio 1887 a Dogali, in Eritrea, dove morirono cinquecento soldati italiani.

[231] Si veda la nota alla lettera precedente.

## Tabella - dati delle lettere

| LUOGO E DATA | MITTENTE | SEGNATURA FILE E CARTACEA | CONSISTENZA (n. di pagine) | LINGUA |
|---|---|---|---|---|
| **1.** Roma 30 aprile 1871 | Clericetti | 089-20878 (17545) | 4 | italiano |
| **2.** s.l. 27 luglio 1872 | Clericetti | 052-12072a (8759) | 1 | italiano |
| **3.** Tremezzo 8 agosto 1872 | Clericetti | 089-20879 (17546) | 3 | italiano |
| **4.** Tremezzo 9 ottobre 1872 | Clericetti | 089-20880 (17547) | 2 | italiano |
| **5.** s.l. 24 maggio 1873 | Clericetti | 089-20882 (17549) | 2 | italiano |
| **6.** Lucino 3 settembre 1873 | Clericetti | 089-20883 (17550) | 4 | italiano |
| **7.** Lucino 13 settembre 1873 | Clericetti | 051-11207 (7895) | 4 | italiano |
| **8.** Milano 17 ottobre 1873 | Clericetti | 089-20881 (17548) | 3 | italiano |
| **9.** Lucino 18 ottobre 1873 | Clericetti | 089-20884 (17551) | 4 | italiano |
| **10.** Milano 26 ottobre 1873 | Clericetti | 089-20885 (17552) | 4 | italiano |
| **11.** Milano 6 dicembre 1873 | Clericetti | 089-20886 (17553) | 4 | italiano |
| **12.** Milano 28 dicembre 1873 | Clericetti | 089-20887 (17554) | 3 | italiano |
| **13.** Milano 27 gennaio 1874 | Clericetti | 089-20888 (17555) | 4 | italiano |
| **14.** Milano 7 agosto 1874 | Clericetti | 089-20889 (17556) | 4 | italiano |
| **15.** Lucino 9 ottobre 1874 | Clericetti | 089-20890 (17557) | 4 | italiano |
| **16.** Milano 19 novembre 1874 | Clericetti | 089-20891 (17558) | 4 | italiano |
| **17.** Milano 20 marzo 1875 | Clericetti | 089-20892 (17559) | 3 | italiano |
| **18.** Milano 21 luglio 1875 | Clericetti | 089-20893 (17560) | 4 | italiano |
| **19.** Milano 20 luglio 1877 | Clericetti | 051-11201 (7889) | 3 | italiano |
| **20.** Monticello 17 settembre 1877 | Clericetti | 051-11200 (7888) | 4 | italiano |
| **21.** Milano 13 giugno 1878 | Clericetti | 051-11202 (78909 | 2 | italiano |
| **22.** Milano 20 dicembre 1878 | Clericetti | 051-11203 (7891) | 4 | italiano |
| **23.** Milano 12 maggio 1879 | Clericetti | 051-11204 (7892) | 3 | italiano |
| **24.** s.l. 20 maggio 1879 | Cremona | 052-11923 (8612) | 4 | inglese |
| **25**. s.l. 5 giugno 1879 | Cremona | 052-11924 (8611) | 4 | inglese |
| **26.** s.l. 15 giugno 1879 | Cremona | 052-11925 (8612) | 4 | inglese |
| **27.** Milano 24 giugno 1879 | Clericetti | 051-11205 (7893) | 4 | inglese |
| **28.** s.l. 5 luglio 1879 | Cremona | 052-11943 (8630) | 3 | inglese |
| **29.** Milano 30 luglio 1879 | Clericetti | 051-11206 (7894) | 4 | inglese |
| **30.** Milano 22 gennaio 1880 | Clericetti | 051-11208 (7896) | 4 | italiano |
| **31.** s.l. 26 gennaio 1880 | Cremona | 052-11926 (8613) | 2 | inglese |
| **32.** Milano 14 febbraio 1880 | Clericetti | 051-11209 (7897) | 4 | inglese |
| **33.** Milano 27 febbraio 1880 | Clericetti | 051-11210 (7898) | 3 | italiano |
| **34.** s.l. 3 marzo 1880 | Cremona | 052-11927 (8614) | 3 | inglese |
| **35.** Roma, 2 aprile 1880 | Cremona | 052-11928 (8615) | 4 | inglese |
| **36.** s.l. 4 giugno 1880 | Cremona | 052-11929 (8616) | 2 | inglese |
| **37.** Milano 8 giugno 1880 | Clericetti | 051-11211 (7899) | 4 | inglese |
| **38.** Milano 18 giugno 1880 | Clericetti | 051-11212 (7900) | 3 | inglese |
| **39.** s.l. 17 novembre 1880 | Cremona | 052-11930 (8617) | 3 | inglese |
| **40.** s.l. e s.d. | Cremona | 052-11931 (8618) | 1 | inglese |
| **41.** Milano 16 marzo 1881 | Clericetti | 051-11213 (7901) | 3 | italiano |
| **42.** Milano 28 marzo 1881 | Clericetti | 051-11214 (7902) | 4 | italiano |
| **43.** Milano 25 maggio 1881 | Clericetti | 051-11215 (7903) | 4 | italiano |
| **44.** Milano 25 luglio 1881 | Clericetti | 051-11216 (7904) | 4 | inglese |
| **45.** Milano 29 ottobre 1881 | Clericetti | 051-11217 (7905) | 4 | italiano |
| **46.** Milano 21 novembre 1881 | Clericetti | 051-11218 (7906) | 3 | italiano |
| **47.** Milano 26 dicembre 1881 | Clericetti | 051-11219 (7907) | 4 | italiano |
| **48.** Milano 12 gennaio 1882 | Clericetti | 051-11221 (7909) | 4 | italiano |
| **49.** Milano 17 gennaio 1882 | Clericetti | 051-11220 (7908) | 4 | italiano |
| **50.** Milano 6 luglio 1882 | Clericetti | 051-11222 (7910) | 4 | italiano |
| **51.** Milano 24 aprile 1883 | Clericetti | 051-11223 (7911) | 4 | italiano |
| **52.** Nervi 1 gennaio 1885 | Clericetti | 051-11224 (7912) | 4 | italiano |
| **53.** Milano 14 marzo 1885 | Clericetti | 051- 11225 (7913) | 4 | italiano |
| **54.** Civello 27 agosto 1885 | Clericetti | 051-11226 (7914) | 4 | italiano |
| **55.** Milano 14 luglio 1886 | Clericetti | 051-11227 (7915) | 4 | italiano |
| **56.** Civello 3 ottobre 1886 | Clericetti | 051-11228 (7916) | 4 | italiano |
| **57.** Milano 3 gennaio 1887 | Clericetti | 089-20895 (17562) | 4 | italiano |
| **58.** Milano 26 febbraio 1887 | Clericetti | 089-20894 (17561) | 4 | italiano |

## Indice dei nomi citati nelle lettere

**Antonelli** Alessandro (Ghemme 1798 - Torino 1888). Dopo aver frequentato a Milano l'Accademia di Brera, si laureò ingegnere architetto nel 1824. Si occupò della risistemazione del centro di Torino e della costruzione della Mole Antonelliana oltre che di numerosi edifici civili e religiosi. Si dedicò alla vita politica come consigliere nel comune di Torino e anche come deputato nel 1849. Si veda DBI (P. Portoghesi).
*Lettera n° 15.*

**Ascoli** Graziadio Isaia (Gorizia 1929 - Milano 1907). Nato nel contesto plurilingue della sua città da famiglia ebraica, fu linguista e glottologo autodidatta. Viene considerato il fondatore, in Italia, della linguistica comparata che insegnò presso l'Accademia scientifico-letteraria di Milano dal 1861. Fu il primo studioso a classificare i dialetti italiani. Nel 1864 divenne membro effettivo dell'Istituto Lombardo; nel 1875 entrò a far parte dell'Accademia Nazionale dei Lincei come socio nazionale per la classe di scienze morali e, dal 1882 al 1899, fu membro del Consiglio superiore della Pubblica Istruzione. Si veda DBI (T. Bolelli).
*Lettera n° 30.*

**Baccarini** Alfredo (Russi, Ravenna 1826-1890). Dopo essersi iscritto al corso di Matematica e Fisica presso l'Università di Bologna, partecipò come volontario alla difesa di Modena nel 1848. Si laureò in Ingegneria nel 1854 e dopo essersi occupato dello studio di alcuni tracciati ferroviari, si impegnò in politica ricoprendo diverse cariche comunali e provinciali. Nel 1871 divenne direttore dell'ufficio provinciale del Genio Civile a Grosseto. In quel periodo la Toscana era infestata dalla malaria e Baccarini relazionò la situazione al governo con la monografia *Sul compimento delle opere di bonificazione e sulla definitiva regolazione delle acque nelle Maremme Toscane*, che fu premiata con la medaglia d'oro all'esposizione di Vienna e in seguito con la convocazione a Roma presso il Consiglio superiore dei Lavori Pubblici. Nel 1876 venne eletto deputato e dal 1878 al 1883 fu Ministro dei Lavori Pubblici. Si veda DBI (G.P. Nitti) e la pagina dei personaggi illustri sul sito del Comune di Russi.
*Lettera n° 48.*

**Baccelli** Guido (Roma 1830-1916). Si laureò in Medicina (1852) e in Chirurgia (1853) presso l'Università romana "La Sapienza" e insegnò nell'ospedale di S. Spirito. La sua attività scientifica fu rivolta soprattutto alla clinica e agli aspetti terapeutici, in particolare della malaria e delle malattie cardiorespiratorie. Nel 1872 ebbe inizio la sua attività politica come presidente del Consiglio superiore di Sanità, carica che ricoprì per un totale di oltre trent'anni (1872-1877 e 1887-1915). Fu per ben sette volte Ministro della Pubblica Istruzione tra il 1881 e il 1900. Fu anche Ministro all'Agricoltura, Industria e Commercio tra il 1901 e il 1903. Si veda DBI (M. Crespi).
*Lettere n° 44, 47.*

**Basile** Giovanni Battista Filippo (Palermo 1825-1891). Dopo aver compiuto studi classici, si laureò in Architettura presso l'Università di Palermo. Si spostò poi a Roma dove seguì corsi di Matematica idraulica e Costruzioni presso l'Università "La Sapienza" e di Disegno presso l'Accademia di S. Luca. Durante il movimento rivoluzionario di Palermo del 1848 venne incaricato di costruire le fortificazioni e la polveriera del Sacramento; nel 1860 seguì Garibaldi alternando le attività patriottiche con lo studio dei monumenti antichi della Sicilia. Fra le sue numerose opere si ricorda il Teatro Massimo di Palermo iniziato nel 1875 e portato a termine dal figlio Ernesto Basile. Si veda DBI (M. Tafuri).
*Lettera n° 20.*

**Belinzaghi** Giulio (Milano 1818 - Cernobbio, Como 1892). Banchiere, nel 1849 fondò la Banca Belinzaghi. Fu eletto alla Camera dei Deputati nel 1867, ma rinunciò alla carica l'anno successivo in seguito all'elezione a Sindaco di Milano. Ricoprì tale carica fino al 1882 e poi dal 1889 alla morte. Fu presidente della Società ferroviaria mediterranea, della Camera di commercio di Milano e della Banca nazionale. Nel 1872 venne nominato Senatore. Si veda DBI (N. Foà).
*Lettere n° 48, 54.*

**Belloni** Giuseppe, medico chirurgo e ostetrico. Esercitò nella città di Milano almeno dal 1847 al 1889.
*Lettera n° 48.*

**Beltrami** Eugenio (Cremona 1835 - Roma 1900). Studiò Matematica a Pavia dal 1853 al 1856 senza peraltro laurearsi. Nel 1862 fu chiamato all'Università di Bologna come professore di Algebra e Geometria analitica. Nel 1864 si trasferì sulla cattedra di Geodesia, lasciata da Ottaviano Fabrizio Mossotti, presso l'Università di Pisa e in questa città rimase fino al 1866, anno in cui tornò all'Università di Bologna come professore di Meccanica razionale. Nel 1873 si trasferì all'Università di Roma sempre sulla cattedra di Meccanica razionale, ma questioni personali varie lo indussero, nel 1876, a tornare all'Università di Pavia sulla cattedra di Fisica matematica. Nel 1891 ritornò a Roma dove trascorse i suoi ultimi anni. Si vedano [Enea, 2009], [Giacardi-Tazzioli, 2012], DSB (D.J. Struik), DBI (N. Virgopia) e [Tricomi 1962].
Citato in [Cerroni-Fenaroli, 2007], [Canepa-Fenaroli, 2009], [Cerroni-Martini, 2009], [Enea-Gatto, 2009], [Giacardi-Tazzioli, 2013].
*Lettere n° 9, 33, 37, 41.*

**Berra Kramer** Teresa (Milano 1804-1879). Fu implicata nei moti milanesi del 1821 e molto attiva nella congiura delle donne lombardo-venete per l'Unità d'Italia. Dopo aver sposato Carlo Kramer, industriale nel ramo del tessile, si allontanò da Milano per Parigi e la Svizzera; tornò nel 1826 e vi rimase fino al 1851, dedicandosi all'educazione del figlio Edoardo e animando il suo salotto repubblicano. Fu sostenitrice di Mazzini. Nel 1871, dopo la morte del figlio, su consiglio dello stesso Mazzini, fondò la "Pia Fondazione Edoardo Kramer", attiva ancora oggi. Si vedano DBI (D. Di Porto) e www.duevoltiditeresa.altervista.org
*Lettere n° 14, 16, 37, 50.*

**Bert** Paul (Auxerre, Francia 1833 - Hanoi, Vietnam 1886). Iniziò i suoi studi all'École Polytechnique con l'intenzione di diventare ingegnere, ma li abbandonò per dedicarsi alla Giurisprudenza, materia nella quale ottenne anche un dottorato nel 1857. Si rivolse infine alla Medicina: si laureò nel 1864 e divenne professore di Fisiologia a Bordeaux nel 1866 e alla Sorbona nel 1869. Lavorò alla creazione di uno scafandro a regolazione di pressione e, nel 1870, si occupò dei gas anestetici. Fu Ministro della Pubblica Istruzione dal 1881 al 1882 e governatore generale del Tonchino e dell'Annam nel 1885. Divenne membro dell'Académie des sciences nel 1882. Si veda *Enciclopedia Italiana* (M. Camis).
*Lettera n° 49.*

**Bertini** Giuseppe (Milano 1825-1898). Fu un pittore romantico e verista. Figlio di Giovanni, pittore ticinese di vetrate che realizzò quelle dei finestroni dell'abside del Duomo di Milano, entrò a 13 anni all'Accademia di Belle Arti di Brera, mentre seguiva il padre in diversi cantieri. Fu docente all'Accademia, ne divenne direttore dal 1882 e, tra i suoi allievi, ebbe Tranquillo Cremona; fu inoltre il primo direttore e amministratore del Museo Poldi Pezzoli di Milano dal 1881. Si veda DBI (A. Ottino Della Chiesa).
*Lettera n° 38.*

**Betocchi** Alessandro (1823-1909). Ingegnere, Ispettore del Genio civile, fece parte della Giuria dell'Esposizione del 1881, nella IV sezione (Ingegneria e Lavori Pubblici), insieme a Clericetti. Dal 1867 fu socio nazionale per la Meccanica della Reale Accademia dei Lincei.
*Lettera n° 44.*

**Betti** Enrico (Pistoia 1823 - Soiana, Pisa 1892). Dopo aver compiuto studi classici, si laureò in Matematica a Pisa nel 1846. Nel 1848 prese parte alla battaglia di Curtatone con il battaglione universitario. Dopo aver insegnato in un liceo, nel 1857 ebbe la cattedra presso l'Università di Pisa, dove dal 1863 diresse anche la Scuola normale superiore. Divenne socio nazionale dei Lincei dal 1875 , e membro di numerose altre accademie e società scientifiche. Fu deputato dal 1861 al 1867 e dal 1874 al 1876 e fu quasi ininterrottamente membro del Consiglio della Pubblica Istruzione del Regno d'Italia. Nel 1884 venne nominato Senatore. Si vedano [Cerroni-Martini, 2009], [Giacardi-Tazzioli, 2012], DSB (E. Carruccio), DBI (N. Virgopia) e [Tricomi 1962].

Citato in [Cerroni-Fenaroli, 2007], [D'Agostino, 2007], [Canepa-Fenaroli, 2009], [Enea, 2009], [Enea-Gatto, 2009], [Giacardi-Tazzioli, 2013].
*Lettera n° 42.*

**Biffi** Serafino (Milano 1822-1899). Medico, si interessò dapprima dell'anatomia e della fisiologia del sistema nervoso centrale, poi si dedicò allo studio delle malattie psichiche e dei problemi connessi all'assistenza dei malati psichiatrici lavorando in stretta collaborazione con A. Verga. Dal 1864 fu membro effettivo dell'Istituto Lombardo di scienze e di lettere. Si veda DBI (G. Coari).
*Lettere n° 30,33.*

**Biondelli** Bernardino (Zevio, Verona 1804 - Milano 1886). Linguista, professore d'archeologia e numismatica, fu Direttore del Gabinetto Numismatico di Milano. Dal 1854 fu membro effettivo dell'Istituto Lombardo di scienze e di lettere. Si veda DBI (T. De Mauro).
*Lettere n° 30, 55.*

**Bismarck** Otto von (Schönhausen 1815 - Friedrichsruh 1898). Politico tedesco. Si veda *Enciclopedia Italiana.*
*Lettera n° 47.*

**Boito** Camillo (Roma 1836 - Milano 1914). Architetto e scrittore, dal 1860 tenne per quarantotto anni la cattedra di Architettura all'Accademia di Brera. Fu un'autorità indiscussa nel campo del restauro dei monumenti e fece parte della Giunta superiore delle belle arti per un lungo periodo. Si veda DBI (G. Miano).
*Lettera n° 48.*

**Bonghi** Ruggiero (Napoli 1826 - Torre del Greco 1895). Dopo aver partecipato nel 1847-'48 ai tentativi tesi a ottenere la Costituzione dal re Ferdinando II e aver partecipato alla guerra del 1848-'49, si trasferì prima a Firenze e poi a Torino. Nel 1860 rientrò a Napoli dove fondò e diresse il *Nazionale*. Nel 1862, tornato a Torino per insegnare Letteratura greca, fondò *La Stampa* che uscì fino al 1865; in quell'anno si trasferì a Firenze dove insegnò Letteratura latina. Nel 1867 venne nominato a Milano professore di Storia antica e dal 1866 al 1874 diresse, a Milano, *La Perseveranza*. Nel 1871 divenne professore di Storia Antica a Roma. Fu Ministro della Pubblica Istruzione dal 27 settembre 1874 al 24 marzo 1876. Fu membro dell'Accademia dei Lincei e di numerose altre accademie e istituti culturali, presidente dell'Associazione della stampa e della società Dante Alighieri a partire dal 1889. Si veda DBI (P. Scoppola).
Citato in [Cerroni-Martini, 2009], [Giacardi-Tazzioli, 2013].
*Lettere n° 15, 16.*

**Bonolis** Alfonso. Tra il 1873 e il 1912 scrisse alcuni articoli di Algebra, di Teoria dei Numeri e di Meccanica razionale. Fu autore anche di un trattato di Topografia. Prima di venire nominato Professore per la Cattedra di Scienze delle Costruzioni a Napoli aveva insegnato Meccanica Agraria e Costruzione Rurale nell'Istituto Agrario di Caserta.
Citato in [Cerroni-Martini, 2009], [Giacardi-Tazzioli, 2013], [Palladino-Mercurio, 2011].
*Lettera n° 41.*

**Bossi** Benigno (Como 1788 - Ginevra 1870). Patriota italiano, legato al gruppo milanese del *Conciliatore*; nel 1821 rappresentò i liberali lombardi al moto insurrezionale di Torino e nel 1824 fu condannato a morte in contumacia. Aderì alla società segreta dei "Sublimi maestri perfetti" a Ginevra e nel 1848 fu nominato a Londra rappresentante del governo provvisorio di Milano. Si vedano DBI (C. Francovich) e www.duevoltiditeresa.altervista.org
*Lettera n° 50.*

**Brioschi** Francesco (Milano 1824-1897). Laureatosi in Ingegneria a Pavia nel 1845, fu professore di Matematica applicata presso l'Università pavese dal 1852 al 1861. Fu segretario del Ministero della Pubblica Istruzione dal 1861 al 1863 e dal 1870 al 1882 fu nel Consiglio Esecutivo del Ministero della Pubblica Istruzione. Nel 1863 fondò l'Istituto Tecnico Superiore di Milano (l'attuale Politecnico), di cui divenne Direttore e dove insegnò Matematica e Idraulica fino alla morte. Per un breve periodo fu anche Deputato e, dal 1865, Senatore. Dal 1866, insieme a Luigi Cremona, prese la direzione degli *Annali di Matematica pura ed applicata*. I suoi contributi più significativi alla ricerca matematica sono legati alla teoria delle equazioni algebriche e alla teoria dei determinanti. Si vedano DSB (J. Pogrebyssky), DBI (N. Raponi, E. Ferri) e [Tricomi 1962].
Citato in [D'Agostino, 2007], [Cerroni-Martini, 2007], [Cerroni-Martini, 2009], [Canepa-Fenaroli, 2009], [Enea-Gatto, 2009], [Enea, 2009], [Giacardi-Tazzioli, 2012].
*Lettere n° 1, 3, 8, 9, 10, 11, 12, 14, 15, 17, 18, 33, 34, 38, 41, 42, 44, 47, 48, 50, 53, 57.*

**Broglio** Emilio (Milano 1814 - Roma 1892). Segretario del Governo Provvisorio di Milano nel 1848, fu Deputato al Parlamento subalpino (1848-1849). Riparò a Torino quando gli Austriaci tornarono a Milano. Fu Ministro della Pubblica Istruzione nel primo governo Menabrea dal 18 novembre 1867 al 5 gennaio 1868 e nel secondo governo Menabrea dal 5 gennaio 1868 al 13 maggio 1869. Fu vicepresidente della Camera del Regno d'Italia (1869-1870). Si veda DBI (N. Raponi).
Citato in [Cerroni-Martini, 2009].
*Lettera n° 7.*

**Bussone** Francesco, Conte di Carmagnola (Carmagnola 1385 ca. - Venezia 1432). Condottiero; alla sua figura si ispira la prima tragedia di Alessandro Manzoni. Si veda DBI (D.M. Bueno de Mesquita).
*Lettera n° 29.*

**Byron** George Gordon, Lord Byron (Londra 1788 - Missolungi 1824). Poeta e politico inglese. Si veda *Enciclopedia Italiana* (M. Praz).
*Lettere n° 51, 54.*

**Caccia**. Non identificato.
*Lettera n° 50.*

**Camperio** Manfredo (Milano 1826 - Napoli 1899). Viaggiatore, scrittore, uomo politico. Fu educato in casa con il cosiddetto "metodo inglese". Diciottenne, si dedicò a una vivace attività cospirativa antiaustriaca, partecipò alle Cinque giornate di Milano e alla guerra di indipendenza. Dopo aver viaggiato in Europa, nel 1850 si recò in Australia tentando la fortuna come cercatore d'oro. Rientrato in patria combatté ancora contro l'Austria. Nel 1869, dopo essere stato in Norvegia, in Egitto, a Ceylon e in India, divenne consigliere comunale a Milano fino al 1875 e assessore all'edilizia. Nel 1874 fu eletto deputato e, a Roma, entrò in contatto con la Società geografica italiana. Nel 1877 fondò il periodico *L'Esploratore* e negli anni successivi si dedicò alle esplorazioni commerciali in Africa e in Estremo Oriente. Si veda DBI (M. Carazzi).
*Lettera n° 6.*

**Cannizzaro** Stanislao (Palermo 1826 - Roma 1910). Nel 1845 e nel 1846, prima a Pisa e poi a Torino, fu assistente di Raffaele Piria (1815-1865), il chimico che per primo preparò l'acido salicilico. Partecipò ai moti siciliani del 1848 e fu condannato a morte: nel maggio 1849 dovette fuggire a Marsiglia. Nell'ottobre successivo raggiunse Parigi e qui operò presso il laboratorio di Michel Eugène Chevreul. Nel 1851 tornò in Italia e ottenne la cattedra di Chimica e Fisica al Collegio Nazionale di Alessandria. Nel 1855 diventò professore di Chimica all'Università di Genova. Dopo aver insegnato a Pisa e Napoli, occupò la cattedra di Chimica organica e inorganica a Palermo fino al 1871, quando ottenne una cattedra di Chimica all'Università di Roma e diventò Senatore per i suoi meriti scientifici. Nel 1858 pubblicò l'opera Sunto di un corso di filosofia chimica in cui pose le basi del moderno sistema atomico. Fu Senatore e membro del Consiglio superiore della Pubblica Istruzione e svolse un importante ruolo nell'educazione scientifica in Italia. Si vedano DBI (A. Gaudiano e D. Marotta) e DSB (H.M. Leicester).
Citato in [Cerroni-Fenaroli, 2007], [Canepa-Fenaroli, 2009], [Giacardi-Tazzioli, 2012].
*Lettera n° 42.*

**Cantoni Jung** Bice. Moglie di Giuseppe Jung.
*Lettera n° 23.*

**Cantoni** Gaetano (Milano 1815-1887). Si laureò in Medicina e chirurgia e, dopo aver esercitato la professione medica per qualche anno, si dedicò all'agronomia. Compì osservazioni sperimentali e pubblicò molti lavori, tra i quali un corso di agricoltura pratica ad uso dei contadini. Patriota attivista e liberale, partecipò alle Cinque giornate di Milano e, esule in Svizzera insieme al fratello Giovanni, tenne un corso di lezioni di agricoltura presso un liceo di Lugano. Tornato a Milano, nel 1870 fondò e diresse fino alla morte la Scuola d'agricoltura, continuando la sua azione per il miglioramento dell'agricoltura sia nell'aspetto scientifico (trovò ad esempio il metodo per fermare un'epidemia che stava decimando la popolazione dei bachi da seta) che in quello divulgativo. Fu consigliere comunale a Milano dal 1874 al 1878, cooperò con il Consiglio superiore della Pubblica Istruzione partecipando a varie missioni in Italia e all'estero. Fu membro effettivo del R. Istituto Lombardo dal 1874 e socio di altre accademie e società scientifiche. Si veda DBI (R. Giusti).
*Lettere n° 33, 34, 41.*

**Cantoni** Giovanni (Milano 1818-1897). Ingegnere e Senatore, nel 1859 fu incaricato di insegnare Fisica presso la Scuola Reale Superiore di Milano e l'anno seguente fu chiamato a Pavia a occupare l'analoga cattedra. Dal 1874 tenne la direzione del Servizio centrale di Meteorologia presso il Ministero dell'Agricoltura. Ricoprì varie cariche istituzionali: fu membro del Consiglio superiore della Pubblica Istruzione, Preside della Facoltà di Scienze di Pavia e Rettore dell'università. Fu membro effettivo del R. Istituto Lombardo dal 1863. Suo principale interesse fu la fisica sperimentale, ma si occupò anche di Storia della fisica e di divulgazione scientifica. Si veda DBI (G. Gliozzi).
Citato in [Cerroni-Fenaroli, 2007], [Canepa-Fenaroli, 2009], [Giacardi-Tazzioli, 2012].
*Lettere n° 33, 34, 41, 42.*

**Casorati** Felice (Pavia 1835-1890). Laureatosi in Ingegneria a Pavia nel 1856, restò in questa università come assistente. Nel 1858 partecipò, insieme a Enrico Betti e a Francesco Brioschi, a un viaggio scientifico in Francia e Germania che molto contribuì a rompere l'isolamento in cui era vissuta fino ad allora la matematica italiana. Nel 1859, a soli 24 anni, fu nominato professore di Algebra e Geometria analitica nell'Università pavese, ove più tardi (1863) insegnò Analisi infinitesimale e poi anche Analisi superiore. Si vedano DBI (E. Togliatti) e [Tricomi 1962].
Citato in [D'Agostino, 2007], [Enea-Gatto, 2009], [Cerroni-Martini, 2009], [Enea, 2009].
*Lettere n° 9, 30, 33, 37, 41.*

**Cavallini** Achille (Milano 1812-1881). Si laureò in Matematica nel 1833 e si dedicò alla professione di ingegnere-architetto specializzandosi nelle discipline tecnico-giuridiche soprattutto nel campo dell'idraulica. Nel 1864 venne chiamato da Brioschi alla nuova cattedra di Giurisprudenza agricola ed Elementi di diritto amministrativo che tenne fino al 1876. Fu presidente del Collegio degli Ingegneri ed Architetti di Milano.
*Lettera n° 16.*

**Cecchi** Filippo (Borgo Buggiano, Pistoia 1822 - Firenze 1887). Padre Filippo Cecchi, al secolo Giulio Isdegerde, religioso dell'ordine degli Scolopi, fu un fisico, inventore e costruttore di strumenti scientifici. Si occupò di elettromagnetismo, di meteorologia e sismologia. Diresse l'Osservatorio Ximeniano di Firenze dal 1872 al 1887 dotandolo di un importante centro sismologico. Fu membro dell'Accademia Nazionale dei Lincei e fondò la Società Metereologica Italiana. Si veda DBI (N. Janiro).
*Lettera n° 23.*

**Celoria** Giovanni (Casale Monferrato, Alessandria 1842 - Milano 1920). Studiò all'Università di Torino dove si laureò in Ingegneria nel 1863 e poi per tutta la vita fu astronomo all'Osservatorio di Brera a Milano, dove entrò sotto la guida di G. Schiaparelli. Insegnò Geodesia al Politecnico di Milano e nel 1902 fu nominato

Presidente della Commissione Geodetica Italiana. Venne nominato Senatore nel 1909. Si veda DBI (N. Janiro).
Citato in [Giacardi-Tazzioli, 2012].
*Lettere n° 33, 37.*

**Ceradini** Cesare (Milano 1844 - Roma 1935). Laureatosi in Ingegneria civile presso il Politecnico di Milano nel 1887, fu nominato professore straordinario di Costruzioni stradali e idrauliche a Palermo. Nel 1872 fu chiamato a Roma alla cattedra di Statica delle costruzioni e nel 1875, sempre a Roma, fu nominato ordinario per la Meccanica applicata alle costruzioni. Nel 1909 successe a Valentino Cerruti nella direzione della Scuola per gli ingegneri di Roma e in tale carica rimase fino al 1922 quando già da tre anni aveva lasciato l'insegnamento per limiti di età. Ebbe dallo stato numerosi incarichi di fiducia. Nel campo della tecnica sono notevoli i suoi studi sulla teoria delle capriate e delle vòlte, sulla spinta delle terre, e sulla trave continua. Suo è un trattato di meccanica applicata alle costruzioni. Si veda DBI (E. Ferri).
*Lettere n° 8, 10, 18, 46.*

**Ceriani** Antonio Maria (Uboldo, Varese 1828 - Milano 1907). Fu sacerdote della diocesi di Milano, biblista e studioso di lingue orientali. Dal 1857 entrò a far parte del Collegio dei dottori della Biblioteca Ambrosiana. Fu membro effettivo dell'Istituto Lombardo di Scienze e di Lettere dal 1862. Si veda DBI (F. Parente).
*Lettera n° 30.*

**Cerruti** Valentino (Crocemosso 1850-1909). Laureatosi in Ingegneria civile al Politecnico di Torino nel 1873 con una tesi sulla statica dei sistemi articolati, il 15 dicembre del medesimo anno fu nominato assistente di Idraulica presso la Scuola d'Applicazione per gli Ingegneri dell'Università di Roma, diretta da Luigi Cremona. Successivamente divenne insegnante privato dei figli di Quintino Sella. Quattro anni dopo ottenne la cattedra di Meccanica Razionale presso la Scuola di Ingegneria di Roma dove nell'ottobre del 1877 vinse il concorso a cattedra. La sua attività non si esaurì tuttavia nell'insegnamento: nel 1880 collaborò con Luigi Cremona alla riorganizzazione della Biblioteca Alessandrina dal Ministro Coppino che poi lo volle - nel 1886 - Segretario Generale del Ministero della Pubblica Istruzione. Tra il 1888 e il 1892 fu Rettore dell'Università di Roma, nel 1901 fu nominato Senatore e nel 1903 fu chiamato a dirigere - come successore di Luigi Cremona - la Scuola per gli Ingegneri. Fu socio di diverse Accademie straniere e italiane e la sua produzione scientifica comprende soprattutto lavori di Teoria dell'elasticità. Si vedano DBI (E. Pozzato) e [Tricomi 1962].
Citato in [Cerroni-Fenaroli, 2007], [Cerroni-Martini, 2009], [Enea-Gatto, 2009].
*Lettere n° 17, 44.*

**Clericetti** Cecilia, moglie di Celeste.
*Lettere n° 1, 3, 4, 6, 7, 10, 11, 12, 13, 14, 15, 16, 17, 18, 19, 20, 21, 22, 24, 25, 26, 28, 30, 31, 35, 37, 39, 45, 46, 47, 48, 50, 51, 52, 54, 55, 56, 57, 58.*

**Clericetti** Emilio, figlio di Celeste.
*Lettere n° 1, 3, 6, 14, 19, 24, 28, 31, 35, 44, 48, 50, 51, 52, 54, 55, 56.*

**Clericetti** Emilio, fratello di Celeste.
*Lettere n° 27, 28, 50, 55.*

**Clericetti** Guido, figlio di Celeste.
*Lettere n° 1, 3, 6, 10, 11, 13, 14, 19, 22, 24, 28, 31, 35, 44, 46, 48, 50, 52, 55, 56, 58.*

**Clericetti** Pietro, fratello di Celeste.
*Lettere n° 7, 18, 52, 55.*

**Codazza** Giovanni (Milano 1816 - Como 1877). Laureatosi in Ingegneria e Architettura a Pavia nel 1837, insegnò nelle scuole medie. Nel 1842 fu nominato professore di Geometria descrittiva nell'Università di

Pavia, di cui fu anche Rettore nel 1857-1858. Nel 1862-1863 fu Sindaco della città. Nel 1863 si trasferì al R. Istituto Tecnico Superiore di Milano per insegnarvi Fisica Tecnologica e nel 1868 accettò la nomina a docente di Fisica industriale al Regio Museo Industriale di Torino che diresse dal 1870 al 1877. I suoi studi furono essenzialmente dedicati alla Fisica matematica, alla Geometria descrittiva, alla Fisica tecnologica. Fu socio dell'Accademia Nazionale dei Lincei e dell'Istituto Lombardo dal 1854. Si vedano DBI (R. Ferola) e [Tricomi 1862].
Citato in [Cerroni-Martini, 2009], [Giacardi-Tazzioli, 2012].
*Lettera n° 12.*

**Colombo** Giuseppe (Milano 1836-1921). Laureatosi a Pavia in Matematica, insegnò Meccanica dal 1857 al 1883 alla Società d'incoraggiamento di arti e mestieri di Milano e, a partire dal 1865, anche al Politecnico di Milano. Si occupò di questioni legate allo sviluppo industriale relative soprattutto alla produzione e alla distribuzione dell'energia elettrica italiana creando la Società Edison. Nel 1897, alla morte di Brioschi, divenne direttore del Politecnico, carica che conservò fino al 1911. Tra le numerose pubblicazioni la più nota è il *Manuale dell'Ingegnere* (1877), a lungo utilizzato come testo di base dagli ingegneri italiani. Si dedicò anche attivamente alla vita politica: nel 1899-1900 fu presidente della Camera dei deputati e nel 1900 venne nominato Senatore. Ebbe anche due brevi incarichi di governo, nel 1891 come ministro delle Finanze e nel 1896 come ministro del Tesoro. Si veda DBI (R. Cambria).
*Lettere n° 9, 13, 19, 30, 33, 37, 38, 44, 50, 53, 55, 56.*

**Conti** Pietro, tenente colonnello. Non meglio identificato.
*Lettera n° 18.*

**Comotto** Paolo (Brianzè, Vercelli 1824 - Roma 1897). Si laureò a Torino in ingegneria e architettura civile. Nel 1861 lavorò alla prima aula del Parlamento italiano dove Vittorio Emanuele II pronunciò lo storico discorso dopo l'incontro di Teano. Nel 1870 venne incaricato di progettare un'aula provvisoria all'interno di Montecitorio per ospitare la Camera dei Deputati che venne terminata nel novembre dell'anno successivo. Nel 1833, in qualità di ispettore del Genio civile fece parte di una Commissione per le prescrizioni edilizie ad Ischia istituita dopo il terremoto. Si veda DBI (F. Quinterio).
*Lettera n° 1.*

**Cornalia** Emilio (Milano 1824-1882). Prima Curatore, divenne Direttore del Museo Civico di Storia Naturale di Milano nel 1869. Dal 1851 al 1866 partecipò alla fondazione della Società Entomologica Italiana e fu autore di importanti opere di Entomologia applicata. Nel 1875 divenne socio nazionale dell'Accademia dei Lincei. Si veda DBI (F. Di Trocchio).
Citato in [Giacardi-Tazzioli, 2012].
*Lettere n° 33, 34, 37, 41, 42, 50, 55.*

**Corradi** Alfonso (Bologna 1833 - Pavia 1892). Si laureò a Bologna in Medicina nel 1856 e nel 1859 ebbe la cattedra di Patologia generale all'università di Modena. Nel 1863 si trasferì a Palermo sulla cattedra di Patologia generale e per quattro anni si dedicò anche agli studi chimici e farmacologici sotto la guida di Cannizzaro. Ottenne quindi il trasferimento all'Università di Pavia. Nel 1875 fu preside della facoltà di Medicina e nel 1876 divenne rettore dell'Università. Fu membro effettivo dell'Istituto Lombardo di Scienze e di Lettere dal 1874. Si veda DBI (B. Zanobio, G. Armocida).
*Lettera n° 33.*

**Cremona Perozzi** Elena (Pavia 1856-?), figlia maggiore di Luigi.
*Lettere n° 5, 7, 14, 17, 19, 20, 22, 42, 45.*

**Cremona Cozzolino** Itala (Bologna 1865 - Genova 1939), figlia di Luigi.
*Lettere n° 7, 11, 14, 17, 19, 22, 45, 51, 54, 55, 56, 57.*

**Cremona** Luigi (Pavia 1830 - Roma 1903). Si vedano le diverse commemorazioni su questo sito e DBI (U. Bottazzini, L. Rossi).
*Lettera n° 33.*

**Cremona** Tranquillo (Pavia 1837 - Milano 1878). Pittore, fratello di Luigi. Si veda DBI (A. Pino Adami).
*Lettera n° 21.*

**Cremona** Vittorio (Bologna 1861-?), figlio di Luigi.
*Lettere n° 6, 7, 11, 12, 13, 14, 19, 20, 22, 28, 45, 50, 51, 54, 55, 56, 57, 58.*

**Crispi** Francesco (Ribera, Agrigento 1818 - Napoli 1901). Politico, fu una figura di spicco del Risorgimento. Fu uno degli organizzatori della Rivoluzione siciliana del 1848 e l'ideatore della spedizione dei Mille alla quale partecipò. Fu quattro volte presidente del Consiglio: dal 1887 al 1891 e dal 1893 al 1896. Nel primo periodo fu anche ministro degli Esteri e ministro dell'Interno, nel secondo anche ministro dell'Interno. Si veda DBI (F. Fonzi).
*Lettera n° 58.*

**Culmann** Karl (Bergzabern, Germania 1821 - Zurigo 1881). Ingegnere, seguì i lavori ferroviari attraverso il massiccio del Fichtelgebirge fino al 1847 occupandosi soprattutto della costruzione dei ponti. Nel 1849 ottenne una cattedra al Politecnico di Zurigo, del quale fu anche direttore dal 1872 al 1875. Può essere considerato il fondatore della Statica grafica, da lui ridotta a un corpo omogeneo di dottrina. Si veda *Enciclopedia Italiana* (G. Albenga).
*Lettera n° 19.*

**Curioni** Giovanni (Invorio Inferiore, Novara 1831 - Torino 1887). Laureatosi a Torino in Ingegneria e architettura, iniziò la carriera di docente, dapprima quale assistente di Geometria pratica e in seguito alla cattedra di Costruzioni di architettura e geometria pratica nell'istituto di insegnamento tecnico (poi Scuola d'applicazione per ingegneri). Nel 1861 pubblicò il *Corso di topografia* con il quale poté accedere alla cattedra di Geometria pratica, costruzioni ed estimo nell'istituto industriale e professionale di Torino. Nel 1862 superò l'esame di dottore aggregato in scienze fisiche e matematiche trattando l'argomento "Sulla misura di una base geometrica". Fu tra i promotori della Società degli ingegneri ed industriali di Torino di cui anni dopo divenne presidente. Nel 1865 fu nominato direttore della Scuola di applicazione, nel 1866 divenne professore straordinario di Costruzioni e nel 1868 divenne ordinario. Si veda DBI (B. Signorelli).
*Lettera n° 42.*

**Cusani** Luigi. Nel 1881 fu assessore nel comune di Milano.
*Lettera n° 48.*

**Da Vinci** Leonardo. Si veda DBI (P. Marani).
*Lettera n° 29.*

**De Cristoforis** Malachia (Milano 1832-1915). Si laureò in Medicina e chirurgia nel 1856 a Pavia. Durante gli anni universitari si dedicò attivamente alla causa risorgimentale; poi si arruolò nei Cacciatori delle Alpi e svolse attività come ufficiale medico nel corso di vari combattimenti raggiungendo anche i Mille di Garibaldi in Sicilia nel 1860. Nel 1866 si arruolò nuovamente e venne insignito della croce militare di Savoia. Il suo interesse scientifico si era sempre orientato verso l'ostetricia e la ginecologia e di ciò si occupò presso l'ospedale Maggiore di Milano. Scrisse un corposo trattato sulle malattie delle donne e nel 1883 divenne libero docente presso l'Università di Napoli. Prese parte attiva alla vita politica del comune di Milano, in particolare fu presidente della Società di cremazione, fu massone e nel 1905 fu nominato Senatore. Si veda DBI (G. Armocida, G. Bock Berti).
*Lettera n° 38.*

**De Marchi** Odoardo, fu assistente di Costruzione di ponti, opere marittime e lavori in terra all'Istituto Tecnico di Milano dal 1881 al 1888 sostituendo Clericetti durante la malattia. Fece parte del Collegio degli Ingegneri ed Architetti di Milano.
*Lettera n° 57.*

**Errera** Alberto (Venezia 1841 - Napoli 1894). Patriota ed economista, si laureò a Padova nel 1867 e divenne insegnante di Economia e diritto commerciale presso l'Istituto tecnico di Venezia. Si trasferì a Milano nel 1874 dove insegnò Statistica ed economia e Istituzioni civili e morali presso il R. Istituto tecnico di S. Marta. Ebbe una feconda produzione di articoli nei quali esponeva gran quantità di dati e notizie su diverse attività economiche del Veneto che si estese, a partire dal 1874, all'ambito nazionale. Dal 1877 fu a Napoli nel locale Istituto tecnico. Si veda DBI (A. Polsi).
*Lettera n° 17.*

**Fairbairn** William (Kelso, Scozia 1789 - Moor Park, Inghilterra 1874). Ingegnere civile e strutturale, costruì diverse navi per la Peninsular and Oriental Steam Navigation Company. Socio della Royal Society e nel 1861 presidente della British Science Association. Si veda ODNB.
*Lettera n° 49.*

**Fano** Enrico (Milano 1834-1899). Si laureò in Giurisprudenza a Pavia e intraprese la professione forense specializzandosi in "questioni bancarie". Fu molto attivo politicamente e lavorò per promuovere lo sviluppo delle società operaie e del mutualismo. Venne eletto deputato nel 1867 e nel 1870 si recò a Londra come membro della commissione italiana per la Mostra internazionale operaia. Si veda DBI (N. Dell'Erba).
*Lettera n° 6.*

**Farini** Domenico (Montescudo, Rimini 1834 - Roma 1900). Politico, ufficiale del Genio, venne eletto deputato di Ravenna dal 1864 al 1866. Fu segretario alla Camera schierandosi con il gruppo di centrosinistra, ma rifiutò funzioni politiche, incarichi governativi e diplomatici. Fu nominato Senatore nel 1886 e presidente del Senato l'anno successivo. Si veda DBI (F. Bartoccini).
*Lettere n° 43, 58.*

**Favero** Giovanni Battista (Crespano Veneto, Treviso 1832 - Roma 1906). Laureatosi a Padova in Ingegneria civile edile fu addetto al servizio trazione delle Ferrovie Nord austriache con sede a Vienna fino al 1865. Tornato in Italia, lavorò alle Strade ferrate meridionali di cui fu anche direttore. Nel 1875 fu nominato straordinario, e poi ordinario nel 1878, della cattedra di Ponti e strade presso l'Università di Roma. Dal 1878 al 1902 tenne la cattedra di Ponti e strade nella Scuola d'applicazione per gli ingegneri in Roma. Divenne socio nazionale dei Lincei nel 1898. Si veda DBI (E. Pozzato).
*Lettere n° 24, 26.*

**Favre** Louis (Chêne-Bourg, Svizzera 1826 - Göschenen, Svizzera 1879). Architetto e ingegnere autodidatta svizzero, fu titolare dell'impresa alla quale venne affidata la costruzione del traforo ferroviario del San Gottardo. Favre morì di infarto durante un sopralluogo in cantiere. Si veda DSS (C. Courtiau).
*Lettera n° 29.*

**Ferrari Cremona** Elisa (Genova 1826 - Roma 1882). È la prima moglie di Luigi.
*Lettere n° 1, 3, 4, 6, 7, 11, 12, 16, 17, 22, 23, 27, 28, 32, 34, 41, 42, 43, 45, 47, 50.*

**Ferrini** Rinaldo (Milano 1831-1908). Conseguì la laurea in Ingegneria civile e architettura presso l'Università di Pavia nel 1853. A Milano insegnò dapprima Fisica e Matematica presso il Liceo S. Marta e, dal 1860 al 1868, divenne titolare della cattedra di Fisica e Meccanica presso il Politecnico dove l'anno successivo fu docente di Fisica tecnologica. Nel 1873 divenne membro effettivo del R. Istituto Lombardo. Viene considerato uno dei più noti fisici del tempo grazie alle numerose pubblicazioni riguardanti il calore, l'elettricità, il magnetismo, la ionizzazione dei gas e la materia radiante. Si veda DBI (E. Pozzato).
*Lettere n° 33, 37, 53.*

**Finali** Gaspare (Cesena 1829 - Marradi, Firenze 1914). Politico, partecipò al movimento repubblicano. Ricoprì numerose cariche: fu deputato nel 1865-70, venne nominato Senatore dal 1872 e fu Ministro dell'Agricoltura e del Commercio (1873-76), dei Lavori Pubblici (1887-91) e del Tesoro nel 1901. Fu docente di Contabilità dello Stato presso l'Università di Roma dal 1880 e divenne socio dei Lincei dal 1901. Si veda DBI (E. Orsolini).
*Lettera n° 7.*

**Frisiani** Paolo (Milano 1797-1880). Entrò nel 1820 all'Osservatorio Astronomico di Brera in cui, nel 1834, divenne secondo astronomo. Non fece osservazioni, ma si dedicò soprattutto all'insegnamento della Matematica. Fu autore di alcuni lavori di Astronomia teorica e di Analisi matematica. Si vedano DBI (F. Manzotti) e [Tricomi 1962].
Citato in [Canepa-Fenaroli, 2009].
*Lettere n° 30, 33, 38.*

**Frugoni** Carlo Innocenzo (Genova 1692 - Parma 1768). Librettista e poeta. Si veda DBI (G. Fagioli Vercellone).
*Lettera n° 7.*

**Gabet** Luigi (Roma 1823-1879). Studiò presso l'Accademia romana di belle arti dove apprese l'arte architettonica. Diresse importanti restauri di palazzi romani e nel 1870 partecipò a diverse commissioni di lavoro sull'assetto urbanistico di Roma. L'anno successivo ricevette l'incarico di realizzare Palazzo Madama. Nel 1873 fu eletto consigliere comunale di Roma e fece parte della Commissione edilizia. Nel 1874 divenne architetto capo dell'Ufficio tecnico della Provincia di Roma. Si veda DBI (F. Di Marco).
*Lettera n° 1.*

**Galanti** Antonio, docente di Agraria presso il Politecnico di Milano nel 1880.
*Lettera n° 38.*

**Garovaglio** Santo (Como 1805 - Pavia 1882). Si laureò in Chimica e in Medicina presso l'Università di Vienna. Divenne assistente alla cattedra di Botanica all'Università di Pavia e svolse un'intensa attività di studio e di ricerca. Dal 1839 al 1852 fu docente di Fisica, Chimica e Botanica a Pavia e nel 1852 venne nominato professore e direttore dell'Orto botanico della stessa città. Fu membro di molte Società, in particolare del R. Istituto Lombardo di Scienze e Lettere dal 1854. Si veda DBI (M. Rossi).
*Lettera n° 33.*

**Ghezzi**. Studente non meglio identificato.
*Lettere n° 45, 47, 48, 49.*

**Giambastiani** Angelo, ingegnere.
*Lettera n° 51.*

**Giani** (probabilmente si tratta di Lapo Gianni), poeta.
*Lettera n° 7.*

**Gioda** Carlo, nel 1873 fu Regio Provveditore agli Studi di Milano.
*Lettera n° 12.*

**Giovannini** Giuseppe, ingegnere. Segretario del R. Istituto Tecnico di Milano negli anni relativi a questo carteggio.
*Lettere n° 8, 9, 10, 47.*

**Golgi** Camillo (Corteno, Brescia 1843 - Pavia 1926). Laureato in Medicina a Pavia, fu un eminente istologo e patologo. Esercitò come primario ad Abbiategrasso e, in seguito a un'importante scoperta sullo studio delle cellule nervose, nel 1875 divenne docente di Istologia all'Università di Pavia, nel 1979 si spostò sulla cattedra di Anatomia a Siena, tornò a Pavia come docente di Istologia e di Patologia generale e nel 1893 divenne Rettore di tale Ateneo. Ricevette il Premio Nobel per la Medicina e la Fisiologia nel 1906. Fu socio dei Lincei e dal 1882 del R. Istituto Lombardo. Venne eletto Senatore dal 1900. Si veda DBI (G. Cimino).
*Lettera n° 41.*

**Gorini** Paolo (Pavia 1813 - Lodi 1881). Si laureò in Matematica all'Università di Pavia nel 1833; l'anno successivo ebbe a Lodi una cattedra per l'insegnamento delle Scienze naturali, da cui si dimise per ragioni politiche nel 1857. È noto soprattutto come chimico per le sue ricerche volte alla conservazione dei corpi evitandone la putrefazione. Si veda DBI (F. Conti).
Citato in [Giacardi-Tazzioli, 2012].
*Lettera n° 38.*

**Guidi** Camillo (Roma 1853-1951). Si laureò nel 1877 presso la R. Scuola di applicazione degli ingegneri di Roma, dove rimase come assistente presso la cattedra di Statica grafica fino al 1881, quando venne incaricato dello stesso corso presso la R. Scuola d'ingegneria di Torino. Nel 1887 vinse il concorso per la cattedra di Scienza delle costruzioni, sempre a Torino, e la tenne sino al 1928. Si veda DBI (T. Iori).
*Lettera n° 46.*

**Guj** Enrico (Roma 1841-1905). Si laureò in Ingegneria presso l'Università di Roma nel 1864. Affiancò l'attività professionale con quella di docenza: all'Archiginnasio insegnò diverse materie - Geometria grafica, Idrometria, Geometria descrittiva, Architettura statica e idraulica - negli anni dal 1869 al 1871, poi alla nuova Scuola di applicazione per gli ingegneri ricoprì la cattedra di architettura tecnica dal 1873 al 1905. Si veda DBI (F. Di Marco).
*Lettera n° 44.*

**Hajech** Camillo (?-1883). Professore ordinario di Fisica nel R. Liceo Beccaria di Milano. Fu socio del R. Istituto Lombardo dal 1860.
*Lettere n° 33, 34, 37, 55.*

**Jung** Giuseppe (Milano 1845-1926). Si laureò in Matematica nel 1867 a Napoli con Giuseppe Battaglini e subito dopo tornò a Milano dove divenne assistente di Luigi Cremona all'Istituto Tecnico Superiore di Milano. Nel 1876 fu nominato professore straordinario di Geometria proiettiva e Statica grafica e professore ordinario nel 1890. Lasciò alcuni lavori di Geometria e di Statica. Si vedano *Enciclopedia Italiana* (G. Albenga) e [Tricomi 1962].
Citato in [Cerroni-Fenaroli, 2007], [Cerroni-Martini, 2009], [Giacardi-Tazzioli, 2012].
*Lettere n° 11, 23, 24, 44.*

**Körner** Guglielmo (Kassel, Germania 1839 - Milano 1925). Si laureò nel 1860 in Chimica presso l'Università di Giessen dove rimase per quattro anni in qualità di assistente per la Chimica sperimentale. Nel 1867 divenne assistente di Stanislao Cannizzaro a Palermo; nel 1870 si spostò a Milano come professore di Chimica organica sperimentale nella Scuola superiore di agraria, che diresse dal 1899 al 1914 e dove rimase come docente fino al 1922. Dal 1874 insegnò anche al Politecnico. Fu membro effettivo del R. Istituto Lombardo dal 1880. Si veda DBI (G.P. Marchese).
*Lettera n° 41.*

**Lamartine** (de) Alphonse (Maçon, Francia 1790 - Parigi 1869). Poeta, scrittore, storico e politico francese. Si veda *Enciclopedia Italiana* (P.P. Trompeo).
*Lettera n° 3.*

**Lamé** Gabriel (Tours, Francia 1795 - Parigi 1870). Nel 1817 si laureò all'École Polytechnique a Parigi e nella stessa città nel 1820 si laureò in Ingegneria all'École des Mines. Dal 1820 al 1832 visse a San Pietroburgo, in Russia, come professore e ingegnere all'Institut et Corps du Genie des Voies de Communication. Nel 1832 tornò a Parigi e ottenne la cattedra di Fisica all'École Polytechnique. Nel 1836 divenne Ingegnere capo di miniere e progettò la ferrovia da Parigi a Versailles e quella da Parigi a S. Germain, che furono aperte nel 1837. Nel 1843 divenne membro dell'Académie des Sciences e nel 1851 ebbe la cattedra di Fisica matematica e Probabilità alla Sorbona. Si veda DSB (S.L. Greitzer).
Citato in [Cerroni-Martini, 2009], [Enea-Gatto, 2009], [Enea, 2009].
*Lettera n° 49.*

**Leone XIII** (Carpineto Romano, Roma 1810 - Roma 1903). Papa della Chiesa cattolica dal 1878. Si veda DBI (F. Malgeri).
*Lettera n° 47.*

**Lombardini** Elia (Alsazia 1794 - Milano 1878). Tra i più importanti ingegneri idraulici italiani, divenne socio dell'Accademia Nazionale delle Scienze nel 1865 e del R. Istituto Lombardo dal 1844. Si veda *Enciclopedia Italiana* (G. Albenga).
Citato in [Canepa-Fenaroli, 2009], [Enea-Gatto, 2009].
*Lettere n° 22, 55.*

**Loria** Leonardo. Docente di Ferrovie presso il Politecnico di Milano, nel 1899 fu il primo presidente del Collegio Ingegneri Ferroviari Italiani (CIFI).
*Lettere n° 8, 13, 57, 58.*

**Luzzatti** Luigi (Venezia 1841 - Roma 1927). Si laureò in Giurisprudenza presso l'Università di Padova. Fu giurista, politico ed economista; dal 1910 al 1911 fu presidente del Consiglio dei ministri e fondatore della Banca Popolare di Milano. Si veda DBI (P. Pecorari, P. Ballini).
*Lettera n° 58.*

**Maffei** Andrea (Molina di Val Ledro, Trento 1798 - Milano 1885). Si laureò in Giurisprudenza a Pavia nel 1820. Fu poeta e letterato, tradusse Milton, Shakespeare e Schiller. Nel 1879 fu nominato Senatore. Si veda DBI (M. Marri Tonelli).
*Lettera n° 55.*

**Maffei** Clara (Bergamo 1814 - Milano 1886). Contessa, moglie di Andrea Maffei, raccolse nel suo salotto milanese molti artisti legati al movimento patriottico. Per breve tempo fu esule a Locarno dove conobbe Giuseppe Mazzini. Si veda *Enciclopedia Italiana* (A.M. Ghisalberti).
*Lettera n° 55.*

**Maganzini** Italo. Ingegnere, allievo del R. Istituto Tecnico di Milano, divenne Ispettore Superiore del Genio Civile e Presidente di Sezione del Consiglio Superiore dei Lavori Pubblici.
*Lettere n° 17, 44.*

**Maggi** Leopoldo (Rancio Valcuvia, Varese 1840 - Pavia 1905). Nel 1863 si laureò in Scienze naturali e in Medicina e chirurgia presso l'Università di Pavia dove insegnò dapprima Mineralogia e Geologia e, dal 1875, Zoologia e Anatomia comparata. Fu membro effettivo del R. Istituto Lombardo dal 1879. Si veda DBI (F. Barbagli).
*Lettera n° 33.*

**Maglione** Giovanni, ragioniere. Fu il curatore provvisorio del fallimento della società anonima detta "Fabbrica lombarda di prodotti chimici", del cui Consiglio di Amministrazione Brioschi era presidente.
*Lettera n° 53.*

**Mantegazza** Paolo (Monza 1831 - S. Terenzo di Lerici, La Spezia 1910). Si laureò in Medicina a Pavia e iniziò la professione in America latina. Tornato in Italia nel 1860, ricoprì la cattedra di Patologia generale nell'Università di Pavia. Fu igienista, divulgatore e scrittore. Fu membro effettivo del R. Istituto Lombardo dal 1863 e Senatore dal 1876. Si veda DBI (G. Armocida, G.S. Rigo).
*Lettere n° 7, 28, 33, 34, 38.*

**Manzoni** Alessandro. Si veda DBI (P. Floriani).
*Lettere n° 7, 27.*

**Martelli** Giuseppe, docente di Lavori in terra e Disegno di opere stradali ed idrauliche al R. Istituto Tecnico Superiore.
*Lettere n° 6, 19, 44.*

**Massarani** Tullo (Mantova 1826 - Milano 1905). Fu letterato, buon intenditore d'arte e uomo politico. Per un ventennio fece parte del Consiglio provinciale di Milano. Nel 1876 venne nominato Senatore. Si veda DBI (R. Balzani).
*Lettera n° 48.*

**Melzi** Francesco, pittore. Fu allievo di Leonardo che lo nominò erede ed esecutore testamentario lasciandogli manoscritti e altro materiale grafico. Si veda DBI (F. Sorce).
*Lettera n° 29.*

**Mengoni** Francesco (Fontanelice, Bologna 1829 - Milano 1877). Architetto. Si veda DBI (E. Piersensini).
*Lettera n° 29.*

**Menabrea** Luigi Federico (Chambéry, Savoia, 1809 - Saint-Cassin, Chambéry, 1896). Si laureò dapprima in Ingegneria idraulica e poi in Architettura civile a Torino dove insegnò Scienza delle costruzioni e Geometria pratica. Ebbe tra i suoi allievi Quintino Sella. Ricoprì diverse cariche politiche assumendo più volte l'incarico di Ministro. Fu nominato Senatore del Regno nel 1860. Si veda DBI (P.A. Gentile).
*Lettera n° 29.*

**Metastasio** Pietro, poeta. Si veda *Enciclopedia Italiana* (A. Pompeati).
*Lettera n° 7.*

**Minghetti** Marco (Bologna 1818 - Roma 1886), statista. Si veda DBI (R. Gherardi).
*Lettera n° 58.*

**Nazzani** Ildebrando (Parma 1846-1931), ingegnere. Fu volontario nella guerra del 1866. Nel 1876 scrisse il primo trattato italiano di Idraulica, disciplina di cui tenne la cattedra alla Scuola d'Applicazione per gli ingegneri a Roma dal 1878 al 1915. Si veda *Enciclopedia Italiana.*
*Lettere n° 23, 24, 25, 26, 27, 28, 55.*

**Negri** Gaetano (Milano 1838 - Varazze, Savona 1902). Scrittore e uomo politico, fu sindaco di Milano dal 1884 al 1889, Senatore dal 1890 e socio corrispondente dei Lincei dal 1902. Si veda DBI (M. Soresina).
*Lettere n° 51, 57.*

**Olginti**, famiglia di Como, amici dei Clericetti. Non identificati.
*Lettera n° 45.*

**Padoa** Gustavo, studente di Clericetti non meglio identificato.
*Lettera n° 13.*

**Pagani** Luigi (Bergamo 1829 - Milano 1905), scultore. Si veda DBI (P. Bosio).
*Lettera n° 5.*

**Paladini** Ettore (Milano 1848-1930). Si laureò nel 1870 in Ingegneria civile presso il R. Istituto Tecnico di Milano dove insegnò Idraulica e costruzioni idrauliche. Fu presidente della Commissione per la fognatura di Milano nel 1901 e presidente del Collegio degli Ingegneri e Architetti di Milano.
*Lettere n° 8, 9.*

**Parini** Giuseppe, poeta. Si veda DBI (G. Nicoletti).
*Lettera n° 7.*

**Pavesi** Pietro (Pavia 1844 - Asso, Como 1907). Zoologo, insegnò nelle Università di Genova e di Pavia. A Pavia diresse il Museo di storia naturale e fu preside della Facoltà di Scienze. Fu membro effettivo del R. Istituto Lombardo dal 1883. Si veda DBI (A. Volpone).
*Lettera n° 33.*

**Perazzi** Costantino (Novara 1832 - Roma 1896). Si laureò in Ingegneria idraulica e architettura civile presso l'Università di Torino e continuò la sua formazione all'École des Mines dove conobbe Quintino Sella. Ricoprì diversi incarichi tecnici per il governo e intraprese la carriera politica. Venne nominato Senatore nel 1884. Si veda DBI (P. Gentile).
*Lettera n° 48.*

**Perelli** Guido, nel 1877 si laureò in Ingegneria e Architettura presso il R. Istituto Tecnico di Milano.
*Lettera n° 12.*

**Piuma** Carlo Maria (Genova 1837-1912). Laureatosi in Ingegneria all'Università di Genova nel 1860, fu Professore ordinario di Calcolo per la stessa Università dal 1882 al 1905. Si veda [Tricomi 1962].
Citato in [Cerroni-Fenaroli, 2007], [Cerroni-Martini, 2009], [Enea-Gatto, 2009].
*Lettera n° 53.*

**Poli** Baldassarre (Cremona 1795 - Milano 1883). Dopo la laurea in Giurisprudenza a Bologna si dedicò all'insegnamento e nel 1820 ottenne la cattedra di Filosofia al Liceo milanese di Porta Nuova (ora "Parini"'), in cui insegnò per 17 anni. Durante questa sua permanenza a Milano produsse vari saggi filosofici. Trasferitosi a Padova, insegnò Filosofia teorica e pratica, disciplina che comprendeva anche la Storia della filosofia, all'Università, dal 1837 al 1852. Fu Rettore a Padova nel 1849-50; e nel 1852, fu nominato, dall'I. R. Governo, Direttore generale dei ginnasi veneti. Nel 1857 divenne Consigliere scolastico e ispettore generale presso la Luogotenenza di Milano. Una volta lasciato il servizio attivo nel 1859, continuò a scrivere di Filosofia, propugnando la sintesi di razionalismo e di empirismo da lui indicata nei *Supplimenti.* Le sue posizioni furono praticamente ignorate, anche se, negli ultimi anni della vita, fu circondato da rispetto e stima. Fu assiduo alle sedute dell'Istituto Lombardo di Scienze e Lettere, nei cui Atti presentò le sue ultime riflessioni. Si veda [Malusa, 2004].
Citato in [Canepa-Fenaroli, 2009].
*Lettere n° 13, 30, 55.*

**Polli** Giovanni (Milano 1812-1880). Si laureò in Medicina e chirurgia a Pavia nel 1837. Insegnò Chimica in diverse scuole superiori di Milano. Fu membro effettivo del R. Istituto Lombardo dal 1854. Si veda *Enciclopedia Italiana* (A. Lustig).
*Lettere n° 11, 30, 31, 33, 34, 37, 38.*

**Poncelet** Jean Victor (Metz, Francia 1788 - Parigi 1867). Dopo aver frequentato i corsi dell'École Polytechnique a Parigi entrò alla Scuola militare di Metz e seguì Napoleone nella campagna di Russia, come ufficiale del Genio. Successivamente fu fatto prigioniero e restò in Russia fino al 1814. In questo periodo stese le basi del suo trattato di Geometria proiettiva. Al suo rimpatrio in Francia, nel 1825 diventò professore di Meccanica alla École d'application. Si interessò particolarmente alla progettazione di turbine e ruote idrauliche, e nel 1826 disegnò una turbina "Francis" che però verrà realizzata solo nel 1838. Nel suo libro *Meccanica industriale* del 1829 studiò il lavoro e l'energia cinetica. Dal 1815 al 1825 fu ingegnere militare a Metz e dal 1825 al 1835 fu, nella stessa città, professore di Meccanica. Nel 1838 si spostò a Parigi dove insegnò alla Facoltà di Scienze fino al 1848 quando divenne comandante dell'École Polytechnique con il grado di generale. Si ritirò quindi nel 1850 per dedicarsi alla ricerca matematica. Si veda DSB (R. Taton). Citato in [Cerroni-Fenaroli, 2007], [Canepa-Fenaroli, 2009], [Enea-Gatto, 2009], [Enea, 2009].
*Lettera n° 49.*

**Richelmy** Prospero (Torino 1813-1884). Si laureò in Ingegneria presso l'Università di Torino nel 1833 dove insegnò dapprima Matematica come ripetitore e poi Idraulica come titolare. Nel 1860 passò alla nuova Scuola d'Applicazione per Ingegneri in qualità di Direttore e di docente sia di Idraulica che di Meccanica applicata. Si veda G. Curioni, "Commemorazione di Prospero Richelmy", *Atti della Società degli Ingegneri e degli Industriali di Torino*, 1884, pp. 47-49.
*Lettera n° 18.*

**Rignon** Felice (Torino 1829-1914). Politico, fu più volte sindaco di Torino. Venne nominato Senatore nel 1891.
*Lettera n° 10.*

**Riva**, professore di Guido Clericetti al Collegio. Non meglio identificato.
*Lettera n° 11.*

**Robilant** (conte di), Nicolis Carlo Felice (Torino 1826 - Londra 1888). Fu Ministro degli Esteri dal 1885 al 1887 quando venne nominato ambasciatore a Londra. Si veda *Enciclopedia Italiana* (M. Menghini).
*Lettera n° 58.*

**Rodriguez** Francesco, ingegnere. Esule in Svizzera dopo i moti milanesi del 1848, insegnò Geodesia e Meccanica al Liceo cantonale di Lugano e al Ginnasio di Bellinzona fino al 1859. Rientrato in Italia nel 1860 ricoprì l'incarico di Ispettore delle Strade Ferrate Lombarde. Dal 1861 al 1871 fu Preside del R. Istituto Tecnico (poi R. Istituto Tecnico di S. Marta). Nel 1863 fece parte del Consiglio direttivo del nascente R. Istituto Tecnico Superiore insieme a F. Brioschi. Fu preside dell'Istituto Tecnico di Roma dalla sua fondazione nel 1871. Ricoprì varie cariche nella Società Geografica Italiana.
*Lettere n° 12, 15, 16, 17.*

**Rossetti** Francesco (Trento 1833 - Padova 1885). Professore di Fisica all'Università di Padova dal 1866, fu anche Preside della Facoltà di Scienze dal 1876 al 1885. Divenne Socio nazionale dei Lincei nel 1882. Si veda *Enciclopedia Italiana.*
Citato in [Giacardi-Tazzioli, 2012].
*Lettera n° 53.*

**Sacchi** Giuseppe (Milano 1804-1891). Si laureò in Giurisprudenza a Pavia e si dedicò ad attività legate soprattutto alle necessità dell'infanzia: grazie a lui a Milano si aprirono cinque asili. Nel 1860 divenne bibliotecario della Braidense e formò la Società pedagogica italiana residente in Milano. Fu membro effettivo del R. Istituto Lombardo dal 1858. Si veda *Enciclopedia Italiana* (V.B. Brunelli).
*Lettera n° 30.*

**Sangalli** Giacomo (Treviglio 1821 - Pavia 1897). Si laureò in Medicina a Pavia dove insegnò Anatomia patologica. Divenne membro effettivo del R. Istituto Lombardo dal 1868 e Senatore nel 1896.
*Lettere n° 30, 33, 41.*

**Saviotti** Carlo (Calvignano, Pavia 1845-1928). Nel 1870 si laureò in Ingegneria a Milano, ove incominciò a collaborare con Luigi Cremona nell'insegnamento della Statica grafica. Suoi sono per esempio i disegni delle figure che accompagnano nel 1872 la Memoria *Le figure reciproche nella statica grafica* che, pensata per gli studenti, diventerà un'opera classica nella formazione degli ingegneri. Nel 1873 Cremona lo portò con sé a Roma alla R. Scuola per gli Ingegneri dove ricoprì la cattedra di Statica grafica e di Meccanica applicata alle macchine fino al 1927. Un metodo per lo studio delle travature reticolari con aste caricate porta il suo nome. Il suo trattato *La statica grafica* (in tre volumi, Milano, Hoepli, 1888) munito di presentazione di Luigi Cremona, fu molto apprezzato.
*Lettere n° 8, 10.*

**Savoia** (di) Margherita, regina d'Italia. Si veda DBI (D. Adorni).
*Lettera n° 46.*

**Saracco** Giuseppe (Bistagno, Acqui 1821-1907). Avvocato e politico, ricoprì diverse cariche di governo. Divenne Senatore nel 1865. Si veda *Enciclopedia Italiana* (M. Menghini).
*Lettera n° 58.*

**Sayno** Antonio (Milano 1844-1916). Si laureò in Ingegneria civile presso il Politecnico di Milano nel 1867. L'anno successivo fu nominato assistente di Statica grafica e nel 1872 divenne professore straordinario di Scienza delle costruzioni. Dal 1876 fu incaricato dell'insegnamento della Geometria descrittiva e ordinario dal 1894. Fu consigliere e, dal 1896 al 1899, fu sindaco a Monza.
*Lettera n° 42.*

**Sbarbaro** Pietro (Savona 1838 - Roma 1893). Docente di Economia politica all'Università di Modena venne allontanato dalla cattedra per ragioni politiche. Nel 1878 insegnò all'Università di Napoli, ma venne nuovamente destituito dal Ministro Baccelli nel 1883. Si dedicò infine al giornalismo. Si veda *Enciclopedia Italiana* (M. Menghini).
*Lettera n° 47.*

**Schiaparelli** Giovanni Virgilio (Savigliano, Cuneo 1835 - Milano 1910). Laureatosi in Ingegneria a Torino nel 1854, dopo aver studiato con una borsa di perfezionamento a Berlino con Johann F. Encke e in Russia con Otto Struve, fece ritorno in Italia dove nel 1860 entrò a far parte come secondo astronomo dell'Osservatorio di Brera (Milano). Nel 1862, alla morte di Francesco Carlini, gli successe nella direzione, che mantenne fino al 1910. Osservò la superficie di Marte scoprendo delle strutture rettilinee che chiamò canali. Determinò anche i periodi di rotazione di Mercurio e Venere; scoprì l'asteroide Esperia e collegò gli sciami meteorici con le comete. Negli ultimi anni compì fondamentali studi storici sull'Astronomia greca e dell'Antico Testamento. Ottenne, fra gli altri riconoscimenti, il premio Lalande della Académie des Sciences di Parigi e una medaglia d'oro della Royal Astronomical Society di Londra. Fu il maggiore astronomo italiano del secolo e, insieme ad Angelo Secchi, fu uno dei pionieri dell'astrofisica. Si vedano DSB (G. Abetti) e [Tricomi 1962].
Citato in [Cerroni-Fenaroli, 2007], [Cerroni-Martini, 2009], [Enea-Gatto, 2009], [Giacardi-Tazzioli, 2012].
*Lettere n° 33, 37.*

**Schiff**, professore di Ginevra. Non identificato.
*Lettera n° 28.*

**Scialoja** Antonio (Napoli 1817 - Procida, Napoli 1877). Laureatosi in Giurisprudenza a Napoli nel 1841, divenne professore di Economia politica all'Università di Torino nel 1846. Nel 1848 fu Ministro dell'Agricoltura e del Commercio nel Regno delle Due Sicilie durante il governo liberale. Arrestato dopo la repressione del 1849, fu condannato all'esilio perpetuo e quindi costretto a rifugiarsi a Torino. Ritornò nuovamente a Napoli nel 1860, dopo la spedizione dei Mille, per diventare Ministro delle Finanze nel governo provvisorio di Garibaldi. In seguito fu segretario generale al Ministero dell'Agricoltura nel primo governo Ricasoli (dal 12 giugno 1861 al 3 marzo 1862) del Regno d'Italia, consigliere della Corte dei Conti e Senatore dal 1862, Ministro delle Finanze nel secondo Governo La Marmora (dal 31 dicembre 1865 al 20 giugno 1866) e poi nel secondo governo Ricasoli (dal 20 giugno 1866 al 10 aprile 1867). Fu Ministro della Pubblica Istruzione dal 5 agosto 1872 al 6 febbraio 1874. Si veda *Enciclopedia Italiana*.
Citato in [Cerroni-Fenaroli, 2007], [Canepa-Fenaroli, 2009].
*Lettere n° 8, 10, 12.*

**Sciolette** Carlo, ingegnere. Fu assistente nel R. Istituto Tecnico di Roma quando L. Cremona era Direttore.
*Lettera n° 26.*

**Sebregondi** Francesco (Como 1827 - Milano 1888). Compiuti studi classici, si dedicò all'attività letteraria. Fu segretario dell'Accademia di Belle Arti a Milano e assessore al Comune di Milano dal 1868 al 1883. Si veda G. Carotti, *Atti dell'Accademia di Belle Arti*, 1888, pp. 57-58.
*Lettera n° 6.*

**Sinigaglia** Francesco, ingegnere. Fu assistente nel R. Istituto Tecnico di Roma quando L. Cremona era Direttore. Probabilmente poi si dedicò alla libera professione, infatti nel 1881 era Ingegnere-capo presso un'industria di Ancona.
*Lettera n° 44.*

**Stoppani** Antonio (Lecco 1824 - Milano 1891). Da seminarista prese parte ai Moti milanesi del 1848, poi venne ordinato sacerdote. Si dedicò alle scienze naturali e divenne professore di Geologia nell'Università di Pavia per gli anni 1861-62. Si spostò quindi al R. Istituto Tecnico di Milano fino al 1878 dove insegnò Geognosia e Mineralogia applicata. Dal 1878 al 1883 fu sulla cattedra di Geologia e Geofisica all'Istituto di studî superiori di Firenze per poi tornare a Milano dove diresse anche il Museo civico di storia naturale. Fu membro effettivo del R. Istituto Lombardo dal 1862 e Socio nazionale dei Lincei dal 1875. Si veda *Enciclopedia Italiana*.
*Lettere n° 19, 20, 33, 34.*

**Taddeucci** Mary Isabella, insegnò la lingua inglese a L. Cremona durante la sua permanenza a Roma.
*Lettera n° 34.*

**Tagliasacchi** Gioachimo, ingegnere. Fu assessore alla Commissione edilizia a Milano.
*Lettera n° 48.*

**Tagliasacchi** Zaverio (? - Milano 1872). Ingegnere, fece parte della Società italiana di Scienze Naturali.
*Lettera n° 4.*

**Tamanini** Francesco Saverio (Trento 1833-1886). Nel 1856 si laureò in Matematica presso l'Imperial Regia Università di Padova frequentando l'Accademia di Belle Arti di Venezia e quella di Milano. Tornato a Trento esercitò la professione di architetto-ingegnere realizzando molteplici costruzioni e ristrutturazioni. È considerato uno dei più importanti professionisti impegnati nello sviluppo di Trento nella seconda metà dell'Ottocento. Si veda [Burnazzi-Campolongo, 2011, pp. 40-44].
*Lettere n° 22, 27, 28.*

**Taramelli** Torquato (Bergamo 1845 - Pavia 1922). Si laureò in Geologia a Palermo e divenne assistente di Stoppani a Milano fino al 1866. Insegnò Scienze naturali nell'Istituto Tecnico di Udine, poi Geologia nell'Università di Pavia, dal 1876 al 1920. Nel 1880 gli fu assegnato il primo premio per la Geologia e per la Mineralogia dalla Reale Accademia dei Lincei, per i suoi studî sulla geologia delle provincie venete. Fu membro effettivo del R. Istituto Lombardo dal 1880 e Socio nazionale dei Lincei dal 1891. Si veda *Enciclopedia Italiana* (M. Piazza).
*Lettere n° 30, 33.*

**Tardy** Placido (Messina 1816 - Firenze 1914). Compì gli studi tra Milano e Parigi, infine si laureò in Matematica nella R. Università di Messina nel 1841. Qui insegnò Analisi infinitesimale fino al 1847, quando si rifugiò a Firenze. Dal 1851 al 1881 fu a Genova, dapprima presso la Scuola di Marina, poi, dal 1859, presso l'Università della quale fu anche Rettore. Dal 1893 fu Socio nazionale dei Lincei. Si vedano [Cerroni-Fenaroli, 2007], [Cerroni-Martini, 2009], [Canepa-Fenaroli, 2009], [Giacardi-Tazzioli, 2012] e *Enciclopedia Italiana.*
*Lettere n° 4, 53.*

**Tassi** (? - 1874), ingegnere. Non identificato.
*Lettera n° 15.*

**Tatti** Luigi (Como 1808 - Montano Lucino, Como 1881). Si laureò in Matematica a Pavia nel 1829, poi si interessò allo studio dell'Architettura a Milano e a Roma. Progettò diverse strutture e lavorò per le Pubbliche Costruzioni del comune di Milano. Appartenne al Collegio degli Ingegneri e Architetti di Milano e fu membro effettivo del R. Istituto Lombardo dal 1860. Si veda A. Pestalozza, "Cenni biografici sull'Ingegnere Luigi Tatti", *Il Politecnico - Giornale dell'ingegnere architetto civile ed industriale*, v. 14, 1882, pp. 96-104.
*Lettere n° 10, 11, 12, 22.*

**Tecchio** Sebastiano (Vicenza 1807 - Venezia 1886). Avvocato, ricoprì diverse cariche politiche. Venne nominato Senatore nel 1866. Si veda Si veda *Enciclopedia Italiana* (M. Menghini).
*Lettera n° 22.*

**Tenerelli** Francesco (Leonforte, Enna 1839 - Catania 1899). Politico, fu sindaco di Catania tra il 1875 e il 1877. Venne nominato Senatore nel 1886.
*Lettera n° 41.*

**Tresca** Henri Édouard (Dunkerque, Francia 1814 - Parigi 1885). Ingegnere di ponti e strade. A Parigi insegnò Meccanica industriale nel Conservatoire des arts et métiers e Meccanica applicata nel 1859 all'École centrale des arts et manufactures. Nel 1872 entrò all'Accademia delle Scienze a Parigi. Si veda *Enciclopedia Italiana* (G. Albenga).
*Lettera n° 49.*

**Vannucci** Atto (Tobbiana, Pistoia 1810 - Firenze 1883). Fu storico e filologo. Insegnò dapprima latino a Prato, in seguito prese parte ai moti toscani; diresse due riviste a Firenze. Divenne docente di Storia a Lugano, nel 1859 gli fu affidata la direzione della Biblioteca Magliabechiana a Firenze, dove passò a insegnare Letteratura latina nel locale Istituto di studî superiori. Partecipò alla Pia Fondazione Edoardo Kramer di Milano. Fu nominato Senatore nel 1865. Si veda Si veda *Enciclopedia Italiana* (G. Mazzoni).
*Lettera n° 50.*

**Venino**. Non identificata.
*Lettera n° 50.*

**Verga** Andrea (Treviglio 1811 - Milano 1895). Si laureò in Medicina all'Università di Pavia dove divenne assistente alla cattedra di Anatomia. Nel 1842 si trasferì a Milano e cominciò a dedicarsi alla psichiatria, materia della quale è considerato uno dei padri fondatori in Italia. Dal 1852 al 1865 fu direttore dell'Ospedale Maggiore. Fu attivo anche all'interno del consiglio comunale della città, divenne membro effettivo del R. Istituto Lombardo nel 1848 di cui fu presidente negli anni 1857-1858 e 1864-1865. . Si veda Si veda *Enciclopedia Italiana* (A. Palmerini).
*Lettere n° 33, 34.*

**Vimercati** Gaetano, assessore a varie commissioni amministratrici del prestito civico nel Comune di Milano.
*Lettera n° 48.*

**Zanetti**, amici di famiglia dei Clericetti. Non identificati.
*Lettera n° 51.*

# Referenze bibliografiche

[Burnazzi-Campolongo, 2011] = E. Burnazzi, F. Campolongo, *Palazzo Ranzi a Trento, un cantiere alle soglie d'Italia*, Provincia Autonoma di Trento, 2011

E. Canadelli, P. Zocchi (a cura di), *Milano scientifica 1875-1924*, Sironi Editore, Milano, 2008

[Canepa-Fenaroli, 2009] = *Il carteggio Bellavitis-Tardy (1852-1880)*, a cura di G. Canepa e G. Fenaroli (eds), Mimesis 2009, Milano

[Cerroni-Fenaroli, 2007] = *Il carteggio Cremona-Tardy (1860-1886)*, a cura di C. Cerroni e G. Fenaroli (eds), Mimesis 2007, Milano

[Cerroni-Martini, 2009] = *Il carteggio Betti-Tardy (1850-1891)*, a cura di C. Cerroni e L. Martini (eds), Mimesis 2009, Milano

DBI = *Dizionario Biografico degli Italiani*, Istituto dell'Enciclopedia Italiana, 1960- , Roma

DSB = *Dictionary of Scientific Biography*, C.C. Gillespie (ed.), C. Scribner's Sons 1970-1976, New York

DSS = *Dizionario Storico della Svizzera* (www.hls-dhs-dss.ch/)

*Enciclopedia Italiana* = *Enciclopedia Italiana Treccani*, http://www.treccani.it

[Enea, 2009] = *Il carteggio Beltrami-Chelini*, a cura di M.R. Enea (ed.), Mimesis 2009, Milano

[Enea-Gatto, 2009] = *Le carte di Domenico Chelini dell'Archivio Generale delle Scuole Pie e la corrispondenza Chelini-Cremona*, a cura di M. R. Enea e R. Gatto (eds), Mimesis 2009, Milano

[Giacardi-Tazzioli, 2012] = *Le lettere di Eugenio Beltrami a Betti, Tardy e Gherardi*, a cura di L. Giacardi e R. Tazzioli (eds), Mimesis 2012, Milano

[Malusa, 2004] = L. Malusa, "La storiografia filosofica in Italia nella prima metà dell'Ottocento" in: *Storia delle storie generali della filosofia*, G. Santinello, G. Piaia (a cura di), Editrice Antenore 2004, Roma-Padova

ONDB = *Oxford Dictionary of National Biography*, Oxford University Press, 2004-16

[Tricomi 1962] = F.G. Tricomi, "Matematici italiani del primo secolo dello stato unitario", *Memorie dell'Accademia delle Scienze di Torino, classe di Scienze Fisiche matematiche e Naturali*, s. 4, n. 1 (1962), pp. 1-120

## Iscrizioni a Celeste Clericetti

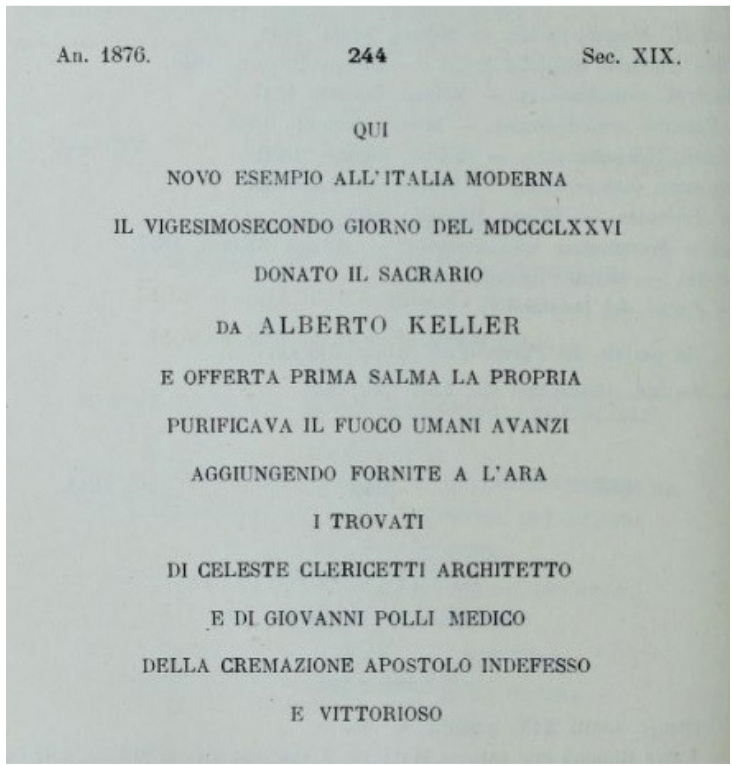

Iscrizione posta nel Tempio crematorio del Cimitero monumentale di Milano in occasione dell'inaugurazione del forno realizzato da C. Clericetti e G. Polli.

V. Forcella (a cura di), *Iscrizioni delle chiese ed altri edifici di Milano*, v. 7, 1891

Iscrizione funebre (Cimitero monumentale di Milano, Giardini di levante fino al 1917, ora Galleria superiore di ponente).

*Ibidem*

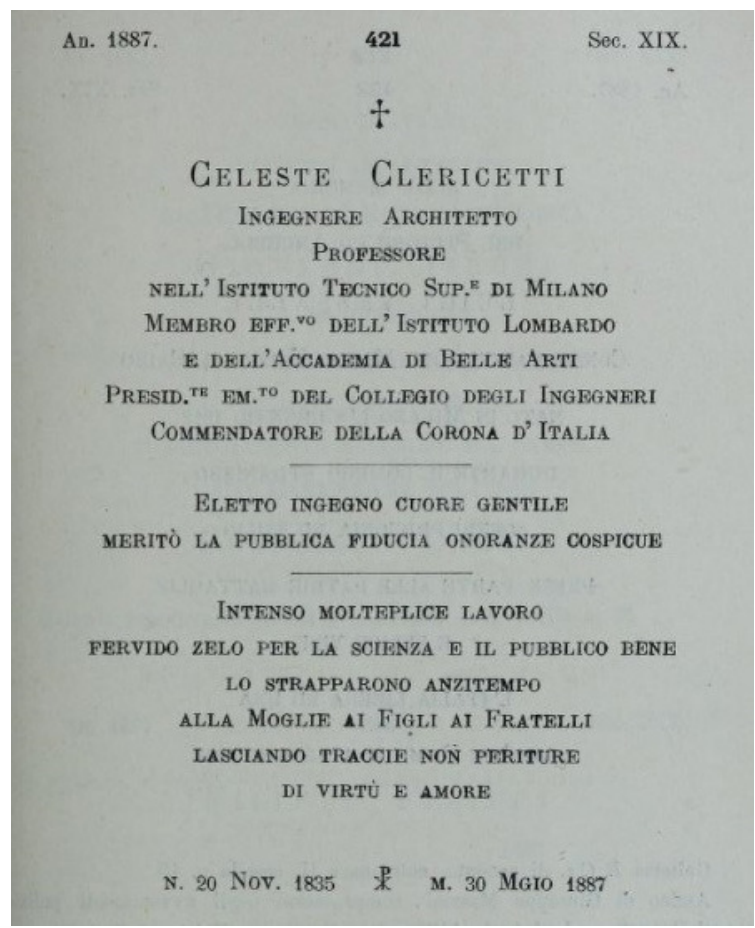

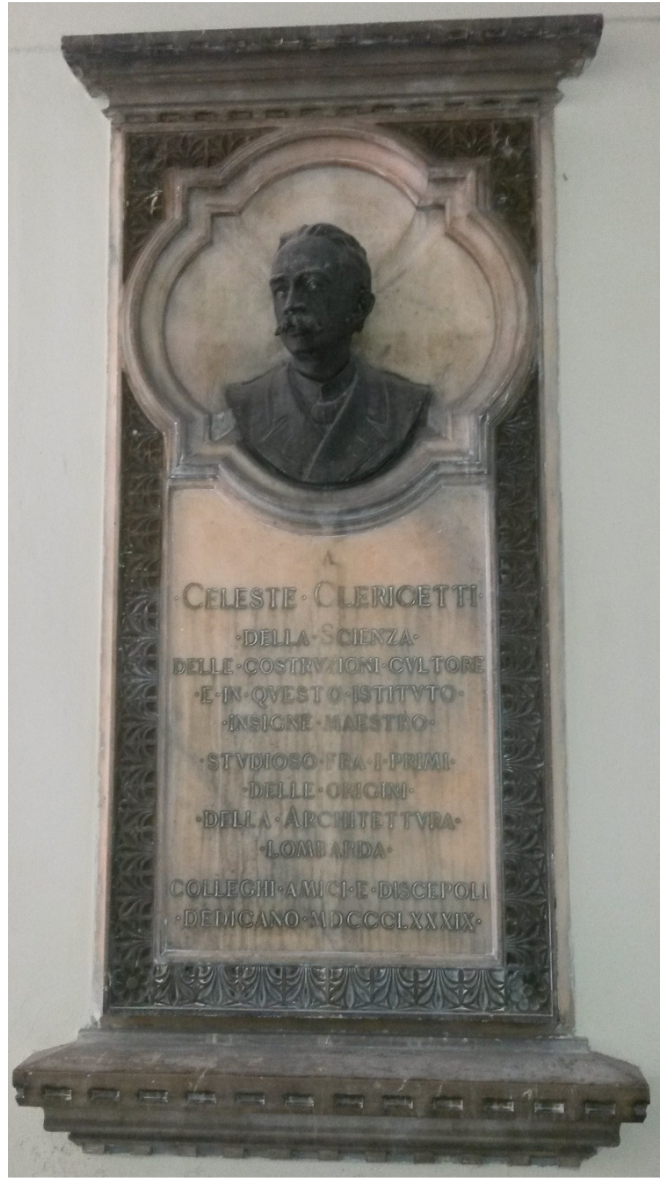

```
An. 1889.        336        Sec. XIX.

              A
      CELESTE · CLERICETTI
         DELLA · SCIENZA ·
   DELLE · COSTRUZIONI · CULTORE ·
      · E · IN · QUESTO · ISTITUTO ·
         INSIGNE · MAESTRO ·
      · STUDIOSO · FRA · I · PRIMI ·
             · DELLE · ORIGINI ·
         · DELLA · ARCHITETTURA ·
                · LOMBARDA ·
    · COLLEGHI · AMICI · E · DISCEPOLI ·
           · DEDICANO · MDCCCLXXXIX ·
```

Busto posto in memoria di C. Clericetti a due anni dalla morte (Politecnico di Milano, ora all'interno della Facoltà di Architettura).

A sinistra, trascrizione della dedica. *Ibidem*, v. 9, 1892